\documentclass[final,leqno,onefignum,onetabnum]{siamltexmm}
\usepackage{showkeys}

\usepackage{amssymb}
\usepackage{bm}
\usepackage{amsmath}

\newcommand{\be}{\begin{equation}}
\newcommand{\ee}{\end{equation}}

\newcommand{\Ah}{\widehat{A}}
\newcommand{\qt}{\widetilde{q}}

\newcommand{\gt}{\widetilde{g}}

\newcommand{\bD}{{\bf D}}
\newcommand{\bE}{{\bf E}}

\newcommand{\bI}{{\bf I}}
\newcommand{\bL}{{\bf L}}
\newcommand{\bM}{{\bf M}}

\newcommand{\bQ}{{\bf Q}}
\newcommand{\bR}{{\bf R}}
\newcommand{\bS}{{\bf S}}
\newcommand{\bU}{{\bf U}}
\newcommand{\bV}{{\bf V}}

\newcommand{\bY}{{\bf Y}}
\newcommand{\bZ}{{\bf Z}}

\newcommand{\bBt}{\widetilde{\bf B}}
\newcommand{\bFt}{\widetilde{\bf F}}
\newcommand{\bPt}{\widetilde{\bf P}}

\newcommand{\bb}{{\bf b}}

\newcommand{\bfe}{{\bf e}}

\newcommand{\bn}{{\bf n}}

\newcommand{\bu}{{\bf u}}
\newcommand{\bv}{{\bf v}}
\newcommand{\bw}{{\bf w}}

\newcommand{\by}{{\bf y}}
\newcommand{\bz}{{\bf z}}

\newcommand{\buh}{\widehat{\bf u}}

\newcommand{\bzero}{\boldsymbol{0}}

\newcommand{\balpha}{\boldsymbol{\alpha}}
\newcommand{\bbeta}{\boldsymbol{\beta}}

\newcommand{\cA}{{\mathcal A}}
\newcommand{\cB}{{\mathcal B}}
\newcommand{\cC}{{\mathcal C}}

\newcommand{\cI}{{\mathcal I}}
\newcommand{\cK}{{\mathcal K}}
\newcommand{\cM}{{\mathcal M}}

\newtheorem{algorithm}[theorem]{Algorithm}

\graphicspath{{./FIGS/}}

\title{A nonlinear method for imaging with acoustic waves via reduced order model backprojection} 

\author{Vladimir Druskin\footnote{Schlumberger-Doll Research Center, 
1 Hampshire St., Cambridge, MA 02139-1578 USA (\email{druskin1@slb.com})},
Alexander V. Mamonov\footnote{University of Houston, Department of Mathematics,
3551 Cullen Blvd., Houston, TX 77204-3008 USA (\email{mamonov@math.uh.edu})} and
Mikhail Zaslavsky\footnote{Schlumberger-Doll Research Center, 
1 Hampshire St., Cambridge, MA 02139-1578 USA (\mbox{\email{mzaslavsky@slb.com}})}}

\begin{document}
\maketitle
\newcommand{\slugmaster}{}

\begin{abstract}
We introduce a novel nonlinear imaging method for the acoustic wave equation based on data-driven 
model order reduction. The objective is to image the discontinuities of the acoustic velocity, a coefficient 
of the scalar wave equation from the discretely sampled time domain data measured at an array of 
transducers that can act as both sources and receivers. We treat the wave equation along with 
transducer functionals as a dynamical system. A reduced order model (ROM) for the propagator of 
such system can be computed so that it interpolates exactly the measured time domain data. 
The resulting ROM is an orthogonal projection of the propagator on the subspace of the snapshots of 
solutions of the acoustic wave equation. While the wavefield snapshots are unknown, the projection 
ROM can be computed entirely from the measured data, thus we refer to such ROM as data-driven. 
The image is obtained by backprojecting the ROM. Since the basis functions for the projection subspace 
are not known, we replace them with the ones computed for a known smooth kinematic velocity model. 
A crucial step of ROM construction is an implicit orthogonalization of solution snapshots. It is a nonlinear 
procedure that differentiates our approach from the conventional linear imaging methods (Kirchhoff 
migration and reverse time migration - RTM). It resolves all dynamical behavior captured by the 
data, so the error from the imperfect knowledge of the velocity model is purely kinematic. This allows 
for almost complete removal of multiple reflection artifacts, while simultaneously improving the 
resolution in the range direction compared to conventional RTM.
\end{abstract}



\pagestyle{myheadings}
\thispagestyle{plain}
\markboth{V. Druskin, A.V. Mamonov, M. Zaslavsky}{Acoustic nonlinear imaging with ROMs}

\section{Introduction}

We consider here a problem of estimating a coefficient of a scalar wave equation, an acoustic velocity, 
in some domain of interest, from the measurements of its solutions, typically referred to as wavefields, 
on an array of transducers located on the boundary of said domain. Two aspects of this problem are 
typically considered. First, a quantitative reconstruction of the coefficient is the subject of the so-called 
full waveform inversion (FWI) \cite{crase1990robust, bunks1995multiscale, virieux2009overview}. 
Second, one may seek a qualitative estimate of the acoustic velocity, say its discontinuities. 
The discontinuities in the acoustic velocity lead to reflections of waves and thus may be referred to as 
reflectors or scatterers. The problem of estimating the reflectors (scatterers) from the measurements 
of the reflected (scattered) wavefields is often referred to as the imaging problem 
\cite{hagedoorn1954process, claerbout1971toward, edler2004use}.
It is the imaging of reflectors with acoustic waves that is of interest to us here. 

The dependency of wavefields on the acoustic velocity is highly nonlinear, which makes both the 
quantitative and qualitative estimation of the velocity very challenging.  Most state-of-the-art imaging
methods employed in seismic geophysical exploration and medical ultrasound diagnostics are linear
in the recorded data. When more than one reflector is present, first reflections (primaries) from them 
are typically mixed together with multiply-reflected wavefields (multiples). This creates additional events 
in the data compared to linearized case with primaries only. Linear imaging methods may interpret
these events as extra reflectors that are not present in the actual medium. This creates image artifacts 
known as multiples in geophysics \cite{verschuur1992adaptive, weglein1997inverse} or reverberation 
(``comet tail'') artifacts in medical ultrasound \cite{kremkau1986artifacts, ziskin1982comet}. 
Consequently, complex (and often computationally expensive) data pre-processing techniques are applied 
to filter the multiple reflection events out of the data before feeding it to imaging algorithms. 

Many conventional methods for imaging with waves rely on the time-reversibility of the wave equation. 
In one form or another they back-propagate the data from the transducers using certain assumptions 
on the kinematics of the unknown medium. Such approaches are commonly known as migration 
\cite{schneider1978integral, docherty1991brief, baysal1983reverse, chang1987elastic} or time-reversal 
\cite{derode1995robust, borcea2002imaging} methods. Here we propose a completely different 
approach. Instead of propagating the data, we use the theory of model order reduction 
\cite{antoulas2001survey,antoulas2010interpolatory,grimme1997krylov} to construct an approximate 
wave equation propagator. The approximant is derived to satisfy certain data interpolation conditions,
as explained later. If the data is discretely sampled in time, 
then the reduced order model (ROM) interpolating the data is an orthogonal projection of the propagator 
on the (Krylov) subspace spanned by the wavefield snapshots taken precisely at the data sampling instants.
Computing the projection requires an orthogonalization of wavefiled snapshots to be performed.  Since 
the wavefields are not known inside the domain of interest, the orthogonalization is performed implicitly.
It is a nonlinear operation that is crucial to our approach, as it allows to suppress the effects of multiple
reflections and probe the medium of interest with localized wavefields. 

The use of ROMs in inversion and imaging is not an entirely new idea. Resistor network ROMs were used
successfully to obtain quantitative reconstructions of conductivity in two-dimensional electrical impedance 
tomography (EIT) both directly and iteratively
\cite{borcea2008electrical,borcea2010circular,borcea2010pyramidal,borcea2011study,borcea2011resistor}. 
They were also employed in \cite{borcea2016discrete} to reconstruct Schr\"{o}dinger potentials.
While the resistor networks interpolate the EIT or Schr\"{o}dinger data, they are not ROMs of projection 
type, unlike the ones considered here. Projection-based ROMs were first studied in the context of inversion for 
the coefficient of a parabolic equation in \cite{druskin2013solution,borcea2014model}. Finally, in 
\cite{druskin2015direct} reduced models interpolating the time domain wave equation data were 
introduced to obtain quantitative reconstructions of acoustic velocity. Here we consider ROMs of the 
same type as in \cite{druskin2015direct}, and their construction follows a similar procedure. However, 
the use of ROMs is different. While in \cite{druskin2015direct} the information about the acoustic velocity 
is extracted directly from the matrix entries of the ROM, here we perform a backprojection operation on 
the ROM to go from the reduced model space back to the ``physical'' space of the domain of interest.  
This method was first introduced by the authors in a conference abstract \cite{mamonov2015rombackproj}, 
and we present here a thorough exposition of the approach briefly outlined therein.

The paper is organized as follows. We begin in section \ref{sec:formulation} with a forward data model 
and a formulation of inversion and imaging problems. To simplify the subsequent ROM construction,
we symmetrize the data model in section \ref{sec:sym}. Then we define the time-domain interpolatory 
ROM and discuss some of its properties in section \ref{sec:rom}.  An algorithm for computing the ROM
from the data is given in section \ref{sec:romcompute}. After the ROM is computed, we introduce its
backprojection and the corresponding imaging functional. Both questions are addressed in section 
\ref{sec:backproj}.  An interpretation of the imaging functional that explicitly separates the kinematic
and reflective behavior is given in section \ref{sec:interp}. A thorough numerical study is performed in 
section \ref{sec:num}. In particular, the orthogonalized wavefield snapshots that span the projection 
subspace are discussed in section \ref{subsec:orthsnap}. Approximations of delta functions by the 
outer products of orthogonalized snapshots are considered in section \ref{subsec:approxdelta}. 
We compare the proposed ROM backprojection imaging method to conventional reverse time 
migration (RTM) in section \ref{subsec:imgcomprtm}. Stability of ROM computation in the presence 
of noise and a regularization algorithm are addressed in section \ref{subsec:noisereg}. The numerical 
study is concluded in section \ref{subsec:largescale} with two large scale synthetic examples from 
geophysical exploration and medical ultrasound diagnostics. Finally, we conclude in section 
\ref{sec:conclude} and discuss the three main directions of future research.

\section{Problem formulation}
\label{sec:formulation}

Consider a scalar wave equation for the acoustic pressure $p(x,t)$
\be
p_{tt} = c(x)^2 \Delta p + q_t(t) \phi(x) e(x), \quad -\infty < t < \infty, \quad x \in \Omega,
\label{eqn:wavep}
\ee
where $\Omega$ is the domain of interest in $\mathbb{R}^d$, $d = 2,3$ . 
Domain $\Omega$ can be either infinite or finite. Here we assume $\Omega$ is finite with boundary 
$\cB = \partial \Omega$ split into the \emph{accessible} $\cB_A$ and \emph{inaccessible}  
$\cB_I = \cB \setminus \cB_A$ parts respectively. We also assume that the \emph{acoustic velocity} 
$c(x)$ is piecewise smooth in $\Omega$.

Wave equation (\ref{eqn:wavep}) is driven by the source term comprised of a wavelet $q(t)$,
source weighting function $\phi(x)$ and a distribution $e(x)$. Source distribution $e(x)$ is 
supported on the accessible boundary
\be
\mbox{supp}\; e \subset \cB_A,
\label{eqn:suppBA}
\ee
and is defined so that
\be
\int_{\mathbb{R}^d} f(x) e(x) dx = \int_{\Omega} f(x) e(x) dx = 
\int_{\mbox{supp}\; e} f(y) d\Sigma_y,
\label{eqn:distribe}
\ee
for any smooth test function $f$, where $d\Sigma_y$ is the surface element 
(arc length in the case $d=2$) at point $y \in \cB$.

Since (\ref{eqn:wavep}) is defined for all $t \in (-\infty, \infty)$, we enforce a decay
condition for large negative times
\be
\lim_{t \to -\infty} p(x, t) = 0.
\ee

We prescribe reflective boundary conditions on the accessible boundary
\be
\bn(x) \cdot \nabla p(x, t) = 0, \quad x \in \cB_A,
\label{eqn:bdrypA}
\ee
where $\bn(x)$ is the unit outward normal at $x \in \cB_A$. We also prescribe zero Dirichlet
conditions on the inaccessible boundary
\be
p(x, t) = 0, \quad x \in \cB_I,
\label{eqn:bdrypI}
\ee
although zero Neumann conditions can be used instead.

Following \cite{druskin2015direct} we assume that the source wavelet $q$ is an even, sufficiently
smooth approximation of $\delta(t)$ with a non-negative Fourier transform
\be
\qt(s^2) = \int_0^\infty 2 \cos(ts) q(t) dt.
\label{eqn:qfourier}
\ee
To simplify the exposition we assume that different sources correspond to varying $\phi(x)$ and $e(x)$, 
while the wavelet $q(t)$ remains the same, i.e. each source emits the same waveform.

All possible sources and receivers on $\cB_A$ provide the knowledge of the 
\emph{boundary measurement operator}
\be
\cM_p(\phi e) = \left.p(x, t)\right|_{x \in \cB_A}, \quad -\infty < t < \infty,
\label{eqn:measp}
\ee
for all $\phi$ and $e$ satisfying (\ref{eqn:suppBA}).

To simplify the exposition it is beneficial to move the source term to the initial condition using a form
of Duhamel's principle, as explained in \cite{druskin2015direct}. Using (\ref{eqn:qfourier}) and the 
fact that $q(t)$ is even, we can work instead of the pressure $p$ with its even part
\be
u(x, t) = p(x,t) + p(x, -t), \quad t>0.
\ee
It satisfies a homogeneous wave equation
\be
u_{tt} = c(x)^2 \Delta u, \quad t>0, \quad x \in \Omega,
\label{eqn:waveu}
\ee
with initial conditions
\be
u(x, 0) = [\qt(-A) \phi e] (x), \quad u_t(x, 0) = 0, \quad x \in \Omega,
\label{eqn:initu}
\ee
where we denote by $A$ the spatial part of the PDE operator in (\ref{eqn:waveu}):
\be
A = c(x)^2 \Delta.
\label{eqn:continuumA}
\ee

The homogeneous boundary conditions (\ref{eqn:bdrypA})--(\ref{eqn:bdrypI}) remain unchanged
\be
\left. \left[ \bn(x) \cdot \nabla p(x, t) \right] \right|_{x \in \cB_A}= 0, \quad 
\left. p(x, t) \right|_{x \in \cB_I} = 0.
\label{eqn:bdryu}
\ee

The knowledge of boundary measurement operator $\cM_p$ implies the knowledge of another
operator
\be
\cM_u(\phi e) = \left. u(x, t) \right|_{x \in \cB_A} = [\cM_p(\phi e)](x,t) + [\cM_p(\phi e)](x,-t), 
\quad t>0, \quad x \in \cB_A.
\label{eqn:measu}
\ee
In practice, we do not have the full knowledge of the boundary measurement operator $\cM_u$. 
Instead, we take a finite sampling of it in both space and time. 

First, consider the spatial measurement sampling. We place $m$ transducers on the accessible 
boundary that act simultaneously as sources and receivers. The corresponding source-receiver
distributions are $e_j(x)$, $j=1,2,\ldots,m,$ with disjoint supports
\be
\mbox{supp}\; e_j \cap \mbox{supp}\; e_i = \varnothing, \quad j \neq i.
\ee
Since the supports are disjoint, we fix single source and receiver weighting functions $\phi(x)$ and 
$\psi(x)$ respectively. Their particular choice is discussed in section \ref{sec:sym}.
Then the \emph{continuum time data} that we can in principle measure with such transducer
set-up is an $m \times m$ matrix-valued function of time
\be
\bD_{i,j}(t) = \int_{\mbox{supp}\;e_i} \psi(y)  [\cM_u(\phi e_j)](y, t) d\Sigma_y, \quad t>0,
\ee
which can also be written using (\ref{eqn:distribe}) and (\ref{eqn:measu}) as
\be
\bD_{i,j}(t) = \int_{\Omega} \psi(x) e_i(x)  u(x, t) dx, \quad t>0,
\label{eqn:datat}
\ee
where $u(x,t)$ satisfies  initial conditions (\ref{eqn:initu}) with $e = e_j(x)$.
Hereafter we denote matrices and matrix-valued functions by bold letters. Physically, the entry
$\bD_{i,j}(t)$ is the signal recorded at receiver $i$ at time $t$ while the medium is excited with 
source $j$.

Second, in practice the temporal measurements are also sampled discretely at a certain rate. 
As explained later, our method relies on the temporal sampling being uniform. Hence, we introduce
time instants $t_k = \tau k$, $k = 0, 1, 2, \ldots, 2n-1$, where $\tau$ denotes the sampling interval 
and the total number of temporal samples is $2n$.

Once the measurements are sampled in both space and time, we can define the 
\emph{sampled data} as
\be
\bD_{i,j}^k = \bD_{i,j}(t_k), \quad i,j = 1,\ldots,m, \quad k = 0,\ldots,2n-1.
\label{eqn:datas}
\ee
We are now in position to state the following two problems. 

\noindent \textbf{Inversion problem.} Given measured data (\ref{eqn:datas}) estimate the
acoustic velocity $c(x)$ in $\Omega$.

\noindent \textbf{Imaging problem.} Given measured data (\ref{eqn:datas}) and a smooth
\emph{kinematic velocity model} $c_o(x)$ estimate the (location of) discontinuities of the acoustic 
velocity $c(x)$.

In \cite{druskin2015direct} a method to solve the inversion problem in one dimension and its
generalization to higher dimensions was considered. The method relies on the techniques of model 
order reduction. Here we use similar techniques to address the imaging problem. Unlike the 
quantitative inversion problem, where no a priori knowledge of the acoustic velocity is available,
the imaging can be viewed as a problem of qualitative estimation of the velocity discontinuities on 
top of a known kinematic background. By referring to $c_o(x)$ as the kinematic model or background, 
we assume that it captures sufficiently well the kinematics of the wave equation (\ref{eqn:waveu}),
i.e. the wave travel times between pairs of points in $\Omega$ are approximately equal for $c(x)$ 
and $c_o(x)$.

\section{Symmetrized data model}
\label{sec:sym}

Our imaging approach relies on the construction of a reduced order model of the sampled data 
(\ref{eqn:datas}). In order to understand the form which such a model has to take, it is beneficial to 
transform (\ref{eqn:datat}) and (\ref{eqn:datas}) further. 

Denote by $u_j(x,t)$ the solution of (\ref{eqn:waveu}) with initial conditions
\be
u_j (x, 0) = [\qt(-A) \phi e_j] (x), \quad \frac{\partial u_j}{\partial t}(x, 0) = 0, 
\quad x \in \Omega, \quad j=1,\ldots,m,
\label{eqn:inituj}
\ee
i.e. the wavefield induced by source $j$. In what follows it is convenient to work simultaneously with
the wavefields induced by all sources. Thus, we define a row-vector-valued function of time
$\bu: \Omega \times (0, \infty) \to \mathbb{R}^{1 \times m}$ by
\be
\bu(x,t) = \left[ u_1(x,t), \ldots, u_m(x,t) \right],
\ee
which solves
\be
\bu_{tt} = A \bu, \quad t>0, \quad x \in \Omega,
\label{eqn:wavebu}
\ee
where the action of operators on row-vector-valued functions is understood component-wise.

The initial conditions for (\ref{eqn:wavebu}) are
\be
\bu (x, 0) = [\qt(-A) \phi \bfe] (x), \quad \bu_t(x, 0) = \bzero,  \quad x \in \Omega,
\label{eqn:initbu}
\ee
with
\be
\bfe = [e_1(x), \ldots, e_m(x)].
\ee

Using operator functions we can write the solution of (\ref{eqn:wavebu}) with initial conditions
(\ref{eqn:initbu}) as
\be
\bu(\,\cdot\,, t) = \cos \left( t \sqrt{-A} \right) \qt(-A) \phi \bfe.
\label{eqn:propbu}
\ee

In the derivations below it is convenient to use the following notation for products of row-vector-valued 
functions. Suppose $\bv, \bw: \Omega \times (0, \infty) \to \mathbb{R}^{1 \times m}$ then
by $\bv^* \bw: (0, \infty) \to \mathbb{R}^{m \times m}$ we denote
\be
[\bv^* \bw](t) =
\begin{bmatrix}
\left< v_1, w_1 \right>(t) & \left< v_1, w_2 \right>(t) & \ldots & \left< v_1, w_m \right>(t) \\
\left< v_2, w_1 \right>(t) & \left< v_2, w_2 \right>(t) & \ldots & \left< v_2, w_m \right>(t) \\
\vdots & \vdots & \ddots & \vdots \\ 
\left< v_m, w_1 \right>(t) & \left< v_m, w_2 \right>(t) & \ldots & \left< v_m, w_m \right>(t)
\end{bmatrix},
\label{eqn:matprod}
\ee
where
\be
\left< f, g \right>(t) = \int_{\Omega} \overline{f(x, t)} g(x, t) dx.
\label{eqn:l2innerprod}
\ee
Then using (\ref{eqn:matprod}) we can write (\ref{eqn:datat}) in matrix form as
\be
\bD(t) = [(\bfe \psi)^* \bu](t) = \bfe^* \psi \cos \left( t \sqrt{-A} \right) \qt(-A) \phi \bfe.
\label{eqn:datam}
\ee

We notice that in (\ref{eqn:datam}) there is a certain asymmetry between the source and receiver 
terms. Since the sources and receivers are collocated, it is possible to transform (\ref{eqn:datam})
to a symmetric form which makes the subsequent derivation of the reduced model much simpler. 
To achieve this we introduce the symmetrized PDE operator
\be
\Ah = c(x) \Delta c(x),
\ee
a similarity transform of $A$, since $A = c(x) \Ah c(x)^{-1}$. Using the facts that operator functions 
of $A$ commute and that the Fourier transform $\qt$ of the source wavelet is non-negative 
(see section \ref{sec:formulation}) we can rewrite (\ref{eqn:datam}) as
\be
\bD(t) = \bfe^* \psi \qt^{1/2} (- c \Ah c^{-1})
\cos \left( t \sqrt{-c \Ah c^{-1}} \right)
\qt^{1/2} (- c \Ah c^{-1}) \phi \bfe.
\ee
To transform further we observe that similarity transforms commute with analytic operator functions, 
thus
\be
\bD(t) = \bfe^* \psi c \,\qt^{1/2} (-\Ah)
\cos \left( t \sqrt{-\Ah} \right)
\qt^{1/2} (-\Ah) c^{-1} \phi \bfe.
\label{eqn:dataf2}
\ee

To write (\ref{eqn:dataf2}) in a fully symmetric form we need to make the following assumption on the 
source and receiver functions, namely
\be
\psi(x) = c^{-1}(x) \theta(x), \quad \phi(x) = c(x) \theta(x), 
\label{eqn:theta}
\ee
for some known function $\theta$. Since the source and receiver functions are always multiplied by
distributions $e_j(x)$ supported on the accessible boundary, assumption above only requires the 
knowledge of the acoustic velocity $c(x)$ on $\cB_A$ which is not restrictive. Indeed, if one has access 
to $\cB_A$ to place the transducers, it should be possible to determine $c(x)$ there as well.

Under assumption (\ref{eqn:theta}) we obtain the expression for the continuum time data in a fully 
symmetric form
\be
\bD(t) = \bb^* \cos \left( t \sqrt{-\Ah} \right) \bb, \quad t>0,
\label{eqn:datafsym}
\ee
where
\be
\bb(x) = [\qt^{1/2}(-\Ah) \theta \bfe](x),
\label{eqn:bb}
\ee
is the row-vector-valued \emph{transducer function} $\bb: \Omega \to \mathbb{R}^{1 \times m}$.

Note that $\Ah$ is self-adjoint with respect to the standard $L^2$ inner product (\ref{eqn:l2innerprod}),
and so are its operator functions including the one in (\ref{eqn:datafsym}). Since the transducer
function $\bb$ enters (\ref{eqn:datafsym}) in a symmetric fashion, the continuum time data is also 
self-adjoint
\be
[\bD(t)]^* = \bD(t), \quad t>0,
\ee
which is known as reciprocity, i.e. the data measured at receiver $i$ induced by source $j$ is the 
same as data measured at receiver $j$ induced by source $i$.

\section{Time-domain interpolatory reduced order model}
\label{sec:rom}

At the core of our imaging approach is the construction of a reduced order model for the sampled 
data $\bD^k$, $k=0,\ldots,2n-1$. To that end we need to transform the symmetrized data 
(\ref{eqn:datafsym}) using the fact that the data is sampled uniformly in time. Specifically, since
\be
\bD^k = \bD(t_k) = \bb^* \cos \left( k \tau \sqrt{-\Ah} \right) \bb = 
\bb^* \cos \left(k \arccos \left[ \cos \tau \sqrt{-\Ah} \right] \right) \bb,
\ee
we can write
\be
\bD^k = \bb^* T_k(P) \bb,
\label{eqn:cheb}
\ee
where $T_k$ are Chebyshev polynomials of the first kind of degree $k$ and $P$ is the 
\emph{propagator}
\be
P = \cos \left( \tau \sqrt{-\Ah} \right).
\label{eqn:propagator}
\ee

We seek the reduced order model in the form (\ref{eqn:cheb}), namely
\be
\bFt^k = \bBt^* T_k(\bPt) \bBt,
\label{eqn:romcheb}
\ee
where $\bPt \in \mathbb{R}^{mn \times mn}$ is the reduced order propagator and 
$\bBt \in \mathbb{R}^{mn \times m}$ is the reduced order transducer matrix.

To compute the reduced order propagator and transducer matrices we solve the following interpolation
problem
\be
\bFt^k = \bD^k, \quad k=0,\ldots,2n-1.
\label{eqn:interp}
\ee

Note that since there are $m$ transducers, reduced order matrices $\bPt$ and $\bBt$ consist of
$m \times m$ blocks. In particular, $\bPt$ consists of $n \times n$ blocks, while $\bBt$ is a block column
matrix $n$ blocks tall. This allows to use block versions of results in \cite{druskin2015direct} on 
the existence and construction of the solution to (\ref{eqn:interp}). This requires defining the following 
quantities.

First, following the procedure of section \ref{sec:sym} we introduce the row-vector-valued function
of \emph{symmetrized wavefields} $\buh : \Omega \times (0, \infty) \to \mathbb{R}^{1 \times m}$ 
as solutions of
\be
\buh_{tt} = \Ah \buh, \quad t>0, \quad x \in \Omega,
\label{eqn:wavebuh}
\ee
with initial conditions
\be
\buh(x, 0) = \bb(x), \quad \buh_t(x, 0) = \bzero.
\label{eqn:wavebuhinit}
\ee

Second, since the data $\bD^k$ is sampled uniformly in time, we also sample the symmetrized 
wavefields to obtain the \emph{symmetrized snapshots}
\be
\buh^k(x) = \buh(x, t_k) = \cos \left( k \tau \sqrt{-\Ah} \right) \bb = T_k(P) \bb.
\label{eqn:snapshotbuh}
\ee
Obviously, the sampled data $\bD^k$ in terms of the symmetrized snapshots is given by
\be
\bD^k = \bb^* \buh^k.
\label{eqn:databuh}
\ee

Third, we introduce the following block Krylov subspace
\be
\cK_n(P, \bb) = 
\mbox{span}\{ \buh^0, \buh^1,\ldots, \buh^{n-1} \} =
\mbox{span}\{ \bb, P \bb, \ldots, P^{n-1} \bb \},
\label{eqn:krylov}
\ee
where by ``block'' we refer to the fact that each $\buh^k$, $k=0,1,\ldots,n-1$ is a 
row-vector-valued function of $x$ and thus the dimension of $\cK_n(P, \bb)$ is $mn$.
We can consider an orthonormal basis for (\ref{eqn:krylov}) that we denote by 
$\bV: \Omega \to \mathbb{R}^{1 \times mn}$, which we write as
\be
\bV = [\bv^0, \bv^1, \ldots, \bv^{n-1}].
\ee
The basis $\bV$ is orthonormal in the sense of the following ``matrix'' product. Let
$\bY, \bZ: \Omega \to \mathbb{R}^{1 \times mn}$ be
\be
\bY = [\by^0, \by^1, \ldots, \by^{n-1}], \quad \bZ = [\bz^0, \bz^1, \ldots, \bz^{n-1}],
\label{eqn:bYbZ}
\ee
then
\be
\bY^* \bZ =
\begin{bmatrix}
(\by^0)^* \bz^0 & (\by^0)^* \bz^1 & \ldots & (\by^0)^* \bz^{n-1} \\
(\by^1)^* \bz^0 & (\by^1)^* \bz^1 & \ldots & (\by^1)^* \bz^{n-1} \\
\vdots & \vdots & \ddots & \vdots \\ 
(\by^{n-1})^* \bz^0 & (\by^{n-1})^* \bz^1 & \ldots & (\by^{n-1})^* \bz^{n-1}
\end{bmatrix} \in \mathbb{R}^{mn \times mn},
\label{eqn:matprodmn}
\ee
where each $(\by^k)^* \bz^l \in \mathbb{R}^{m \times m}$ is defined by (\ref{eqn:matprod})
with the distinction that $\by^k(x)$, $\bz^l(x)$ do not depend on $t$. The role of time is played by
the indices $k$ and $l$ once the temporal sampling (\ref{eqn:snapshotbuh}) is introduced. Then
we can write the orthonormality conditions for $\bV$ as
\be
\bV^*\bV = \bI_{mn}, \quad \mbox{or} \quad (\bv^i)^* \bv^j = \delta_{ij} \bI_m, \quad
i,j = 0,1,\ldots,n-1,
\label{eqn:bVorth}
\ee
where $\bI_{mn}$ and $\bI_m$ are identity matrices of sizes $mn$ and $m$ respectively.

Now we can address solving interpolation problem (\ref{eqn:interp}) to obtain the reduced model 
$\bPt$, $\bBt$. The solution is characterized by the following lemma.

\begin{lemma}
\label{lemma:proj}
Reduced order propagator $\bPt \in \mathbb{R}^{mn \times mn}$ and transducer matrix
$\bBt \in \mathbb{R}^{mn \times m}$ that solve the interpolation problem
\be
\bD^k = \bBt^* T_k(\bPt) \bBt, \quad k=0,\ldots,2n-1,
\label{eqn:lemmainterp}
\ee
are orthogonal projections 
\be
\bPt = \bV^* (P \bV), \quad \bBt = \bV^* \bb,
\label{eqn:orthproj}
\ee
of the true propagator (\ref{eqn:propagator}) and transducer function (\ref{eqn:bb}) on Krylov
subspace (\ref{eqn:krylov}), where $\bV$ is an orthonormal basis for $\cK_n(P, \bb)$ in the
sense of (\ref{eqn:bVorth}), and the products in (\ref{eqn:orthproj}) are understood in the
sense of (\ref{eqn:bYbZ})--(\ref{eqn:matprodmn}), in particular
\be
\bBt = \begin{bmatrix}
(\bv^0)^* \bb \\
(\bv^1)^* \bb \\
\vdots \\
(\bv^{n-1})^* \bb
\end{bmatrix}.
\ee
\end{lemma}

The proof of the lemma is in Appendix \ref{app:lemmaproj}. We make the following observations 
regarding Lemma \ref{lemma:proj} the implications of which are addressed in detail in the subsequent 
sections.

Not surprisingly, in order to match the data at time instants $t_k = \tau k$, $k=0,\ldots,2n-1,$ we 
project on the subspace of wavefield snapshots $\buh^k$ sampled at the same instants $t_k$. 
However, the Krylov subspace $\cK_n(P, \bb)$  is spanned by only $n$ snapshots, while the data is 
sampled at $2n$ instants. This question is addressed in the proof in Appendix \ref{app:lemmaproj}.

We observe that $\bPt$ and $\bBt$ are not unique as defined by (\ref{eqn:orthproj}). This ambiguity 
comes from the choice of an orthonormal basis in $\cK_n(P, \bb)$. If $\bPt$ and $\bBt$ are solutions of 
(\ref{eqn:lemmainterp}) then $\bQ^* \bPt \bQ$ and $\bQ^* \bBt$ are also solutions for any unitary 
$\bQ \in \mathbb{R}^{mn \times mn}$. We make use of this fact in the proof in Appendix 
\ref{app:lemmaproj}. 

However, not all orthonormal basis choices are equally useful if one wants to use the reduced model 
for imaging (or inversion, see \cite{druskin2015direct}). Inversion and imaging applications require 
the basis functions $\bv^k$ to be localized in $\Omega$ which implies a certain structure of $\bPt$ 
and $\bBt$. The exact structure of $\bPt$ and $\bBt$ is established in Lemma \ref{lemma:blocktri}
in section \ref{sec:romcompute}.

While (\ref{eqn:orthproj}) characterizes the solution of (\ref{eqn:lemmainterp}), it does not provide 
a practical way of computing $\bPt$ and $\bBt$. Indeed, it requires the knowledge of $\bV$ which 
comes from orthogonalization of symmetrized wavefield snapshots $\buh^k$, $k=0,\ldots,n-1$. 
Obviously, when solving inversion and imaging problems we do not have the knowledge of 
$\buh^k(x)$ for $x \in \Omega$. The only measurements we make are those of $\bD^k$, 
$k=0,\ldots,2n-1$. Thus, we must find a way to compute projections (\ref{eqn:orthproj}) from the 
knowledge of the sampled data $\bD^k$ only. Algorithm \ref{alg:romcompute} achieving that is derived
in section \ref{sec:romcompute}. 

\section{Reduced order model computation}
\label{sec:romcompute}

Let us introduce the following row-vector-valued function $\bU: \Omega \to \mathbb{R}^{1 \times mn}$
of symmetrized snapshots
\be
\bU = [\buh^0, \buh^1, \ldots, \buh^{n-1}],
\ee
so that $\cK_n(P, \bb) = \mbox{span}(\bU)$. 

If we had the knowledge of $\bU$, then one way to compute an orthonormal basis $\bV$ for 
$\cK_n(P, \bb)$ is to apply block Gram-Schmidt orthogonalization procedure to $\bU$. In matrix notation 
this can be written as a block QR decomposition
\be
\bU = \bV \bL^*,
\label{eqn:blockqr}
\ee
where $\bL \in \mathbb{R}^{mn \times mn}$ is a block lower triangular matrix
\be
\bL = \begin{bmatrix}
\bL_{0,0} & \bzero & \ldots & \bzero \\
\bL_{1,0} & \bL_{1,1} & \ddots & \bzero \\
\vdots & \vdots & \ddots & \vdots \\
\bL_{n-1,0} & \bL_{n-1,1} & \ldots & \bL_{n-1,n-1} \\
\end{bmatrix},
\ee
with blocks $\bL_{k,l} \in \mathbb{R}^{m \times m}$. Note that multiplication of $\bV$ by a matrix
from the right follows the standard rules of matrix multiplication.

While we do not have the knowledge of $\bU$, relation (\ref{eqn:blockqr}) implies that
\be
\bU^* \bU = \bL \bL^*,
\label{eqn:blockchol}
\ee
which is a block Cholesky decomposition of the \emph{mass matrix}
\be
\bM = \bU^* \bU \in \mathbb{R}^{mn \times mn}.
\label{eqn:massmat}
\ee
The advantage of (\ref{eqn:massmat}) over (\ref{eqn:blockqr}) is that it does not require the 
knowledge of wavefield snapshots in $\Omega$, but only the inner products between them. 

Recall from (\ref{eqn:snapshotbuh}) that the snapshots are expressed in terms of Chebyshev
polynomials as $\buh^k = T_k(P) \bb$. On the other hand, the products of Chebyshev polynomials 
satisfy
\be
T_k(x) T_l(x) = \frac{1}{2} \left( T_{k+l}(x) + T_{|k-l|}(x) \right)
\label{eqn:chebprod}
\ee
and so the inner products between the snapshots are given by
\be
(\buh^k)^* \buh^l = \bb^* T_k(P) T_l(P) \bb = 
\frac{1}{2} \left( \bb^* T_{k+l}(P) \bb + \bb^* T_{|k-l|}(P) \bb  \right).
\ee
Combining (\ref{eqn:snapshotbuh}) and (\ref{eqn:databuh}) the data is 
$\bD^k = \bb^* T_k(P) \bb$, hence the blocks of the mass matrix can be expressed just in terms of 
the sampled data
\be
\bM_{k,l} = (\bU^* \bU)_{k,l} = (\buh^k)^* \buh^l = 
\frac{1}{2} \left( \bD^{k+l} + \bD^{|k-l|} \right),
\label{eqn:massmatdata}
\ee
where $\bM_{k,l} \in \mathbb{R}^{m \times m}$, $k,l=0,\ldots,n-1$. Thus, in order to compute 
the block Cholesky factor $\bL$ of $\bM$ only the knowledge of the sampled data is required. 

Note that the mass matrix notation for the Gramian $\bU^* \bU = \bM$ that we use is not a 
coincidence. Indeed, if we want to approximate (\ref{eqn:snapshotbuh}) on the Krylov subspace 
(\ref{eqn:krylov}) with a Galerkin method, $\bM$ is the mass matrix and we can also define the 
\emph{stiffness matrix}
\be
\bS = \bU^* (P \bU) \in \mathbb{R}^{mn \times mn},
\label{eqn:stiffmat}
\ee
which we can also express in terms of the sampled data. 

Let us consider the action of the propagator $P$ on the snapshot $\buh^k$:
\be
P \buh^k = T_1(P) T_k(P) \bb = \frac{1}{2}\left( T_{k+1}(P)\bb + T_{|k-1|}(P)\bb \right) = 
\frac{1}{2}\left( \buh^{k+1} + \buh^{|k-1|} \right),
\label{eqn:actionP}
\ee
i.e. $P$ takes the average of the symmetrized wavefield propagated by time $\tau$ forward and
backward. Now we combine (\ref{eqn:actionP}) with (\ref{eqn:snapshotbuh}) and 
(\ref{eqn:chebprod}) to express the blocks of the stiffness matrix in terms of the sampled data
\be
\quad \bS_{k,l} = \left( \bU^* (P \bU) \right)_{k,l} = (\buh^k)^* (P \buh^l) = 
\frac{1}{4}\left( \bD^{k+l+1} + \bD^{|k-l+1|} + \bD^{|k+l-1|} + \bD^{|k-l-1|} \right),
\label{eqn:stiffmatdata}
\ee
where $\bS_{k,l} \in \mathbb{R}^{m \times m}$, $k,l=0,\ldots,n-1$. 

Formulas (\ref{eqn:massmatdata})--(\ref{eqn:stiffmatdata}) enable us to compute the reduced 
model $\bPt$, $\bBt$ from the knowledge of the sampled data only. Indeed, for the reduced order 
propagator from (\ref{eqn:orthproj}) and (\ref{eqn:blockqr}) we obtain
\be
\bPt = \bV^* (P \bV) = \bL^{-1} \bU^* (P \bU \bL^{-*}) = \bL^{-1} \bS \bL^{-*}.
\label{eqn:romprop}
\ee
Similarly, for the reduced order transducer matrix
\be
\bBt = \bL^{-1} (\bU^* \bb) = 
\bL^{-1} \begin{bmatrix}
\bD^0 \\ \vdots \\ \bD^{n-1}
\end{bmatrix}.
\ee

We summarize the derivations above in the following algorithm.
\begin{algorithm}[Reduced order model computation from sampled data]
\label{alg:romcompute}~\\
Given the sampled data $\bD^k \in \mathbb{R}^{m \times m}$, $k=0,\ldots,2n-1$, 
to compute the reduced order propagator $\bPt \in \mathbb{R}^{mn \times mn}$ 
and transducer matrix $\bBt \in \mathbb{R}^{mn \times m}$ that solve the interpolation problem
\be
\bD^k = \bBt^* T_k(\bPt) \bBt, \quad k=0,\ldots,2n-1,
\ee
perform the following steps:
\begin{enumerate}
\item Form the mass matrix $\bM \in \mathbb{R}^{mn \times mn}$ consisting of blocks 
$\bM_{k,l} \in \mathbb{R}^{m \times m}$ for $k,l=0,\ldots,n-1,$ according to
\be
\bM_{k,l} = \frac{1}{2} \left( \bD^{k+l} + \bD^{|k-l|} \right).
\label{eqn:massmatdataalg}
\ee
\item Form the stiffness matrix $\bS \in \mathbb{R}^{mn \times mn}$ consisting of blocks 
$\bS_{k,l} \in \mathbb{R}^{m \times m}$ for $k,l=0,\ldots,n-1,$ according to
\be
\bS_{k,l} = \frac{1}{4}\left( \bD^{k+l+1} + \bD^{|k-l+1|} + \bD^{|k+l-1|} + \bD^{|k-l-1|} \right).
\label{eqn:stiffmatdataalg}
\ee
\item Perform block Cholesky decomposition of the mass matrix
\be
\bM = \bL \bL^*,
\label{eqn:blockcholromalg}
\ee
to compute the block lower triangular Cholesky factor $\bL \in \mathbb{R}^{mn \times mn}$.
\item Compute $\bPt$ and $\bBt$ from
\be
\bPt = \bL^{-1} \bS \bL^{-*}, \quad 
\bBt =  \bL^{-1} \begin{bmatrix}
\bD^0 \\ \vdots \\ \bD^{n-1}
\end{bmatrix}.
\label{eqn:projromalg}
\ee
\end{enumerate}
\end{algorithm}

Note that while the regular Cholesky factorization is defined uniquely, there is an ambiguity in defining 
block Cholesky factorization that comes from the computation of diagonal blocks $\bL_{k,k}$. The 
particular version of block Cholesky algorithm that we use to compute (\ref{eqn:blockcholromalg}) 
is given below.

\begin{algorithm}[Block Cholesky decomposition]
\label{alg:blockhol}~\\
To compute block Cholesky decomposition of a matrix $\bM \in \mathbb{R}^{mn \times mn}$
with $m \times m$ blocks, perform the following steps:\\
For $k=0,1,\ldots,n-1:$
\be
\bL_{k,k} = \left( \bM_{k,k} - \sum_{r=0}^{k-1} \bL_{k,r} \bL_{k,r}^* \right)^{1/2},
\label{eqn:lkksqrtm}
\ee
~~~~~~~~~~~~~For $l=k+1,\ldots,n-1:$
\be
\bL_{l,k} = \left( \bM_{l,k} - \sum_{r=0}^{k-1} \bL_{l,r} \bL_{k,r}^* \right) \bL_{k,k}^{-1}.
\ee

\end{algorithm}

Other choices can be used instead of (\ref{eqn:lkksqrtm}) to compute the diagonal blocks $\bL_{k,k}$, 
as long as they satisfy
\be
\bL_{k,k} \bL_{k,k}^* = \bM_{k,k} - \sum_{r=0}^{k-1} \bL_{k,r} \bL_{k,r}^*,
\quad k=0,\ldots,n-1.
\label{eqn:lkkdef}
\ee
For example, one may use eigendecomposition or regular Cholesky factorization of the right hand
side of (\ref{eqn:lkkdef}) to obtain $\bL_{k,k}$. 
No matter the choice, the reduced order propagator $\bPt$ obtained by Algorithm \ref{alg:romcompute}
has the structure established in the following lemma.

\begin{lemma}
\label{lemma:blocktri}
Reduced order propagator $\bPt$ computed by Algorithm \ref{alg:romcompute} has a block tridiagonal
structure
\be
\bPt = \begin{bmatrix}
\balpha^0 & \bbeta^1 & \bzero & \ldots & \bzero \\
(\bbeta^1)^* & \balpha^1 & \bbeta^2 & \ddots & \vdots \\
\bzero & (\bbeta^2)^* & \balpha^2 & \ddots & \vdots \\
\vdots & \ddots & \ddots & \ddots & \bbeta^{n-1} \\
\bzero & \ldots & \ldots & (\bbeta^{n-1})^* & \balpha^{n-1}
\end{bmatrix},
\ee
while the reduced order transducer matrix $\bBt$ has the structure
\be
\bBt = \begin{bmatrix}
\bbeta^0 \\ \bzero \\ \vdots \\ \bzero
\label{eqn:lemmabbt}
\end{bmatrix},
\ee
where $\balpha^k, \bbeta^k \in \mathbb{R}^{m \times m}$, $k=0,\ldots,n-1$.
\end{lemma}

The proof is given in Appendix \ref{app:lemmablocktri}. Here we note that since $\bPt$ defined by
(\ref{eqn:romprop}) is self-adjoint, so are its diagonal blocks $\balpha^k = (\balpha^k)^*$,
$k=0,\ldots,n-1$. 


\section{Backprojection imaging functional}
\label{sec:backproj}

Once the reduced order propagator $\bPt$ is computed from the sampled data with 
Algorithm \ref{alg:romcompute}, we can use it for inversion and imaging. In order to construct an 
imaging functional from $\bPt$, we consider first the Green's function $g(x,y,t)$ for the symmetrized 
problem (\ref{eqn:wavebuh})--(\ref{eqn:wavebuhinit}). Then $g(x,y,t)$ satisfies
\be
g_{tt} = \Ah g, \quad t>0, \quad x \in \Omega,
\label{eqn:greenwave}
\ee
with initial conditions
\be
g(x,y,0) = \delta(x-y) = \delta_y(x), \quad g_t(x,y,0) = 0, \quad y \in \Omega,
\label{eqn:greeninit}
\ee
where we denote by $\delta_y$ a delta function concentrated at $y$.

Then similarly to (\ref{eqn:propbu}) we can formally write the solution to 
(\ref{eqn:greenwave})--(\ref{eqn:greeninit}) at time $t = \tau$ as
\be
g(x,y,\tau) = [P \delta_y](x) = \left< \delta_x, P \delta_y \right>,
\label{eqn:greendotprod}
\ee
in terms of the inner product (\ref{eqn:l2innerprod}). Alternatively, if we treat $\delta_x$ and 
$\delta_y$ as row-vector-valued functions with a single column, then we can use the notation of 
(\ref{eqn:matprod}) to write
\be
g(x,y,\tau) = \delta_x^* [P \delta_y].
\label{eqn:greenxtpy}
\ee

If we could probe the propagator $P$ with delta functions, this would allow us to image the Green's
function (\ref{eqn:greenxtpy}). This is not possible if we only have access to the sampled data $\bD^k$
and thus to $\bPt$. However, we can probe $P$ with other functions. Let us consider 
\be
\bV \bV^* = 
\left[ \bv^0, \bv^1, \ldots, \bv^{n-1} \right]
\begin{bmatrix} (\bv^0)^* \\ (\bv^1)^* \\ \vdots \\ (\bv^{n-1})^* \end{bmatrix},
\ee
a function of two variables
\be
[\bV \bV^*](x,y) = \sum_{k=0}^{n-1} \bv^k(x) (\bv^k(y))^*.
\ee
As we will see later, $[\bV \bV^*](x,y)$ approximates $\delta(x-y)$ in a certain sense. Thus, we 
introduce an approximate Green's function
\be
\gt(x,y,\tau) = [\bV \bV^*](x, \;\cdot\;) \left[ P [\bV \bV^*](\;\cdot\;, y) \right] \approx g(x,y,\tau)
\label{eqn:gtdef}
\ee
by replacing the delta functions in (\ref{eqn:greenxtpy}) with $\bV \bV^*$. Then we can transform
\be
\gt(x,y,\tau) = [\bV(x)] \left( \bV^* [P \bV] \right) [\bV^*(y)] = [\bV(x)] \bPt [\bV^*(y)],
\label{eqn:gtmat}
\ee
or explicitly
\be
\gt(x,y,\tau) = \sum_{k,l = 0}^{n-1} [\bv^k(x)] \bPt_{k,l} [\bv^l(y)]^*,
\label{eqn:gtsum}
\ee
where the standard rules for matrix-vector multiplication apply.

While $\bPt$ can be computed from the sampled data, it is not possible to obtain the orthonormal
basis functions $\bV$ from the knowledge of $\bD^k$ only. Thus, we need to make another 
approximation in (\ref{eqn:gtmat})--(\ref{eqn:gtsum}). Recall that when solving an imaging problem
we are given the approximate knowledge of the travel times between pairs of points in $\Omega$
in the form of a smooth kinematic model $c_o(x)$. For $c_o(x)$ we can compute the symmetrized 
wavefield snapshots $\buh_{(o)}^k(x)$ in the whole domain $\Omega$, as well as their 
orthogonalization
\be
\bV_{(o)}(x) = \left[ \bv_{(o)}^0(x), \bv_{(o)}^1(x), \ldots, \bv_{(o)}^{n-1}(x) \right],
\ee
and the projected propagator
\be
\bPt_{(o)} = \bV_{(o)}^* \left[ P_{(o)} \bV_{(o)} \right].
\ee
Hereafter we use the subscript $(o)$ to refer to all quantities associated with $c_o(x)$ including
the Green's function $g_{(o)}(x,y,\tau)$ and its approximation $\gt_{(o)}(x,y,\tau)$ defined
in (\ref{eqn:gtdef}).

If $c_o(x)$ approximates the kinematics of the problem with unknown desired $c(x)$ sufficiently well,
we may expect that the orthogonalized snapshots for the two problems are also close
\be
\bV(x) \approx \bV_{(o)}(x).
\ee
Particular reasons for expecting such behavior are discussed later. Then we can make another 
approximation in (\ref{eqn:gtdef})--(\ref{eqn:gtmat}) to obtain
\be
\gt(x,y,\tau) \approx  [\bV_{(o)}(x)] \bPt [\bV^*_{(o)}(y)],
\ee
where all terms can be computed either from the knowledge of the sampled data or the kinematic
model.

Finally, we have all the components to define the \emph{ROM backprojection imaging functional}
\be
\cI_{BP}(x) = [\bV_{(o)}(x)] \left( \bPt - \bPt_{(o)} \right) [\bV^*_{(o)}(x)], 
\quad x \in \Omega,
\label{eqn:imgfun}
\ee
that approximates the difference between the Green's functions for the unknown medium $c(x)$ and 
the kinematic model $c_{(o)}(x)$, that is
\be
\cI_{BP}(x) \approx \gt(x,x,\tau) - \gt_{(o)}(x,x,\tau) \approx g(x,x,\tau) - g_{(o)}(x,x,\tau), 
\quad x \in \Omega.
\label{eqn:imgapprox}
\ee
If in the vicinity of a point $x \in \Omega$ the unknown medium $c(x)$ contains a reflector, then the difference
(\ref{eqn:imgapprox}) will be large. Otherwise, $c(x) \approx c_{(o)}(x)$ and $\cI_{BP}(x) \approx 0$.

The use of term ``backprojection'' should not be confused with the classical method in computed 
tomography \cite{kuchment2013radon}. The (filtered) backprojection formula for inversion of 
Radon transform in computed tomography is linear in the input data, while the computation of 
$\cI_{BP}$ from the data $\bD^k$, $k=0,\ldots,2n-1,$ in Algorithm \ref{alg:romcompute} is a 
highly nonlinear procedure due to block Cholesky factorization (\ref{eqn:blockcholromalg}) and 
inversion (\ref{eqn:projromalg}).

\section{Interpretation in terms of Schr\"{o}dinger wave equation}
\label{sec:interp}

Here we present an interpretation of the imaging formula (\ref{eqn:imgfun}) that goes beyond the
approximate Green's function perturbation (\ref{eqn:imgapprox}). Note that the imaging procedure 
in section \ref{sec:backproj} makes a distinction between the kinematic part of $c(x)$ that is assumed
to be well approximated by $c_o(x)$, and the reflective part of $c(x)$ that is imaged by $\cI_{BP}$.
Meanwhile, the operator (\ref{eqn:continuumA}) has a single coefficient that does not allow for a 
clear separation between the kinematic and reflective parts. Thus, instead of the wave equation
(\ref{eqn:waveu}) with operator (\ref{eqn:continuumA}), in a slight abuse of notation we consider 
the wave equation
\be
u_{tt} = A u,
\label{eqn:waveA}
\ee
with operator
\be
A = \sigma(x) c(x) \nabla\cdot \Big[ \frac{c(x)}{\sigma(x)} \nabla \Big],
\label{eqn:Acs}
\ee
where $c(x)$ is the wavespeed and $\sigma(x)$ is the \emph{acoustic impedance}. 

The operator (\ref{eqn:Acs}) is symmetric with respect to the $1/[\sigma(x) c(x)]$-inner product and 
nonnegative definite (positive definite for a bounded $\Omega$). Thus, following section \ref{sec:sym} 
we use $u(x, t) = \sqrt{\sigma(x) c(x)} \widehat{u}(x, t)$ to transform (\ref{eqn:waveA}) into
$ \widehat{u}_{tt} = \widehat{A} \widehat{u}, $ where the symmetrized operator $\widehat{A}$ 
is given by
\be
\widehat{A} = \sqrt{\sigma(x) c(x)} \nabla \cdot \Big[ \frac{c(x)}{\sigma(x)} \nabla \sqrt{\sigma(x) c(x)} \Big].
\ee
Furthermore, we can transform $\widehat{A}$ to Schr\"{o}dinger form
\be
\cA =  \cC - qc \cI,
\label{propop}
\ee
with
\be
\cC = \sqrt{c} \nabla \cdot \Big( c \nabla \sqrt{c} \Big),
\ee
and the Schr\"{o}dinger potential
\be 
q = \sqrt{\sigma} \nabla \cdot \left( c \nabla \frac{1}{\sqrt{\sigma}} \right),
\label{eqn:q}
\ee
where we omit the dependency on $x$ of $c(x)$, $\sigma(x)$ and $q(x)$ to simplify the notation.

To obtain the separation between the kinematic and reflective parts we assume that the wavespeed 
is known exactly $c(x) \equiv c_o(x)$ and we only image the perturbations of the Schr\"{o}dinger potential
\be
\delta q = q - q_o = 
\sqrt{\sigma} \nabla \cdot \left( c \nabla \frac{1}{\sqrt{\sigma}} \right) - 
\sqrt{\sigma_o} \nabla \cdot \left( c_o \nabla \frac{1}{\sqrt{\sigma_o}} \right),
\ee
where $\sigma_o(x)$ is the known background impedance.
Note that the knowledge of $c(x) \equiv c_o(x)$ is a rather reasonable assumption. Discontinuities of 
$c(x)$ do not produce significant reflections observable for moderate apertures 
\cite[Chapter 6]{symes1995mathematics}, so typically the inexactness of $c(x)$ does not introduce 
significant topological errors in the image, and may only result in a slight image deformation.

Similarly to (\ref{eqn:greenwave}) we introduce the Green's functions $g(x,y,t)$ and $g_{(o)}(x,y,t)$ for
the unknown and reference media respectively. They solve
\begin{align}
g_{tt} & = \cA g, \quad g(x,y,0) = \delta(x - y), \quad g_t(x,y,0) = 0, \label {eqn:green} \\
g_{(o),tt} & = \cA_o g_{(o)}, \quad g_{(o)} (x,y,0) = \delta(x - y), \quad g_{(o),t}(x,y,0) = 0, \label {eqn:green0}
\end{align}
where $\cA_o =  \cC_o - q_o c_o \cI$, and $\cC_o = \sqrt{c_o} \nabla \cdot \Big( c_o \nabla \sqrt{c_o} \Big)$.

According to (\ref{eqn:imgapprox}) the functional $\cI_{BP}(x)$ is related to the perturbation 
\be
\delta g(x,y,t) = g(x,y,t) - g_{(o)}(x,y,t).
\ee
It follows from (\ref{eqn:green})--(\ref{eqn:green0}) that it satisfies approximately
\be
\delta g_{tt}(x,y,t) = \cA_o \delta g(x,y,t) - c_o(x) \delta q(x) g_{(o)} (x,y,t), 
\quad \delta g(x,y,0) = \delta g_t(x,y,0) = 0,
\label{eqn:wavedeltag}
\ee
if we assume 
\be
\delta q(x) \; \delta g(x,y,t) \ll 1.
\label{eqn:dqdg}
\ee 
If we also introduce the Green's function $G_{(o)} (x,y,t)$ solving
\be
G_{(o),tt} = \cA_o G_{(o)}, \quad G_{(o)} (x,y,0) = 0, \quad G_{(o),t} (x,y,0) = \delta(x - y), 
\label{eqn:green2}
\ee
we can express the solution of (\ref{eqn:wavedeltag}) as
\be
\delta g(x,y,t) = - \int_{0}^{t} \int_{\Omega} \delta q (z) c_o(z) g_{(o)}(z, y, s) G_{(o)}(x, z, t-s) dz ds,
\ee
and so
\be
\cI_{BP}(x) \approx \delta g(x,x,\tau) \approx 
- \int_0^\tau \int _{\Omega} \delta q(z) c_o(z) g_{(o)}(z, x, s) G_{(o)}(x, z, \tau-s) dz ds.
\label{eqn:Ibpq}
\ee
Thus, due to causality of (\ref{eqn:green0}) and (\ref{eqn:green2}), $\cI_{BP}(x)$ yields $\delta q$ 
averaged in the travel time ball of radius $\tau/2$ centered at $x$.
Note that since (\ref{eqn:Ibpq}) only relies on the solutions of (\ref{eqn:wavedeltag}) for $t < \tau$, 
the assumption (\ref{eqn:dqdg}) is valid for $\tau$ small enough compared to the scale of 
inhomogeneities of $q$ (and thus $\sigma$).



\section{Numerical study and discussion}
\label{sec:num}

The proposed imaging method differs significantly from existing approaches. In this section we 
illustrate numerically and discuss its various aspects that distinguish it from conventional methods.

\subsection{Numerical model}
\label{subsec:nummodel}

To study the performance of the proposed imaging method we consider the following numerical model.
Let $\Omega = [0,3]\; km \times [0,3]\; km$. The inaccessible boundary $\cB_I$ consists of the left, right
and bottom sides of the square, while $\cB_A$ is the top side. The accessible boundary $\cB_A$
supports $m=32$ point-like transducers, i.e. 
\be
e_j(x) = \delta(x - x_j),
\label{eqn:edelta}
\ee 
where $x_j \in \cB_A$, $j=1,\ldots,m,$ are the transducer locations.

\begin{figure}[t!]
\centering
\begin{tabular}{cc}
$c(x)$ & $c_o(x)$ \\
\includegraphics[width=0.4\textwidth]{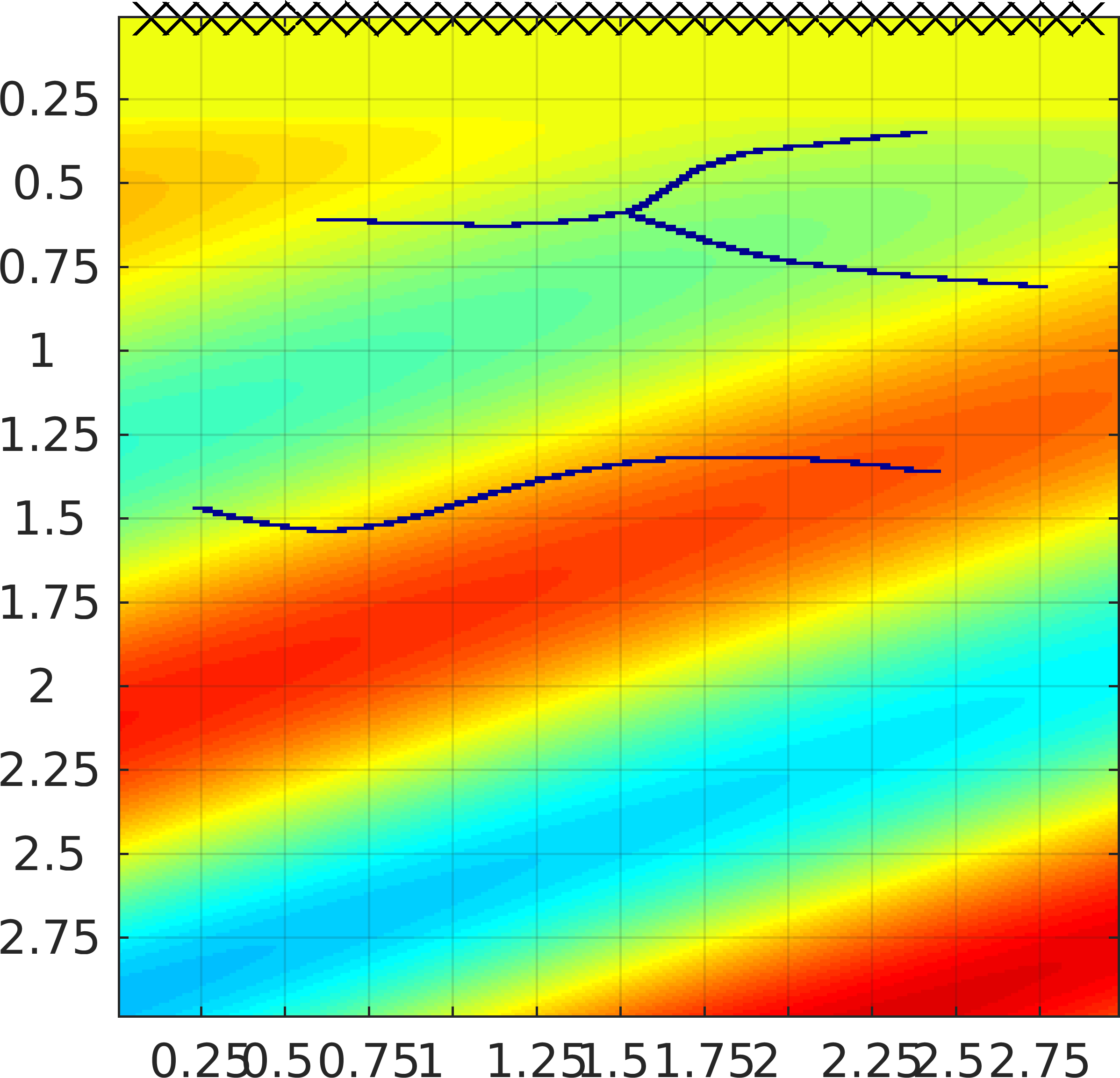} & 
\includegraphics[width=0.4\textwidth]{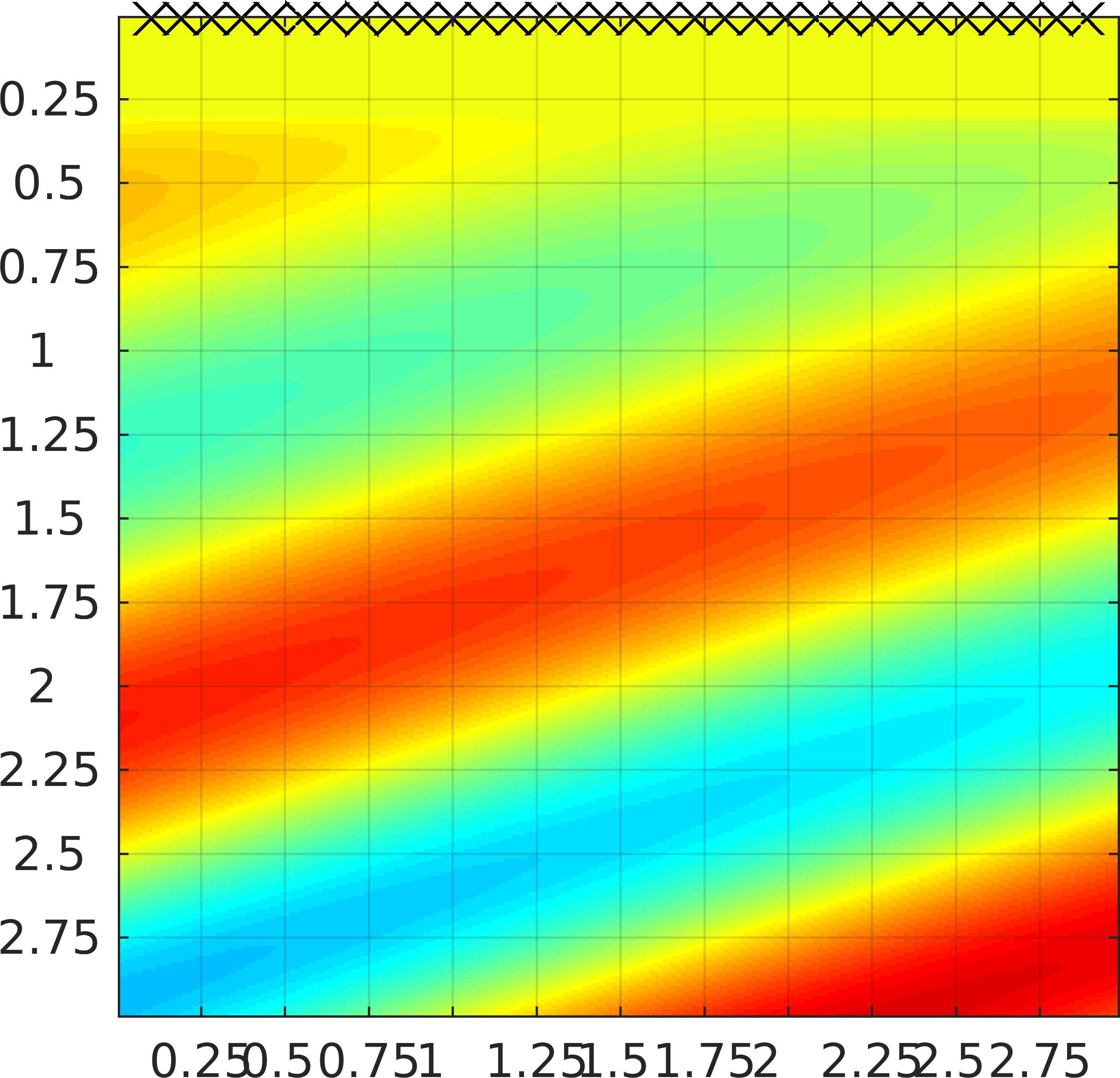}
\end{tabular}
\caption{Acoustic velocity with two extended reflectors $c(x)$ and the kinematic model $c_o(x)$
used in the numerical study. Locations $x_j$, $j=1,\dots,m,$ of $m=32$ point-like transducers are 
black $\times$. All distances in $km$, velocities in $km/s$.}
\label{fig:cfrac4}
\end{figure}

We use an acoustic velocity model $c(x)$ consisting of a smooth background and two extended 
reflectors one of which is branching. The smooth part of $c(x)$ is used as the kinematic model $c_o(x)$.
Both $c(x)$ and $c_o(x)$ are shown in Figure \ref{fig:cfrac4}. Reflectors are relatively high contrast 
with velocity of $1\; km/s$ inside the reflecting inclusion, while the surrounding background velocities
are about $2-3\; km/s$ depending on the location. 

The wave equation is discretized with second order finite differences on a uniform tensor product
grid with $300 \times 300$ nodes spaced $10 \; m$ apart. Each of the two reflectors has thickness 
of $20 \; m$, two grid nodes.

The source wavelet is chosen so that its Fourier transform (\ref{eqn:qfourier}) is
\be
\widetilde{q}(s^2) = e^{- \frac{\sigma^2 s^2}{2}},
\label{eqn:qgauss}
\ee
where $\sigma$ determines the characteristic duration of the source wavelet $q(t)$ in time. 
As discussed in \cite{druskin2015direct}, it is essential for the stability of ROM computation for 
the temporal sampling interval $\tau$ to be consistent with the Nyquist sampling limit for the source 
wavelet. Since the Gaussian (\ref{eqn:qgauss}) does not have a sharp frequency cutoff, its Nyquist 
limit is not precisely defined, but the choice $\tau \approx \sigma$ works well in practice. In particular,
in all numerical examples in this section we take 
$\tau = \sqrt{3} \sigma / 2 \approx 0.87 \sigma = 1.5 \cdot 10^{-2} \; s$. While it is mostly 
pointless to talk about the ``frequency'' of the broadband (infinite-bandwidth) source wavelet
(\ref{eqn:qgauss}), we can characterize the source in terms of the \emph{sampling frequency}
$\omega_\tau = 1 / \tau$. Here $\omega_\tau = 66.66 \; Hz$. The data is measured at $2n = 130$
time instants $t_k = k \tau$, $k = 0,\ldots,2n-1$ with a terminal time of $1.715 \; s$.

\subsection{Orthogonalized wavefield snapshots}
\label{subsec:orthsnap}

\begin{figure}
\begin{tabular}{cccc}
$\buh_{(o),j}^k$ & $\bv_{(o),j}^k$ & $\buh_j^k$ & $\bv_j^k$ \\
\includegraphics[width=0.225\textwidth]{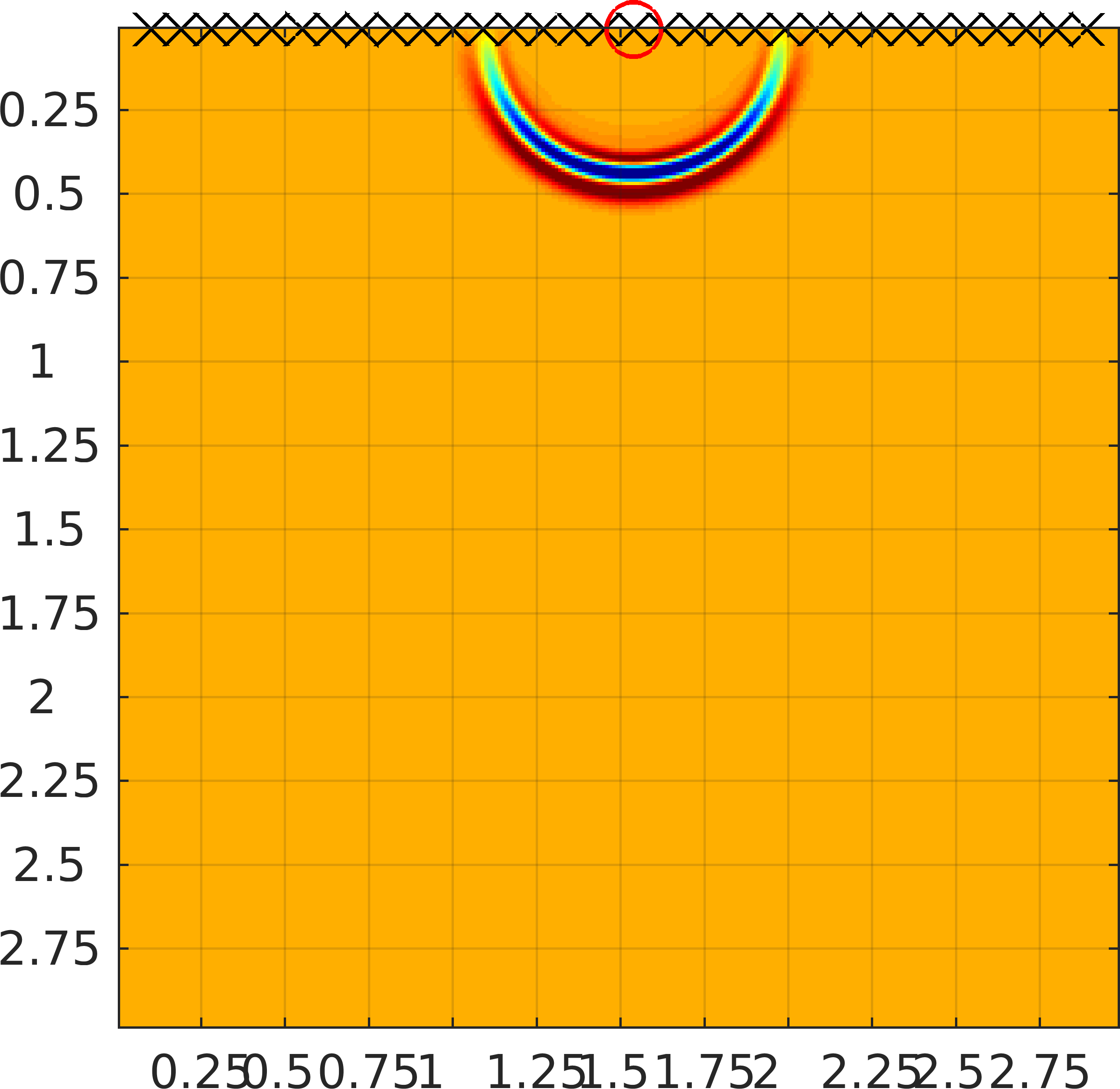} & 
\includegraphics[width=0.225\textwidth]{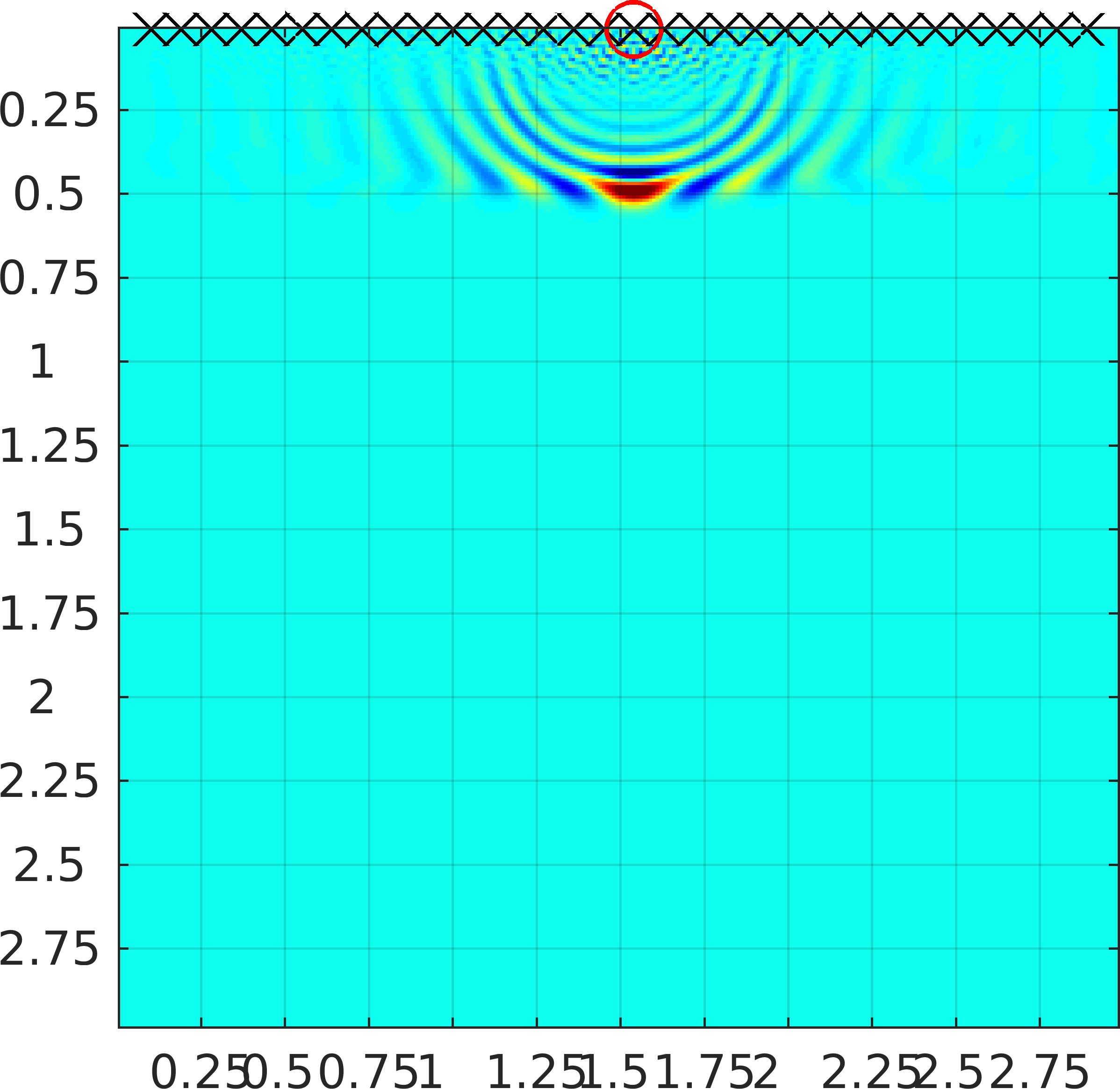} &
\includegraphics[width=0.225\textwidth]{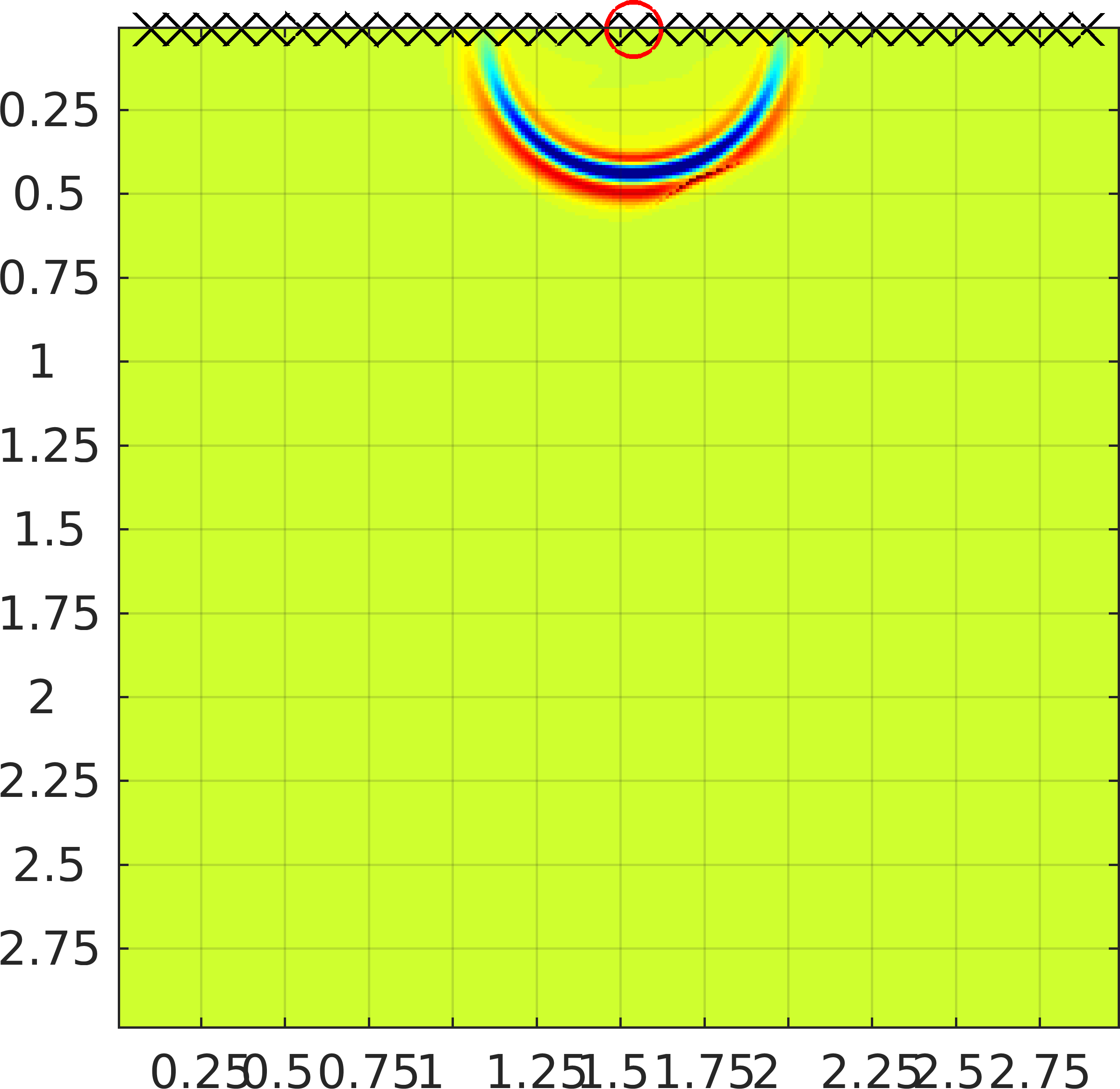} & 
\includegraphics[width=0.225\textwidth]{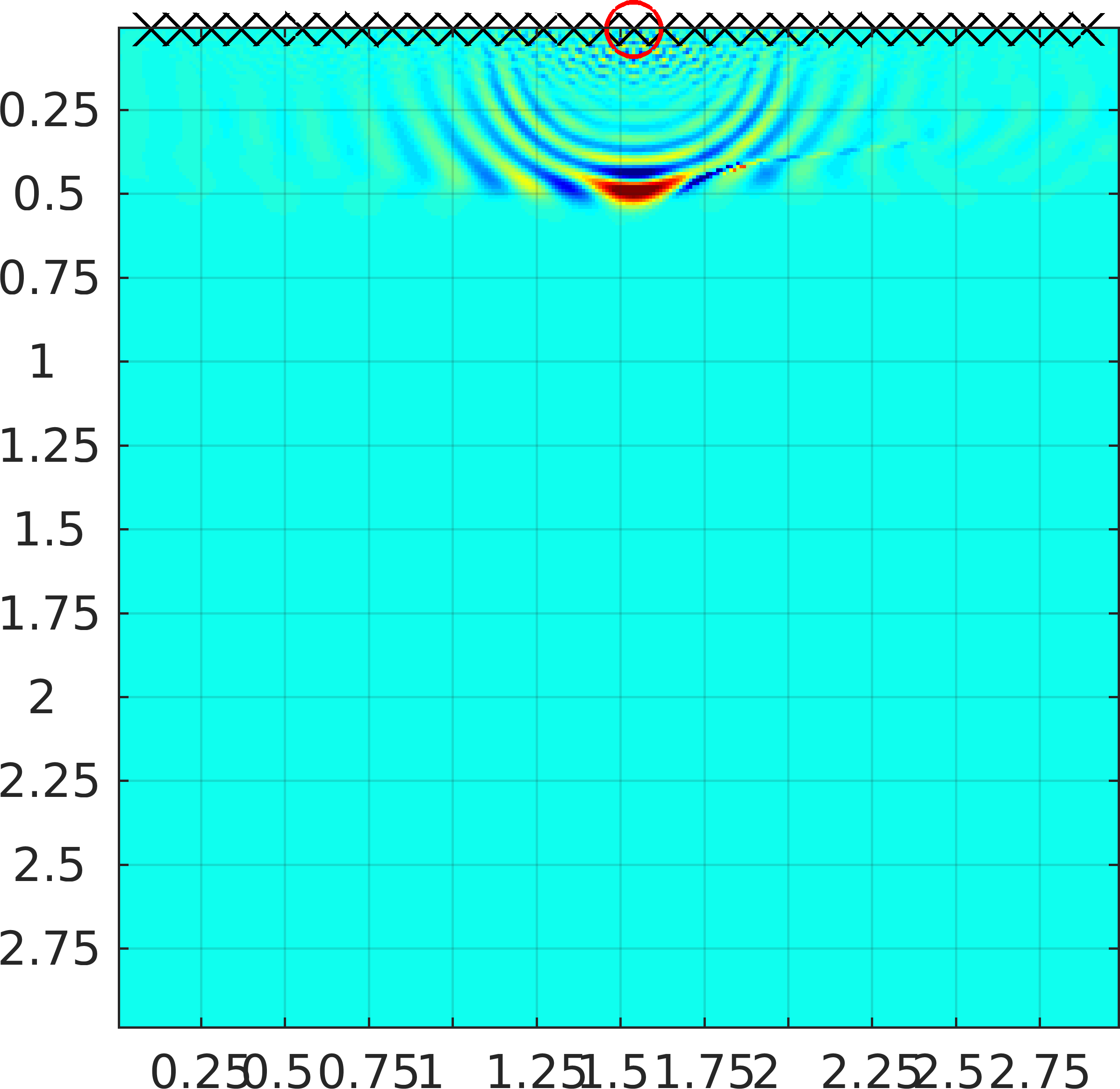} \\
\includegraphics[width=0.225\textwidth]{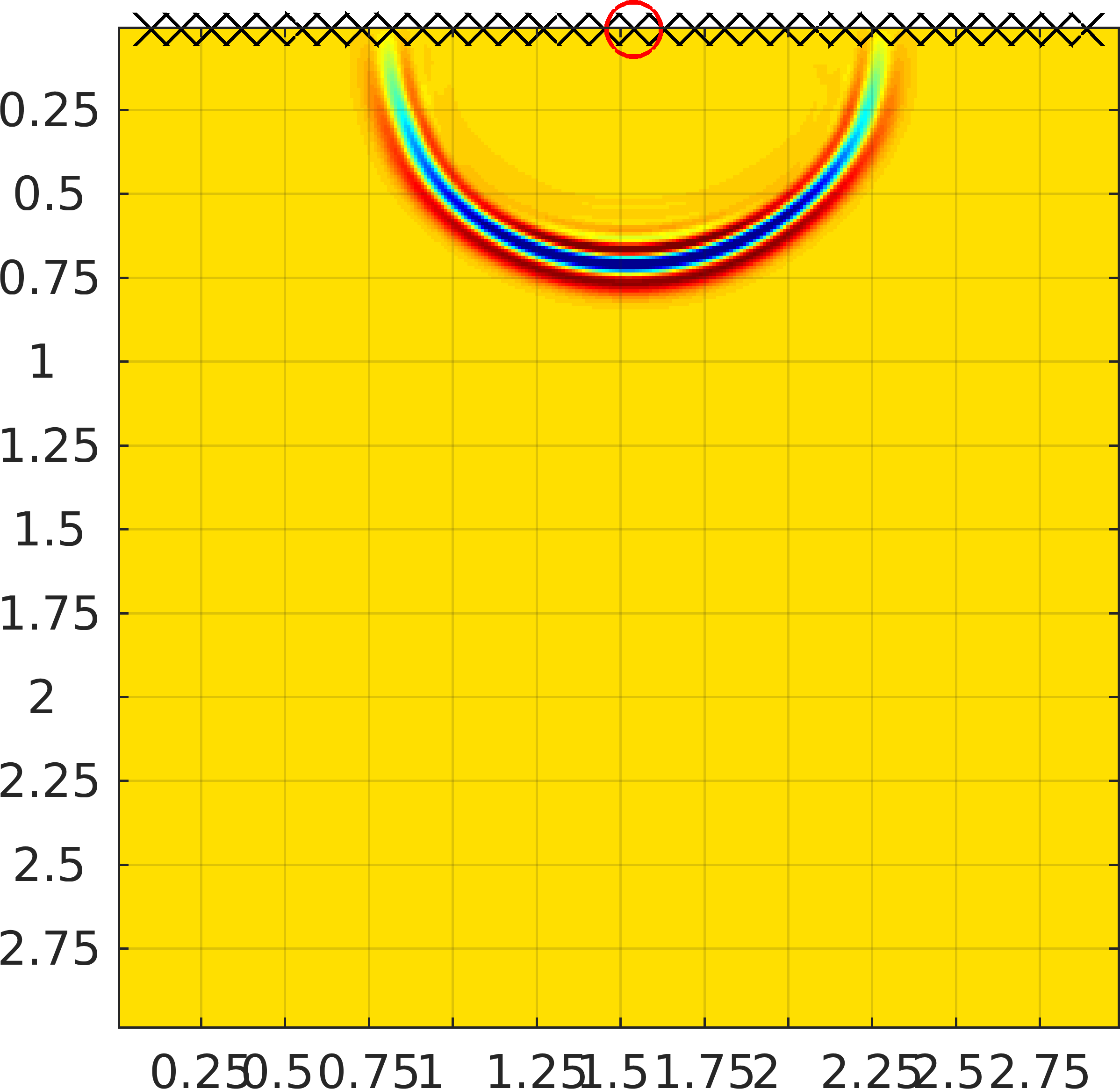} & 
\includegraphics[width=0.225\textwidth]{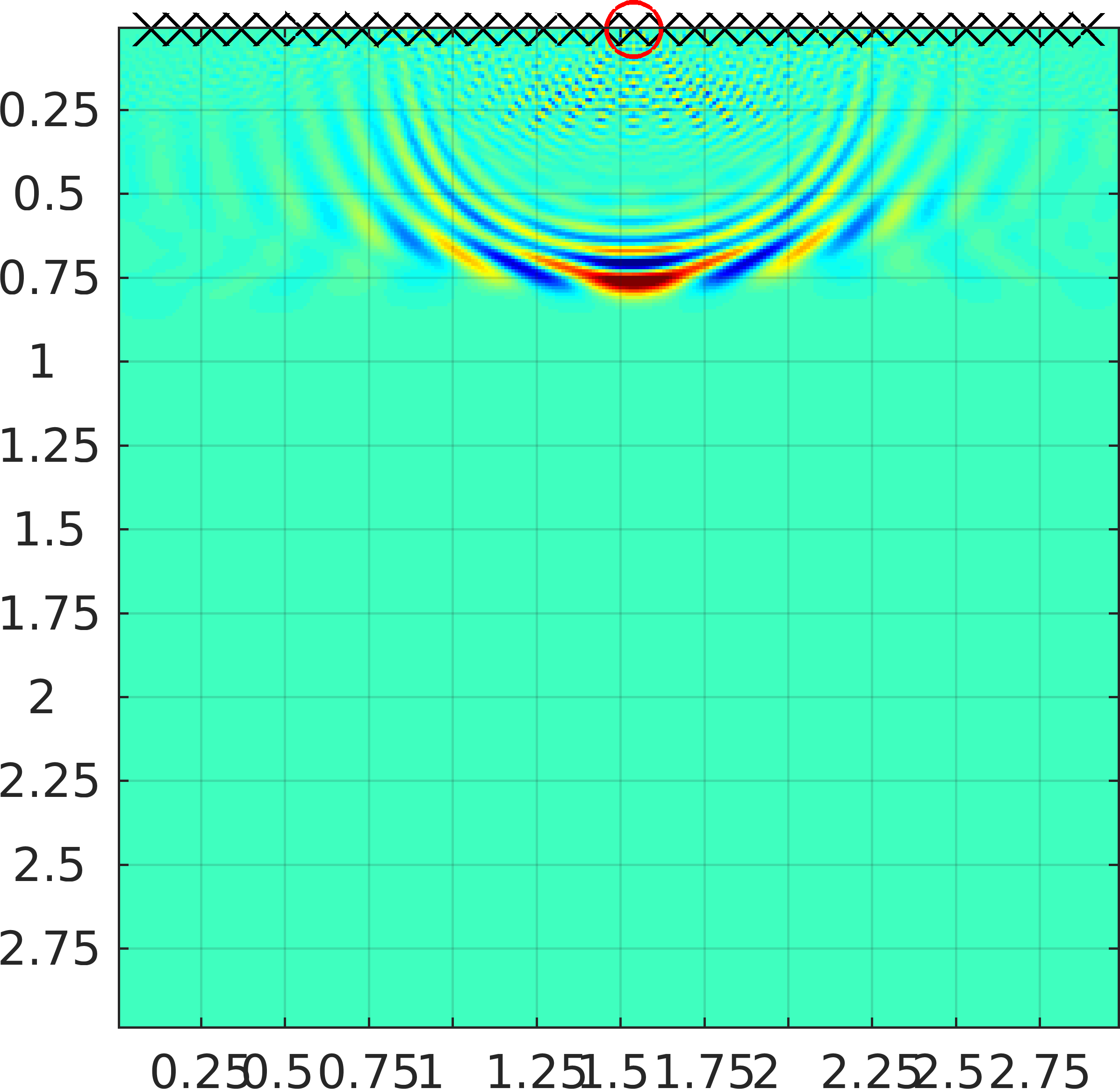} &
\includegraphics[width=0.225\textwidth]{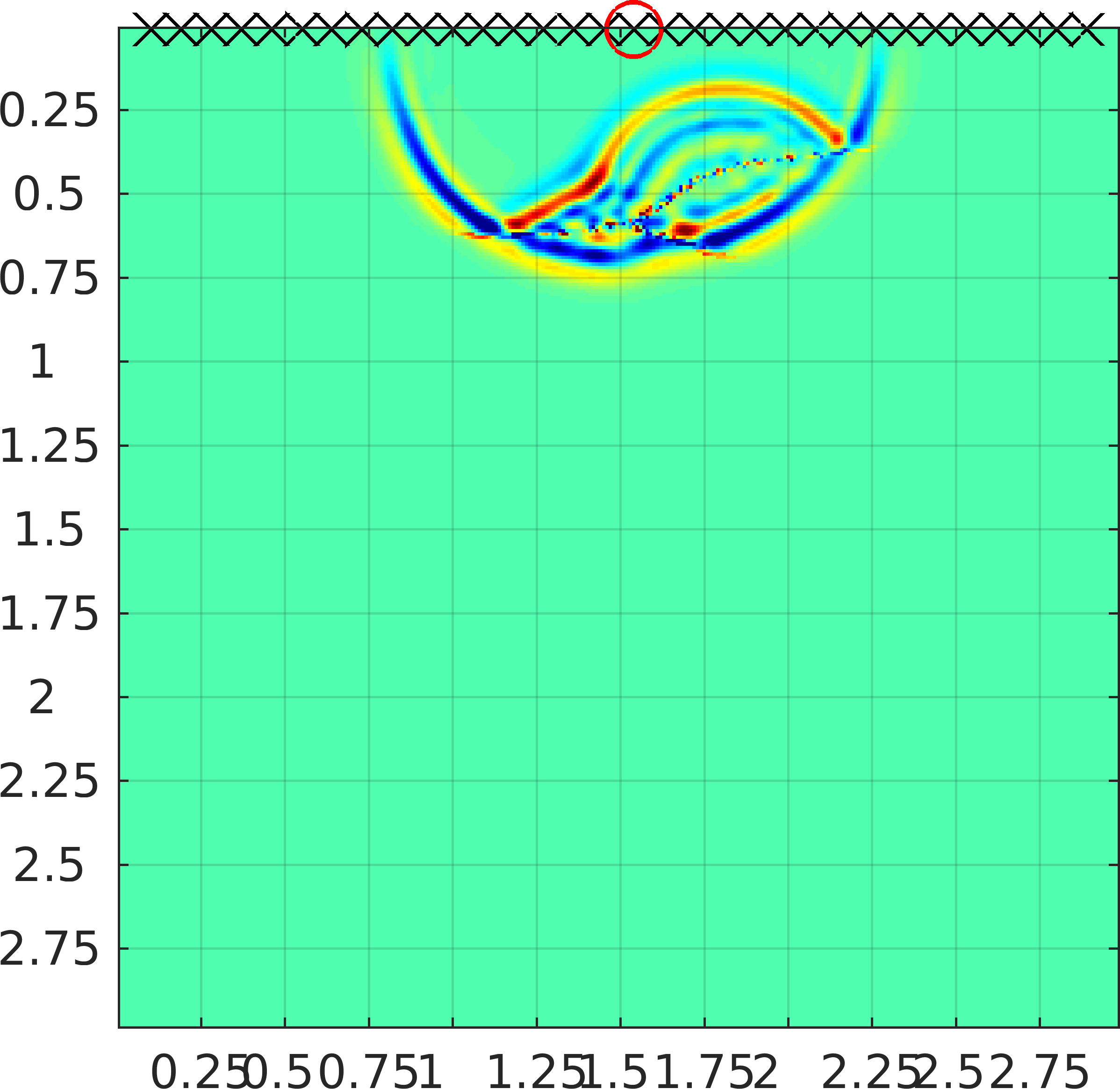} & 
\includegraphics[width=0.225\textwidth]{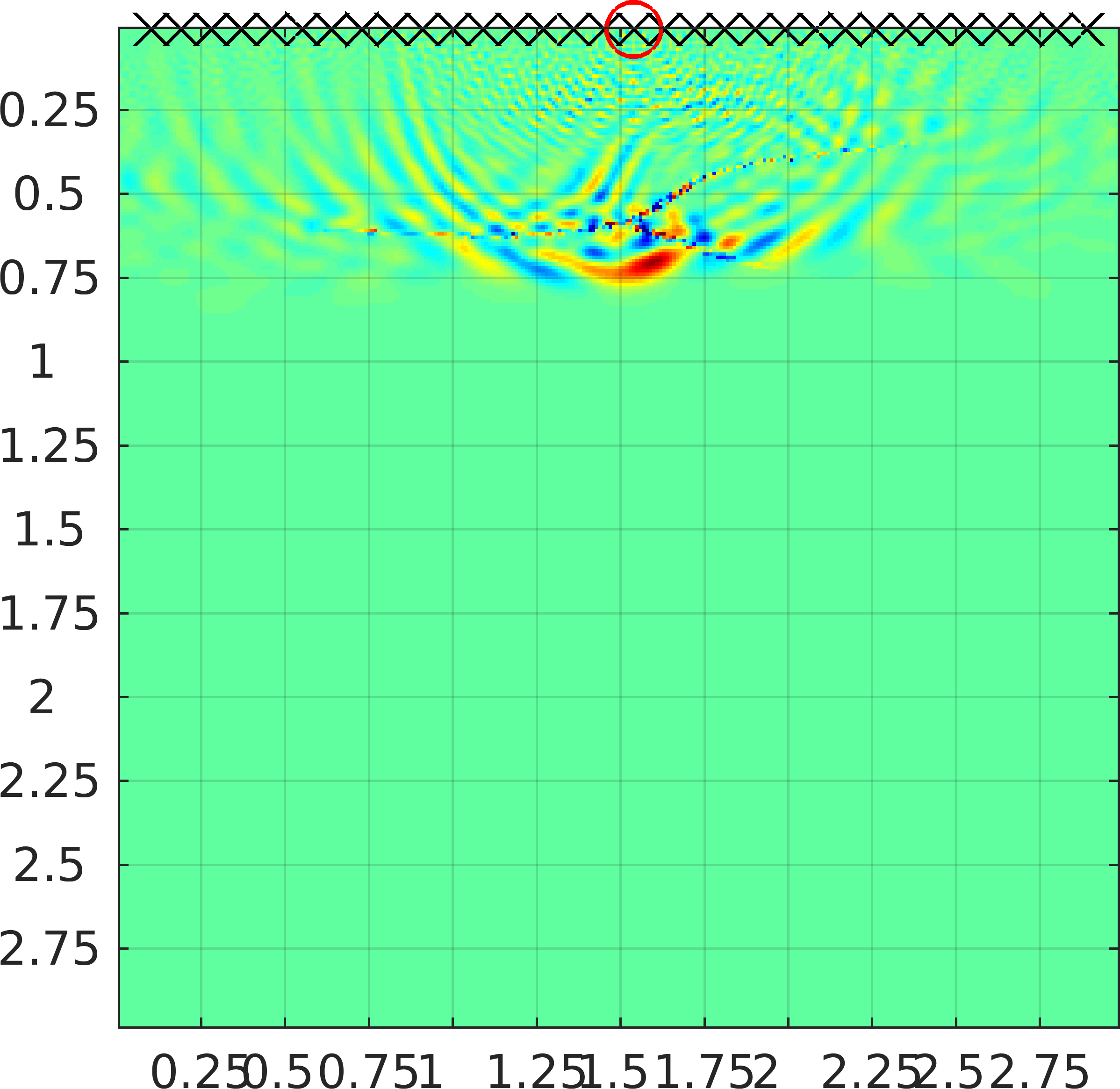} \\
\includegraphics[width=0.225\textwidth]{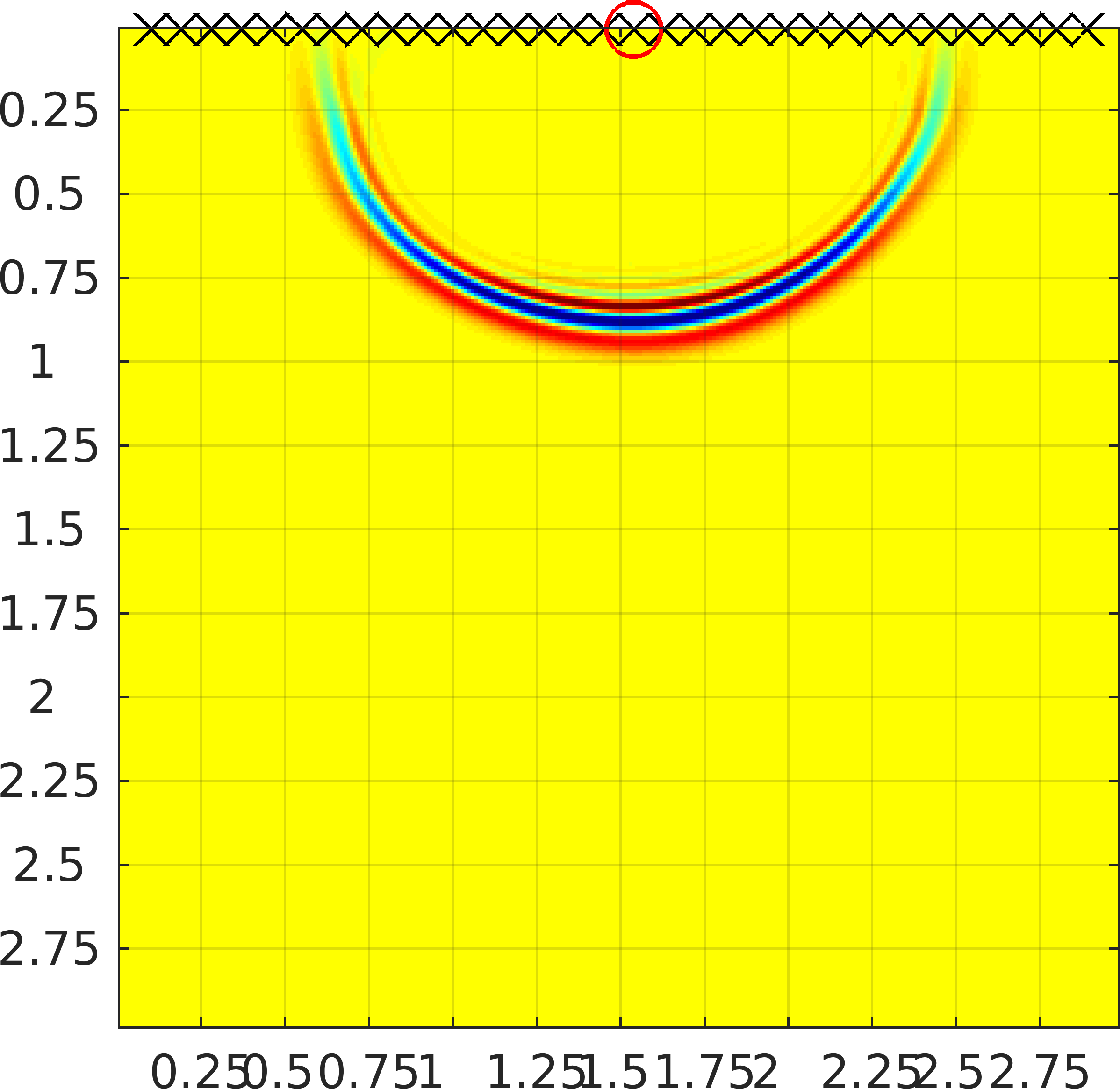} & 
\includegraphics[width=0.225\textwidth]{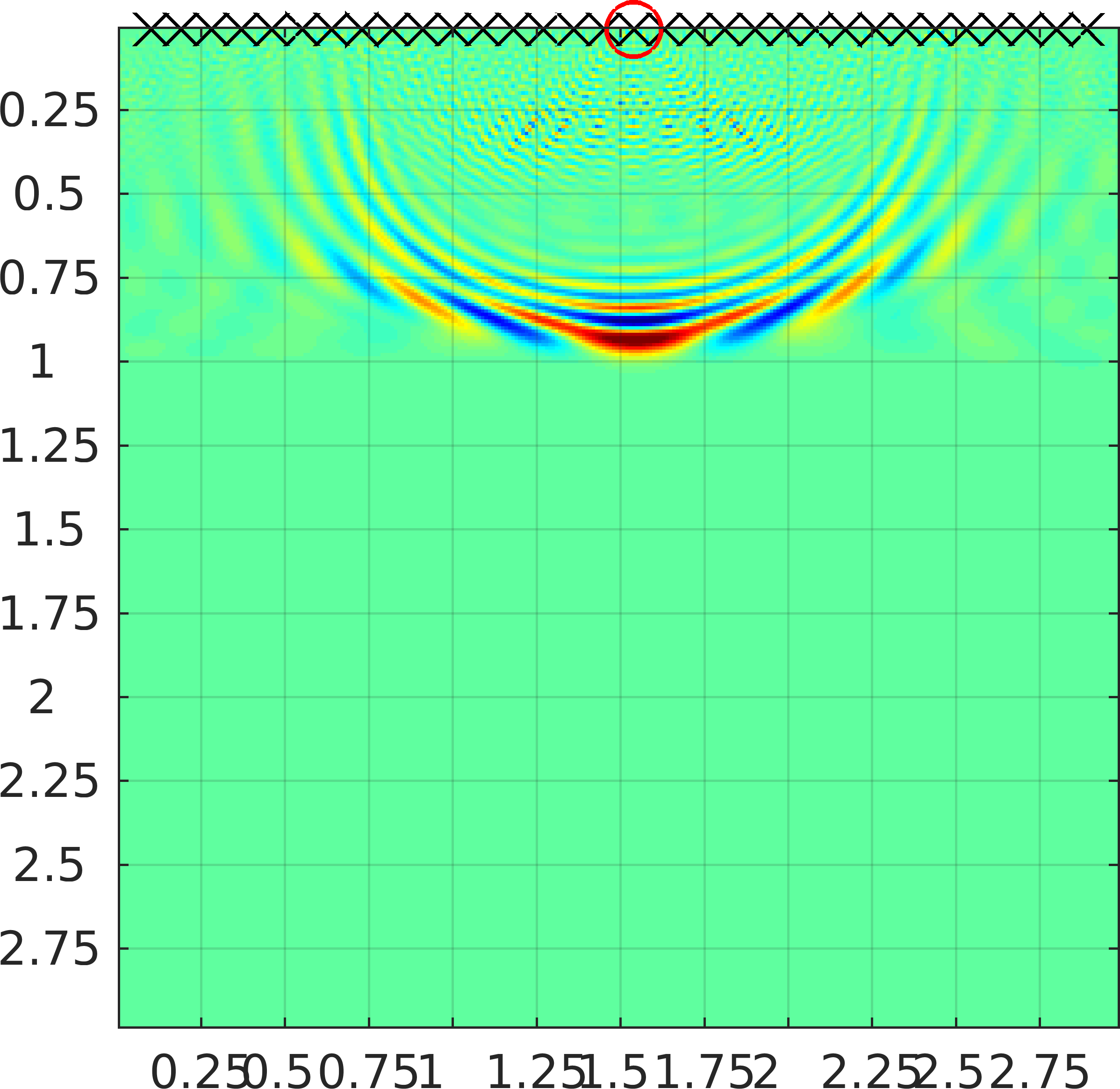} &
\includegraphics[width=0.225\textwidth]{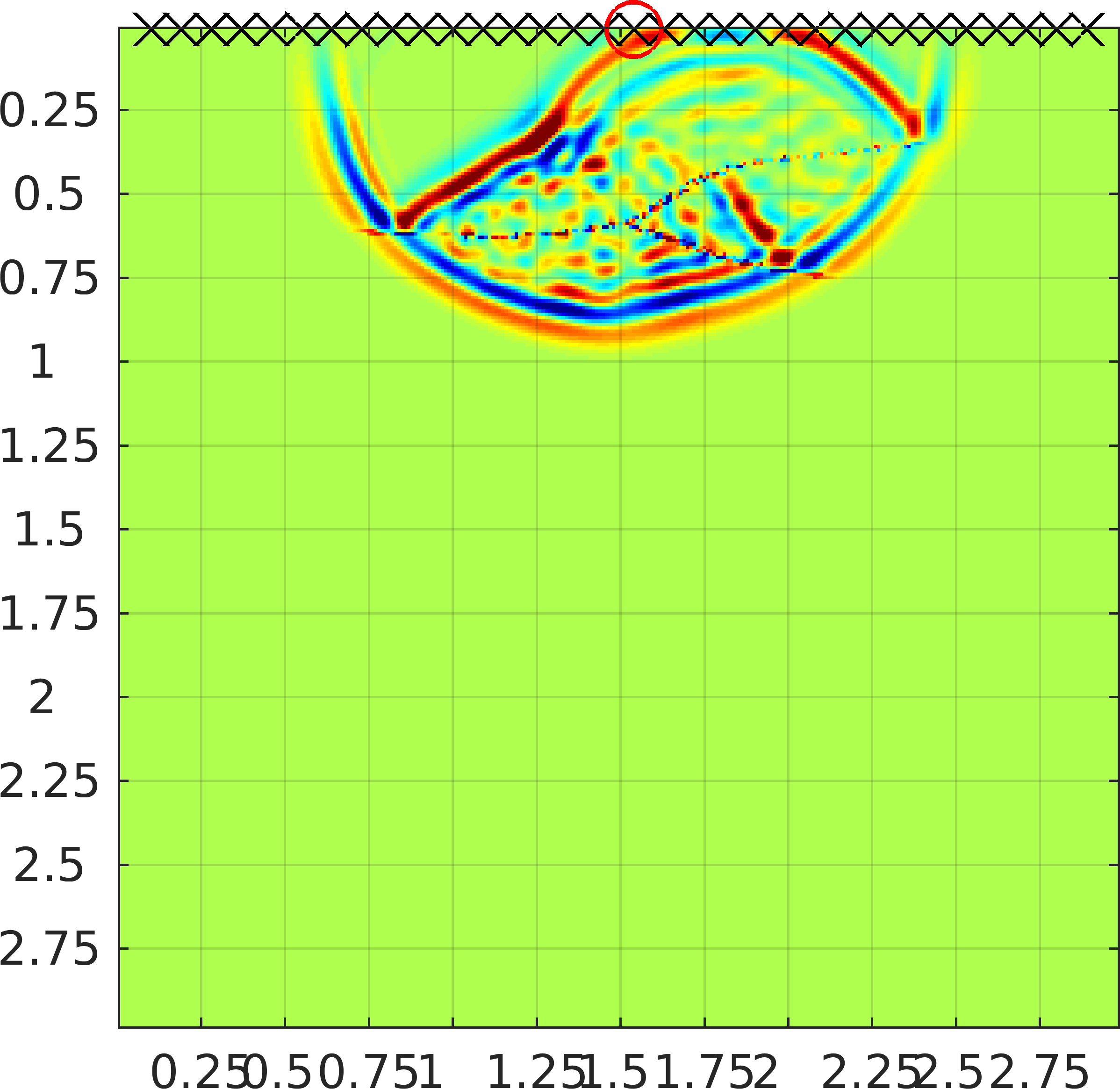} & 
\includegraphics[width=0.225\textwidth]{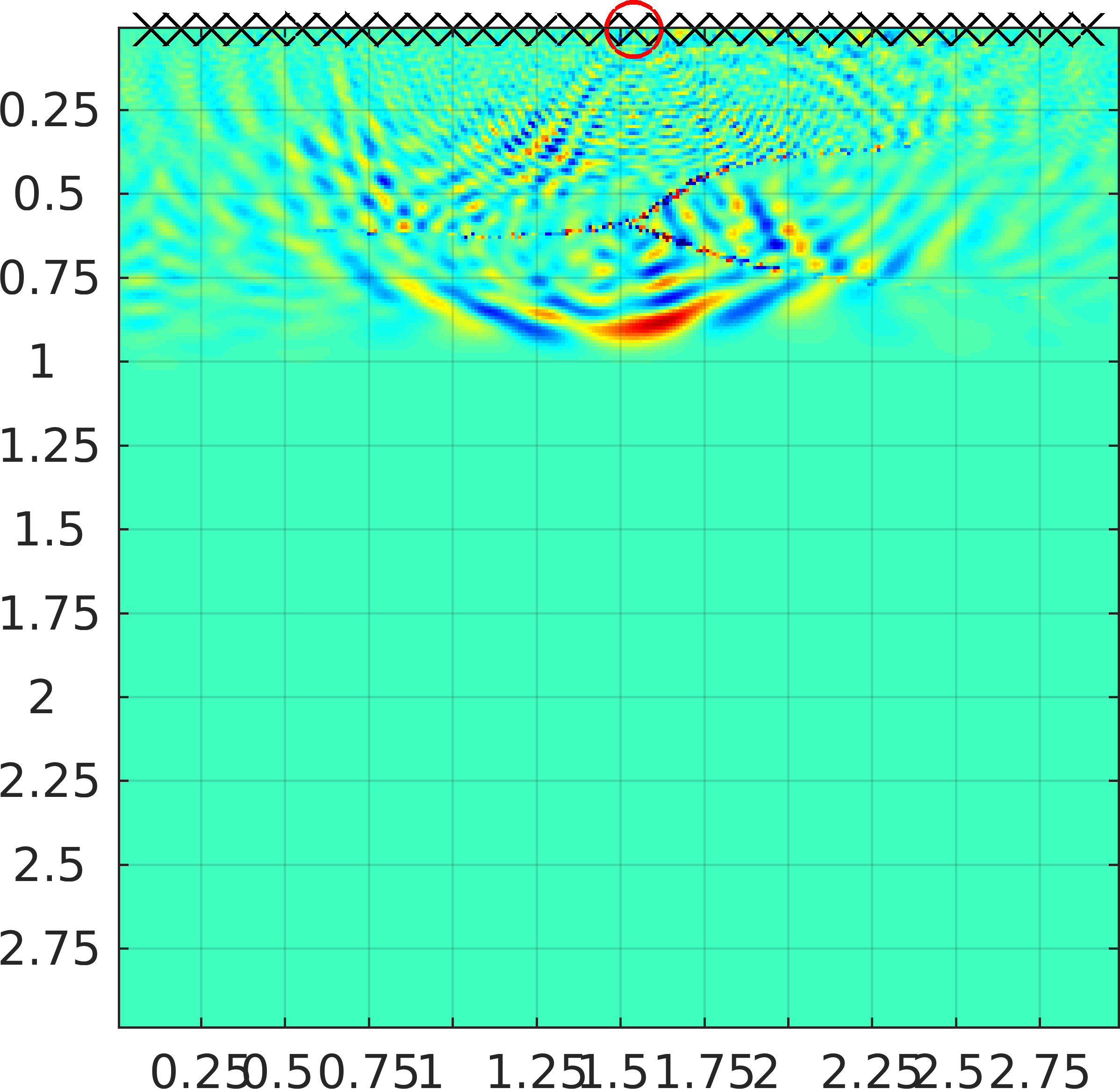} \\
\includegraphics[width=0.225\textwidth]{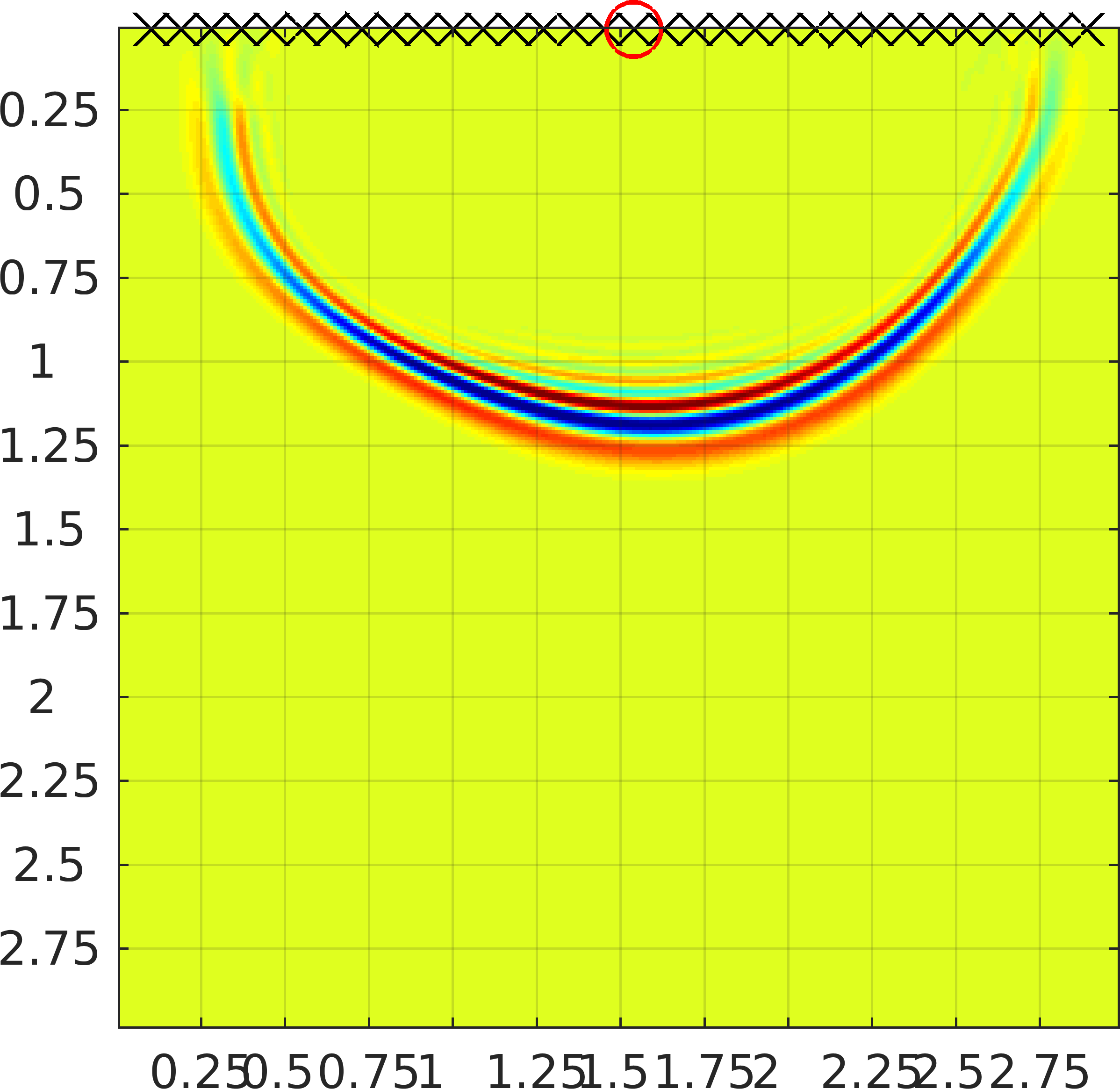} & 
\includegraphics[width=0.225\textwidth]{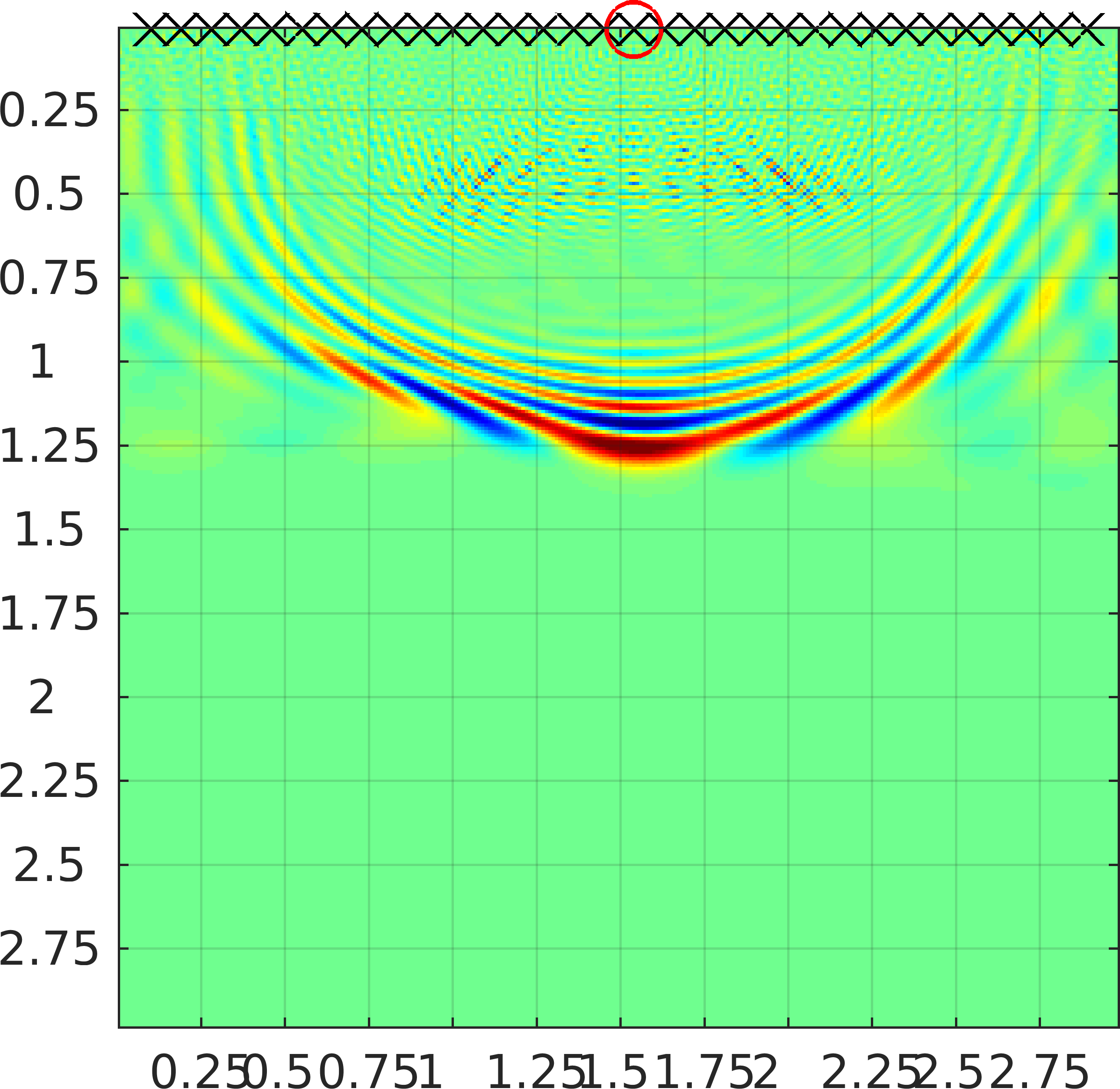} &
\includegraphics[width=0.225\textwidth]{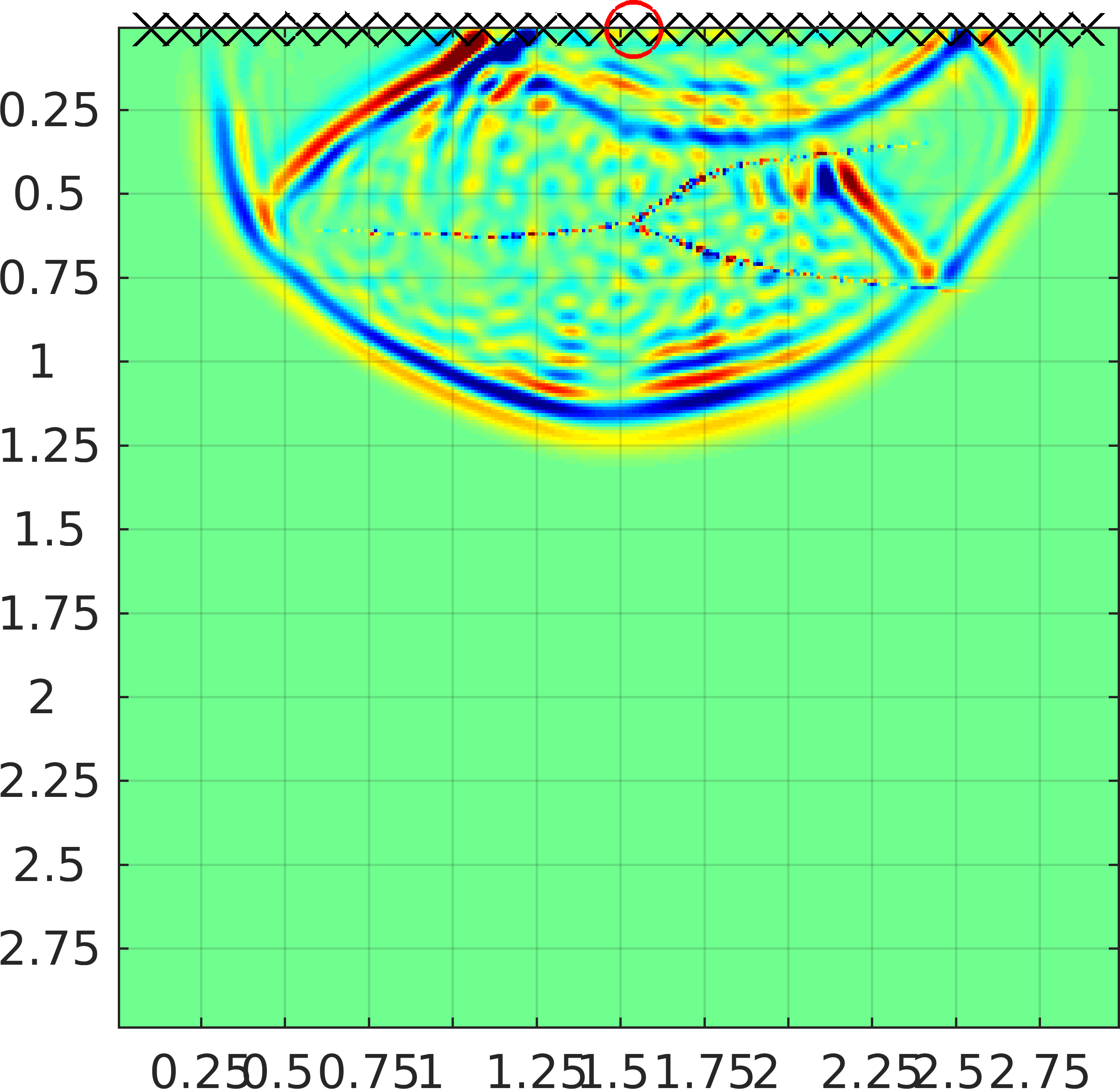} & 
\includegraphics[width=0.225\textwidth]{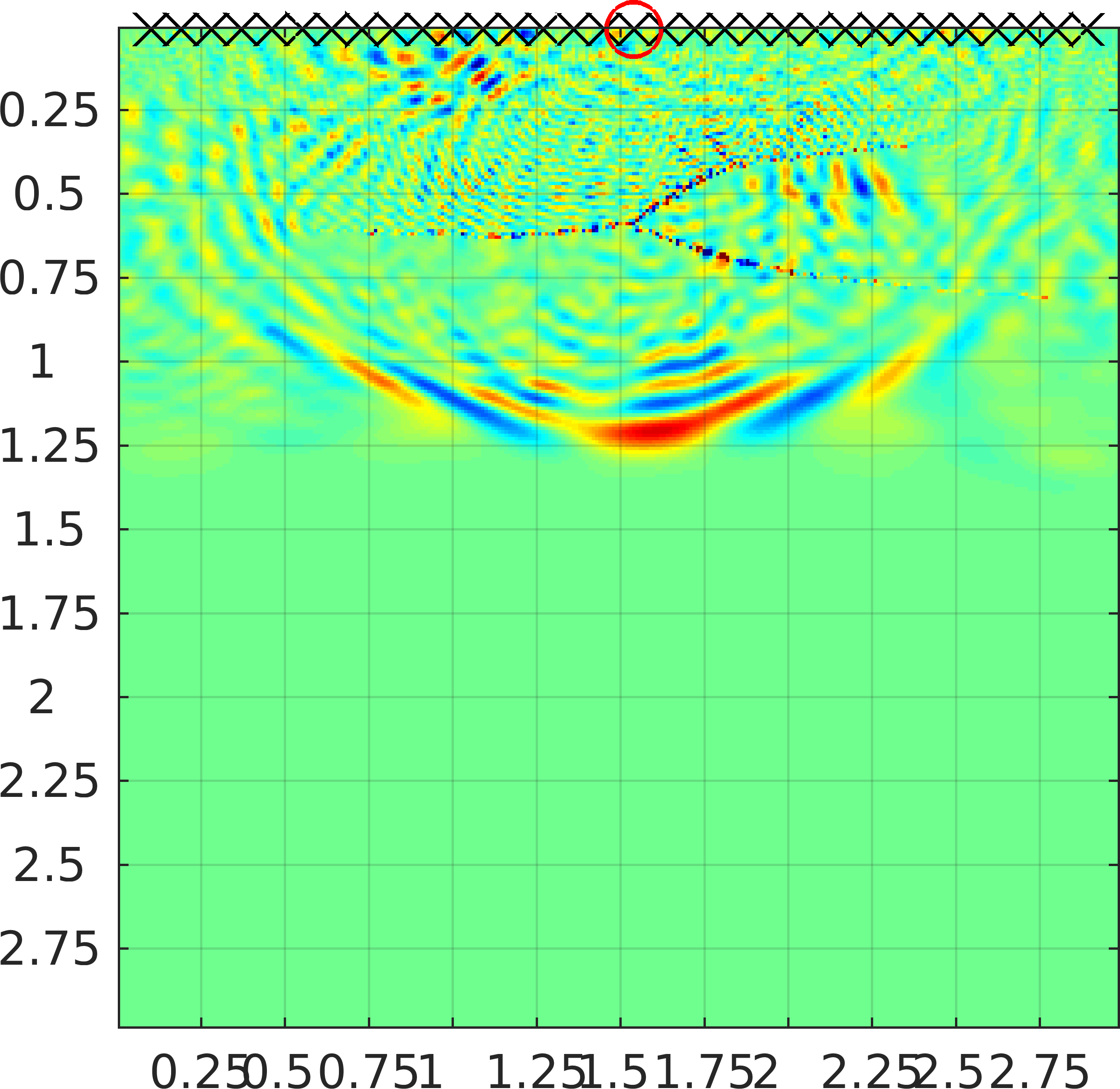} \\
\includegraphics[width=0.225\textwidth]{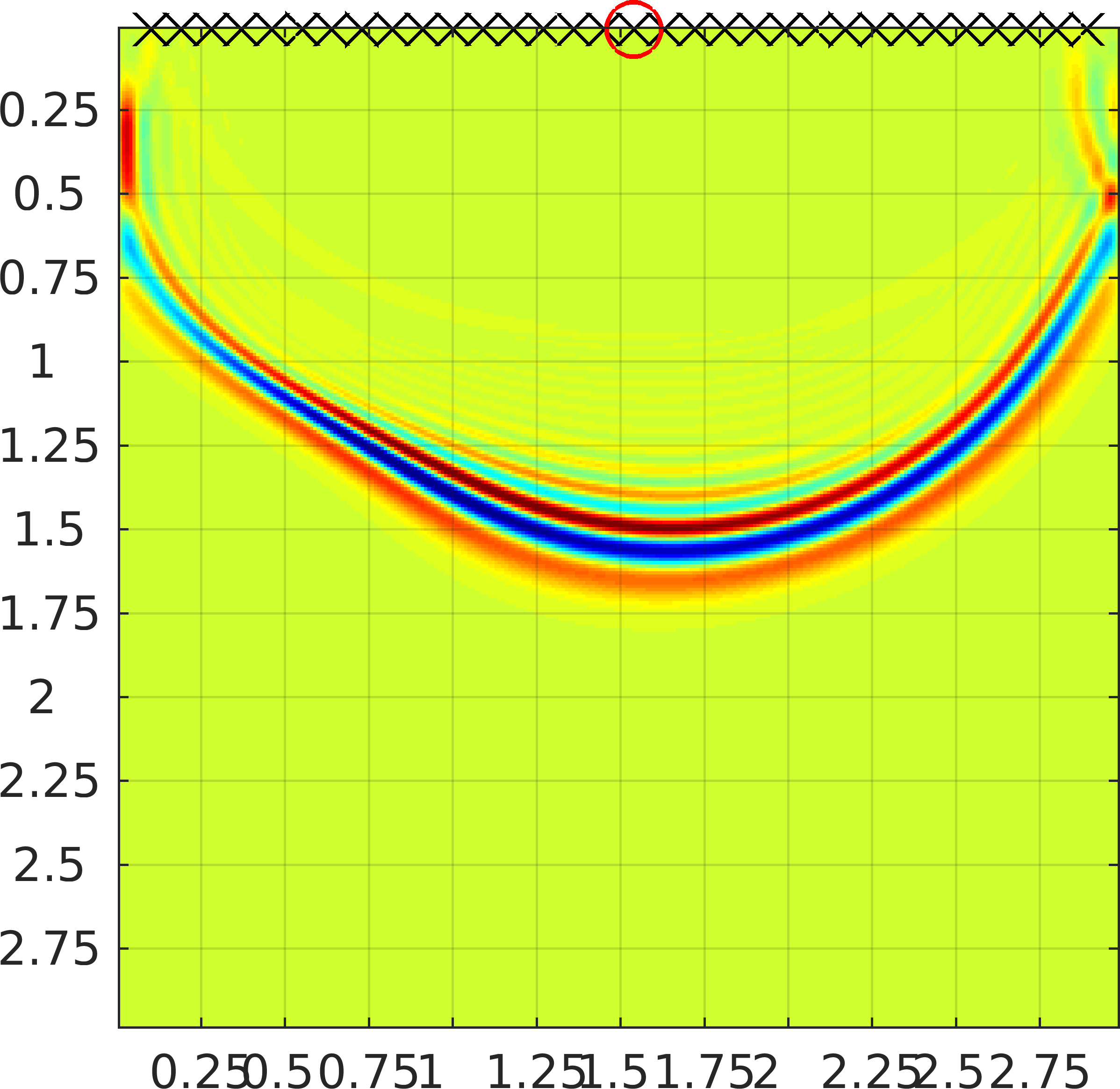} & 
\includegraphics[width=0.225\textwidth]{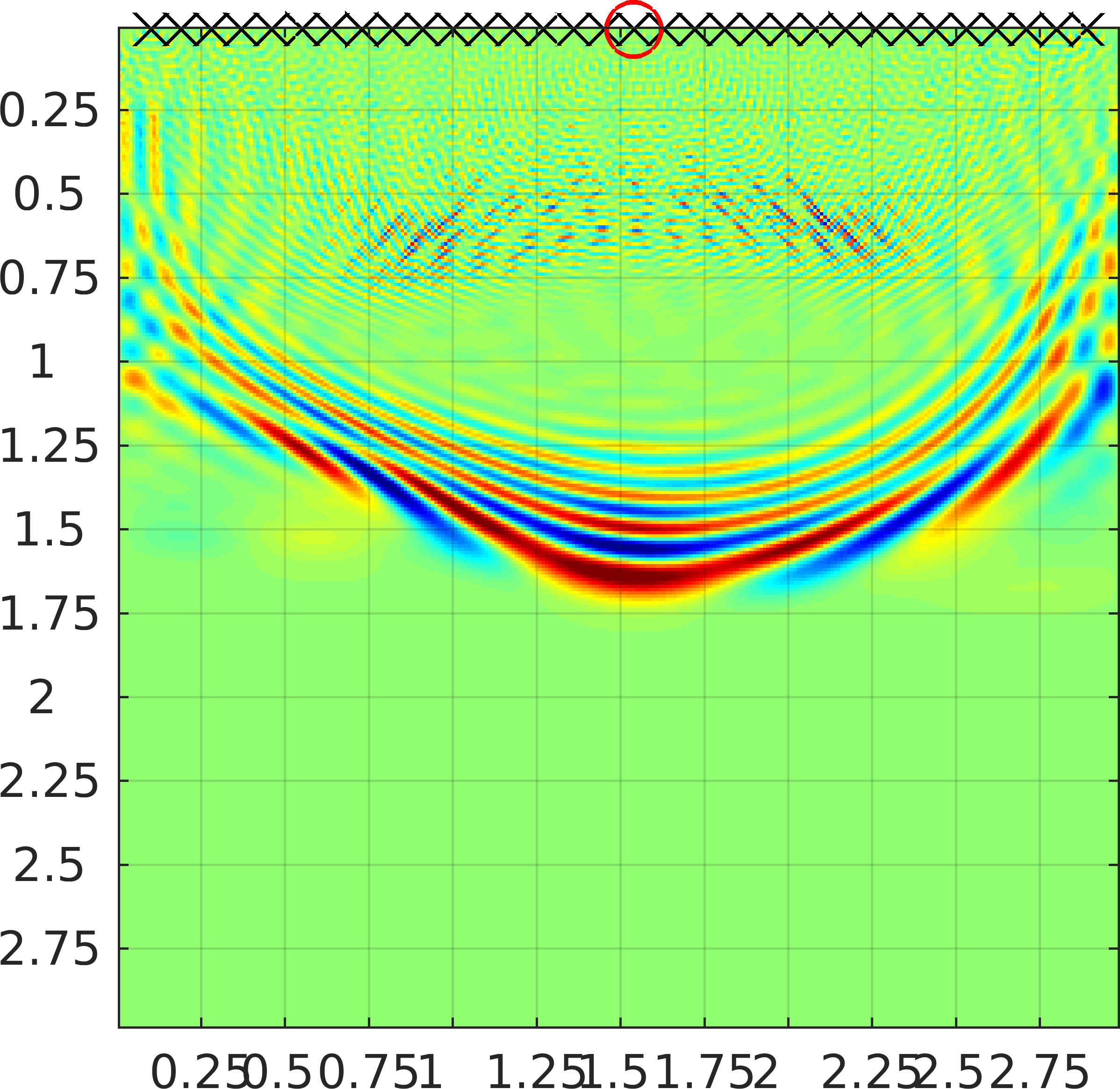} &
\includegraphics[width=0.225\textwidth]{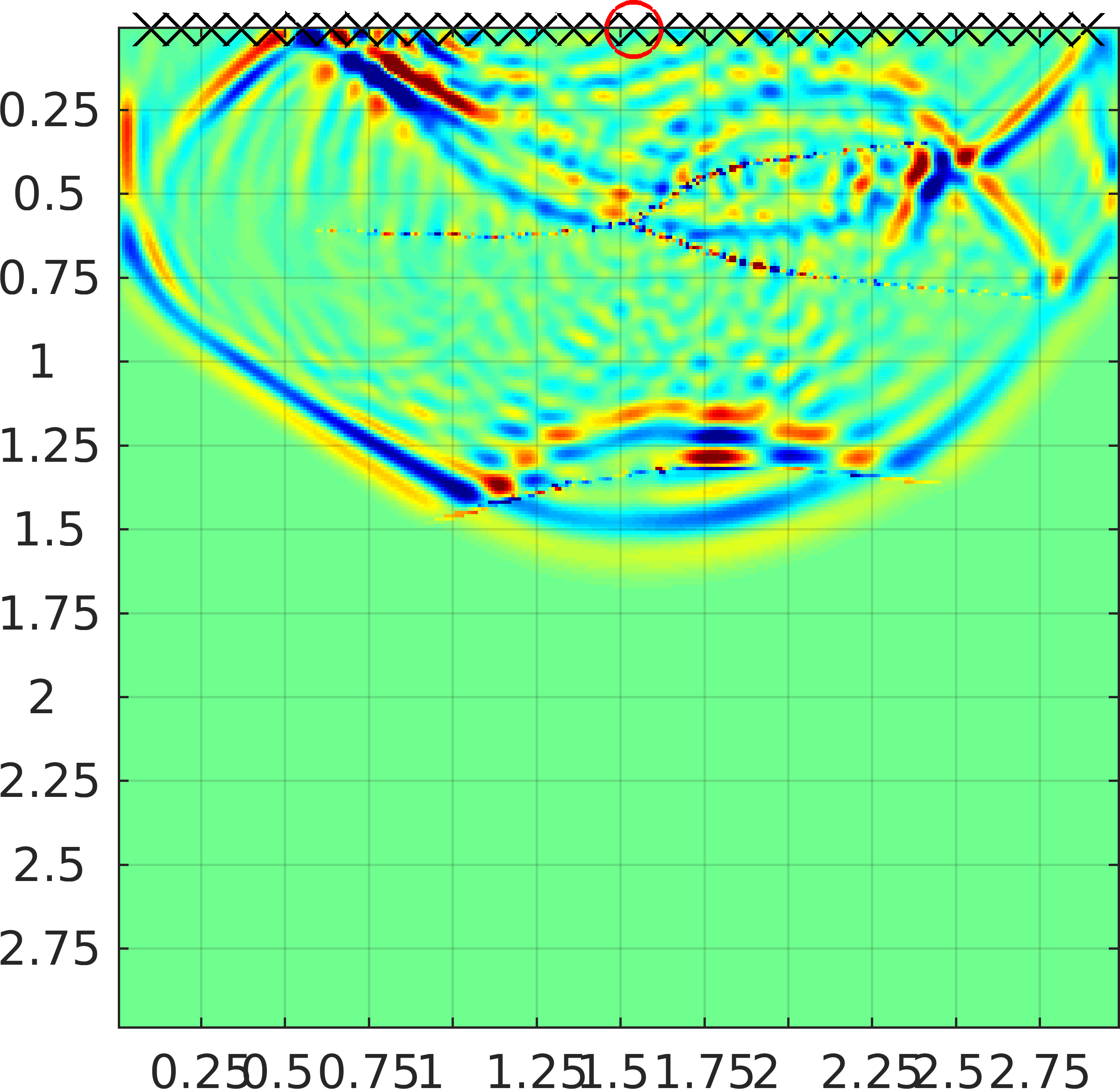} & 
\includegraphics[width=0.225\textwidth]{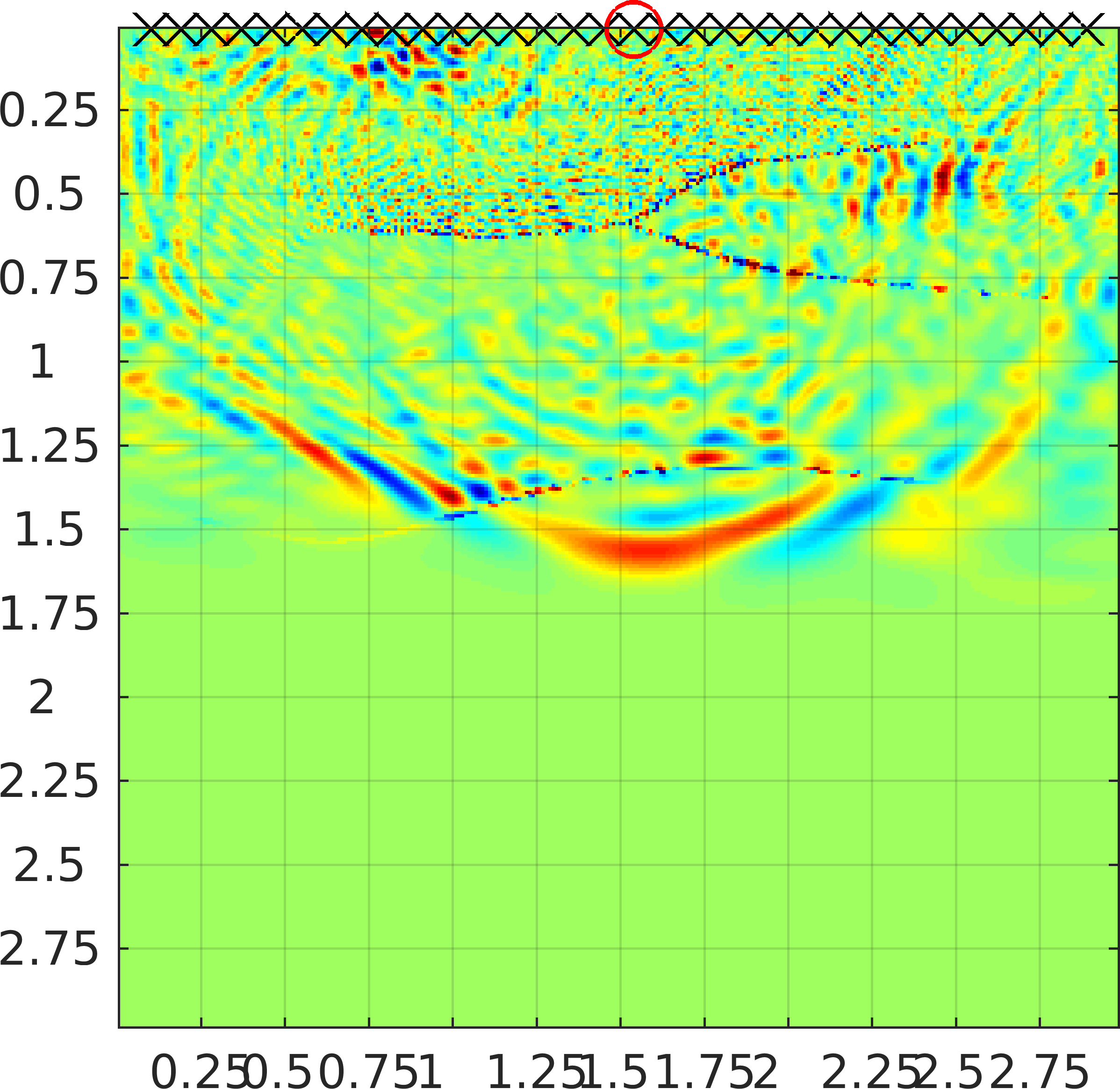} 
\end{tabular}
\caption{Symmetrized wavefield snapshots and their orthogonalized counterparts. 
Snapshots $\buh^k_{(o)}$, $\buh^k$ correspond to the kinematic model $c_o(x)$ and the true 
medium $c(x)$ respectively and their orthogonalized counterparts are $\bv_{(o)}^k$, $\bv^k$.
For every snapshot a single component is plotted corresponding to the source $j=17$ (red $\circ$) 
at the middle of the array of $m=32$ transducers (black $\times$). 
Rows (top to bottom) correspond to discrete times $t_k = k \tau$, $k = 15, 24, 30, 40, 50$.}
\label{fig:UV}
\end{figure}

At the heart of our imaging approach is the implicit orthogonalization of the symmetrized wavefield 
snapshots $\buh^k$, $k = 0, 1, \ldots, n-1$. This is an essentially nonlinear operation (due to the 
nonlinearity of block Cholesky factorization in (\ref{eqn:blockcholromalg}) and Cholesky factor inversion
in (\ref{eqn:projromalg})) that makes possible accounting for such nonlinear phenomena as multiple
reflections, as discussed later. Wavefield snapshots for the true medium $c(x)$ and the kinematic model 
$c_o(x)$ as well as their orthogonalized counterparts are shown in Figure \ref{fig:UV}. The following 
features can be observed.

First, when comparing orthogonalized snapshots $\bv^k$ to $\buh^k$ (and also $\bv^k_{(o)}$ to 
$\buh^k_{(o)}$) we observe the \emph{focusing} in both the range and cross-range. Hereafter
we refer to the direction orthogonal to the transducer array as the range and to the direction along
the array as the range. The orthogonalized snapshots have a well pronounced peak below the 
corresponding source (up to the lateral bending of the wavefront due to $c_{o}(x)$ variations). 
This peak becomes wider as the propagation time increases, i.e. we observe cross-range defocusing 
which is the consequence of the decrease of effective aperture as the wave travels away from the 
transducer array. 

Second, a phenomenon closely related to focusing that can only be observed in the reflecting medium 
$c(x)$ is the \emph{suppression of reflections} in the orthogonalized snapshots $\bv^k$. While the 
energy is spread out over reflections in $\buh^k$, after orthogonalization the wavefield is concentrated 
around the peak of $\bv^k$. This is an essential feature of the approach, since it allows to suppress 
the artifacts caused by multiple reflections that conventional linear imaging methods are unable to
handle. No matter how complex the reflections in the true medium $c(x)$ are, the propagator $P$ is 
probed in (\ref{eqn:romprop}) with focused orthonormal basis functions $\bv^k$ that are localized 
around the same locations as $\bv_{(o)}^k$. Note that the presence of reflections can actually 
improve the focusing, as we can see comparing $\bv^k$ to $\bv^k_{(o)}$ in Figure \ref{fig:UV}.

\subsection{Approximations of delta functions}
\label{subsec:approxdelta}

\begin{figure}
\begin{tabular}{cccc}
$z_1 = 500\;m$ & $z_2 = 800\;m$ & $z_3 = 1000\;m$ & $z_4 = 1330\;m$ \\
\includegraphics[width=0.225\textwidth]{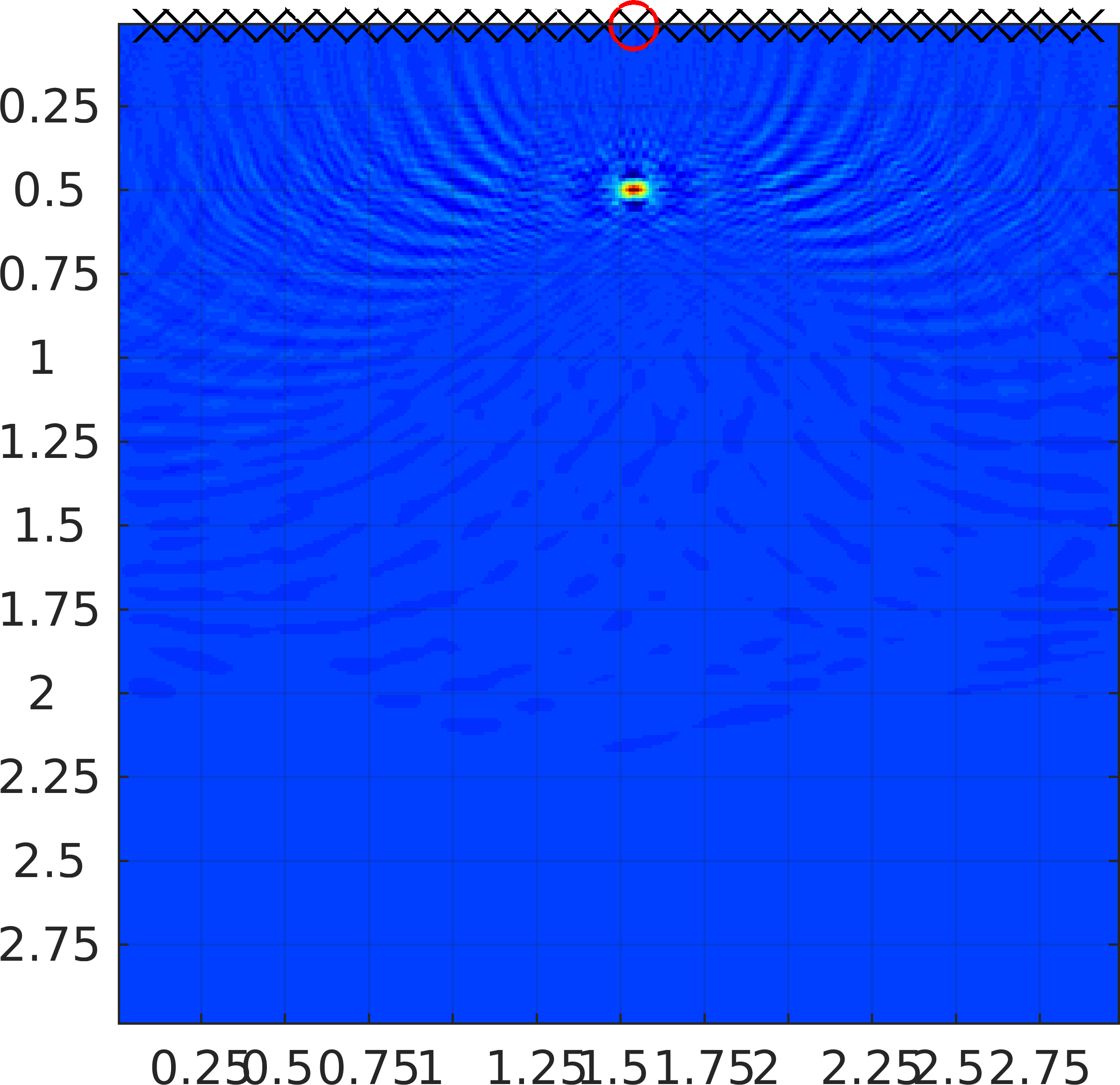} & 
\includegraphics[width=0.225\textwidth]{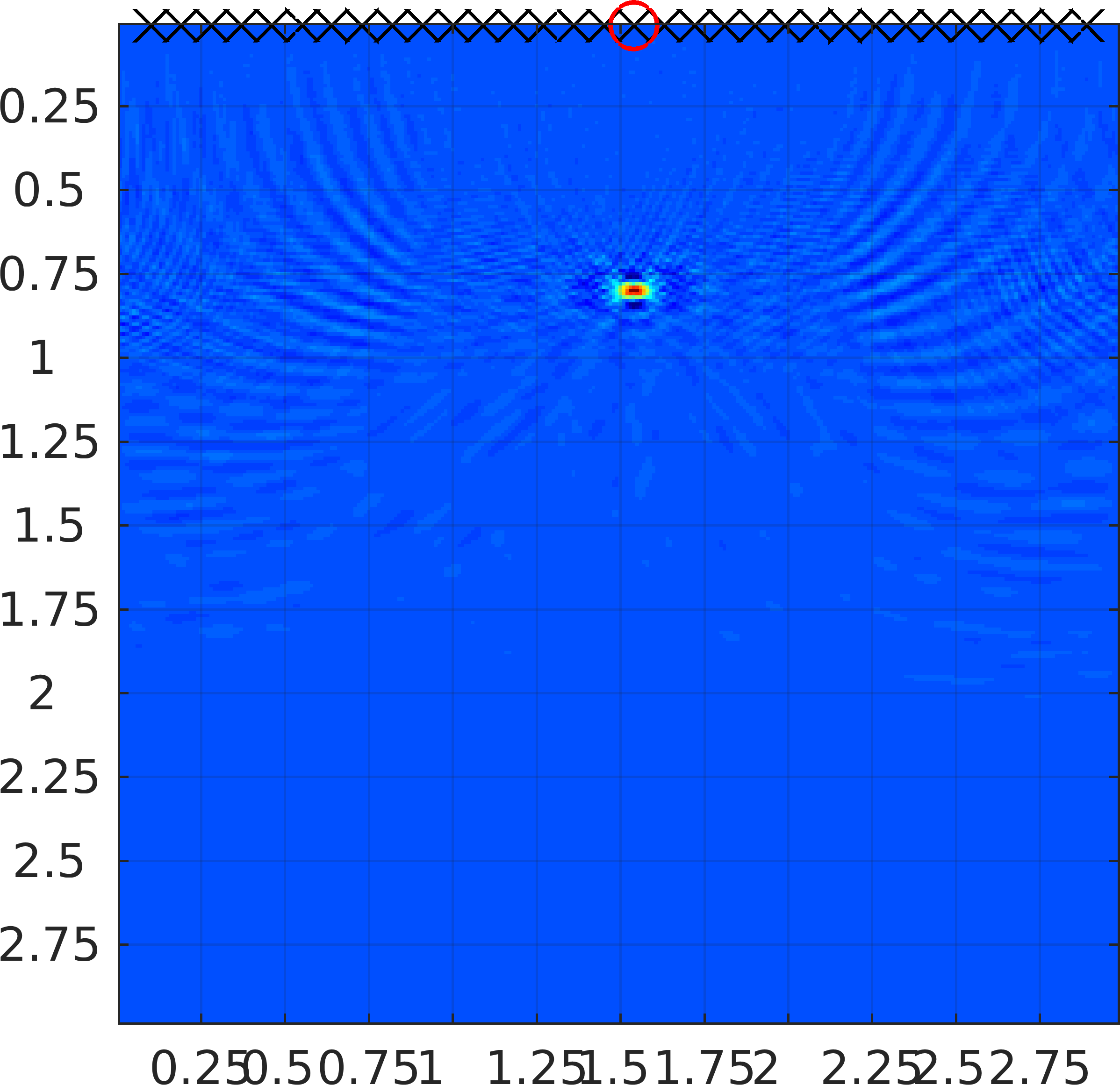} & 
\includegraphics[width=0.225\textwidth]{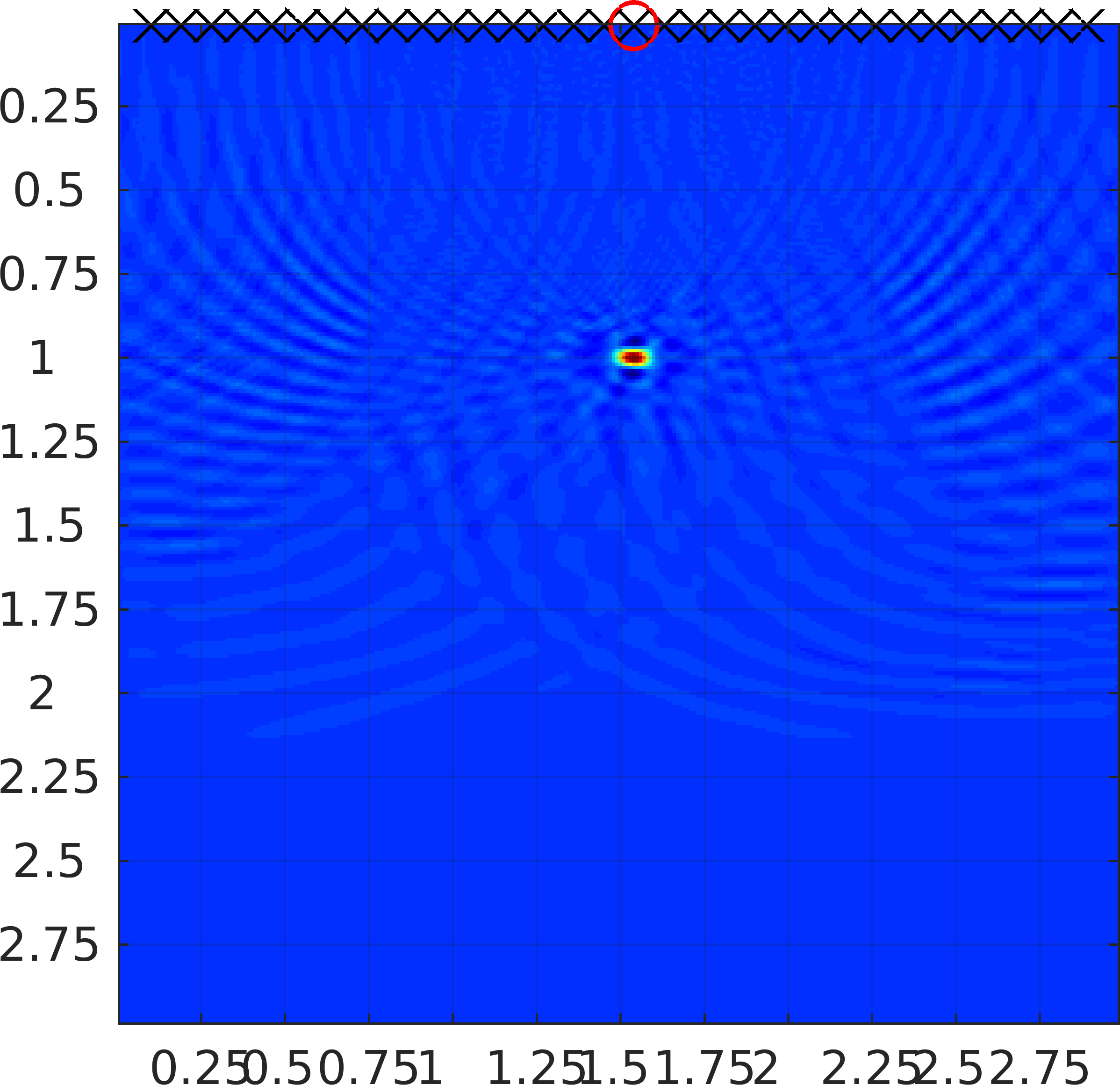} & 
\includegraphics[width=0.225\textwidth]{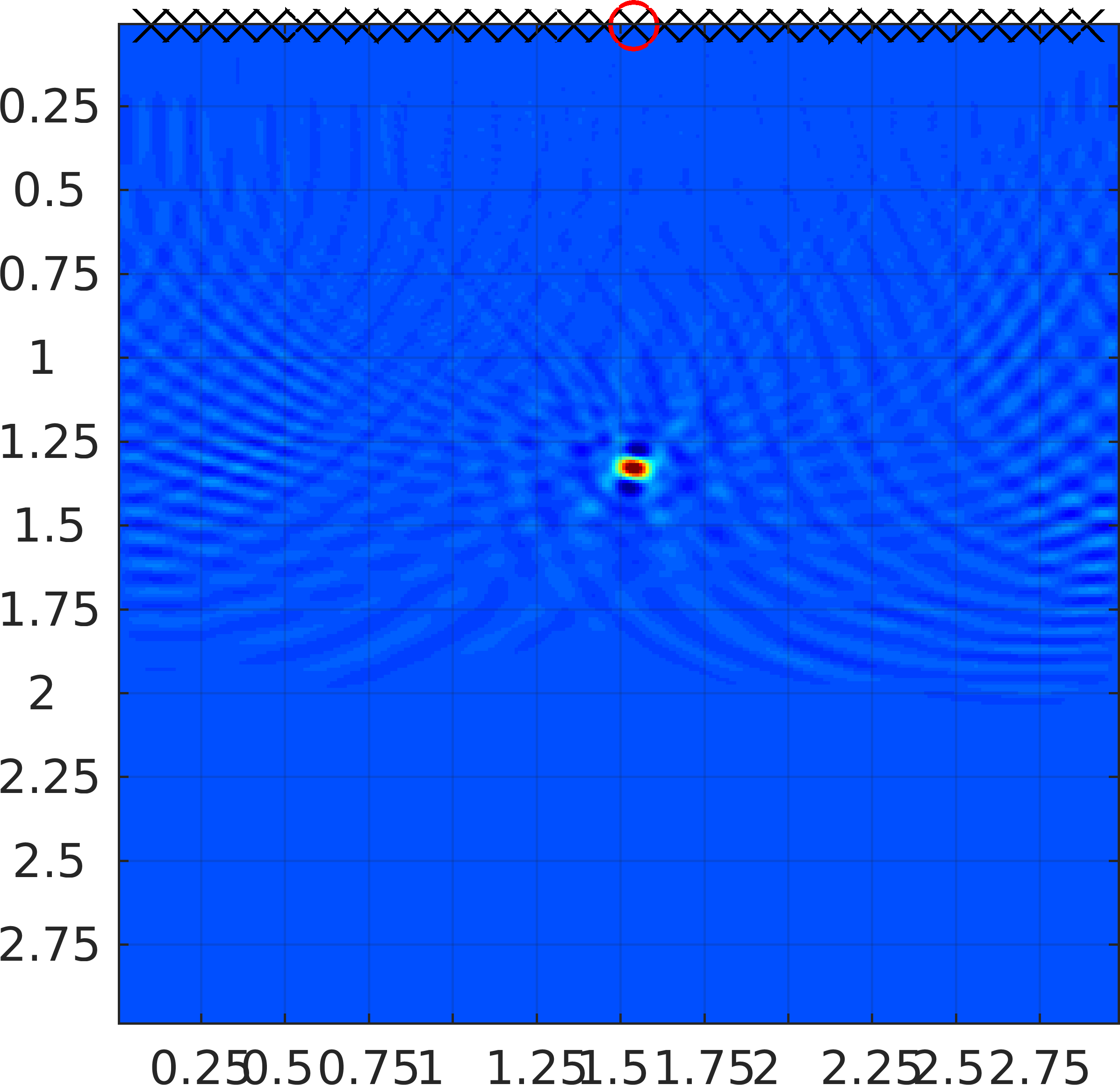} \\
\includegraphics[width=0.225\textwidth]{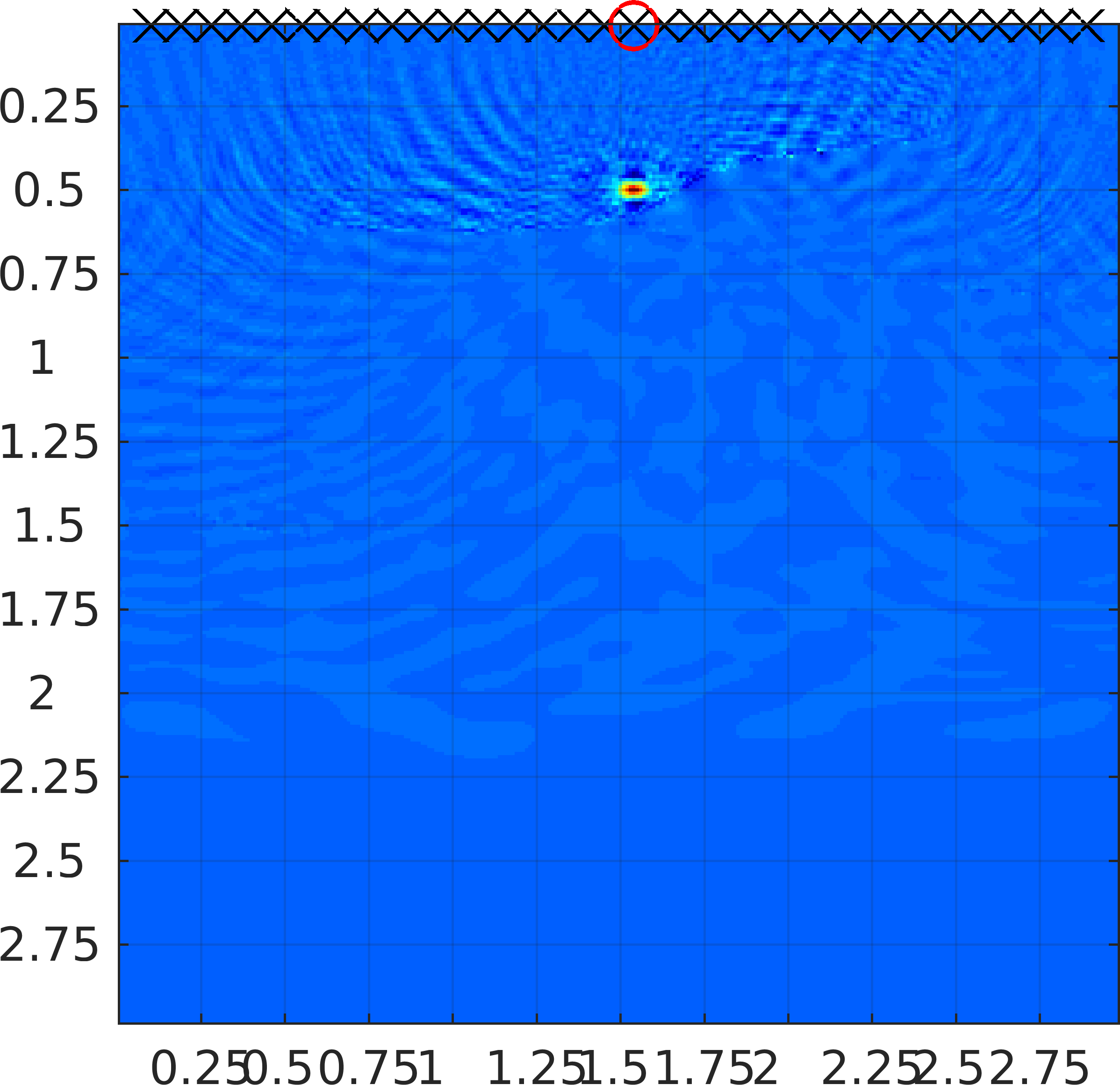} & 
\includegraphics[width=0.225\textwidth]{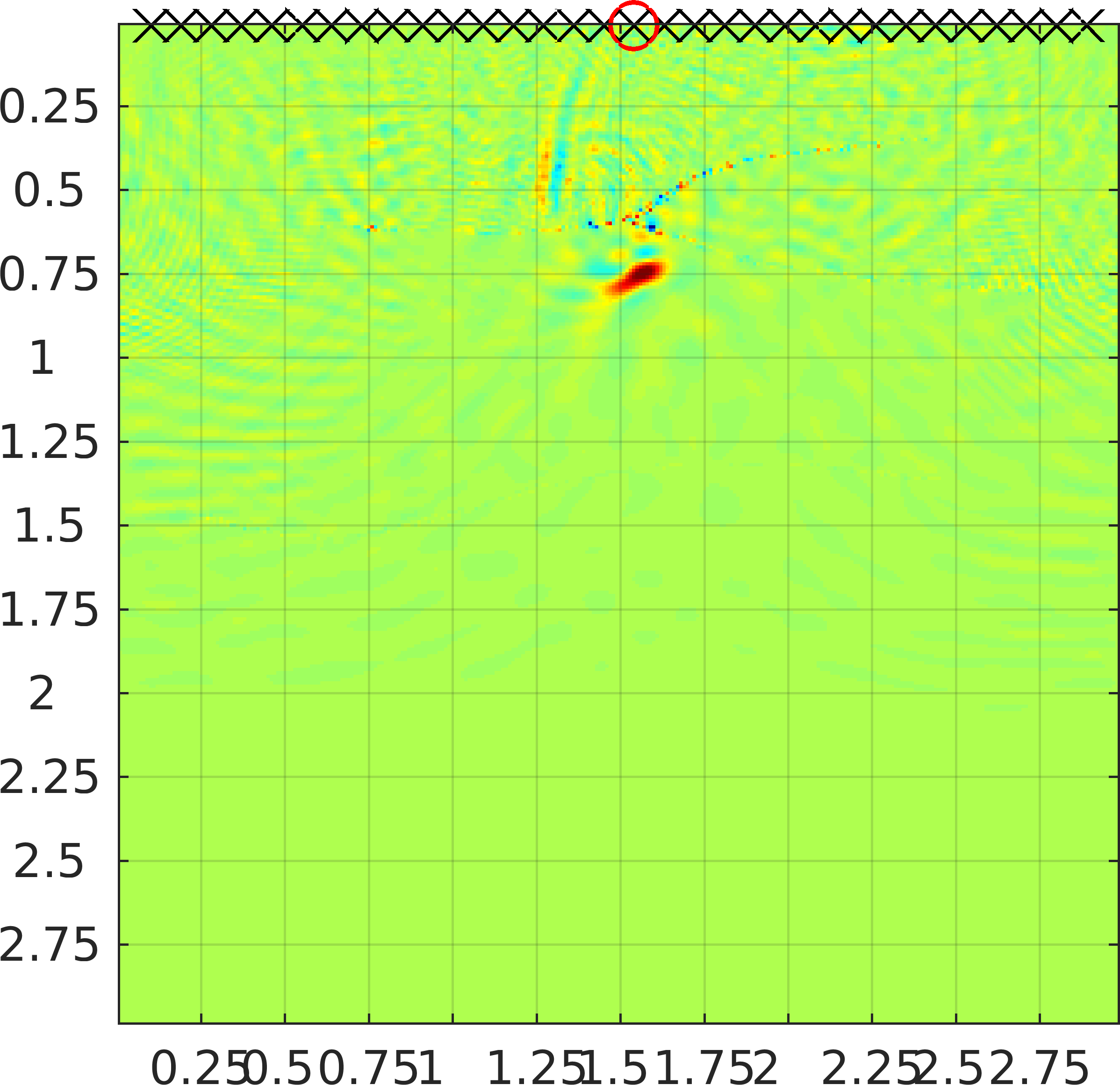} & 
\includegraphics[width=0.225\textwidth]{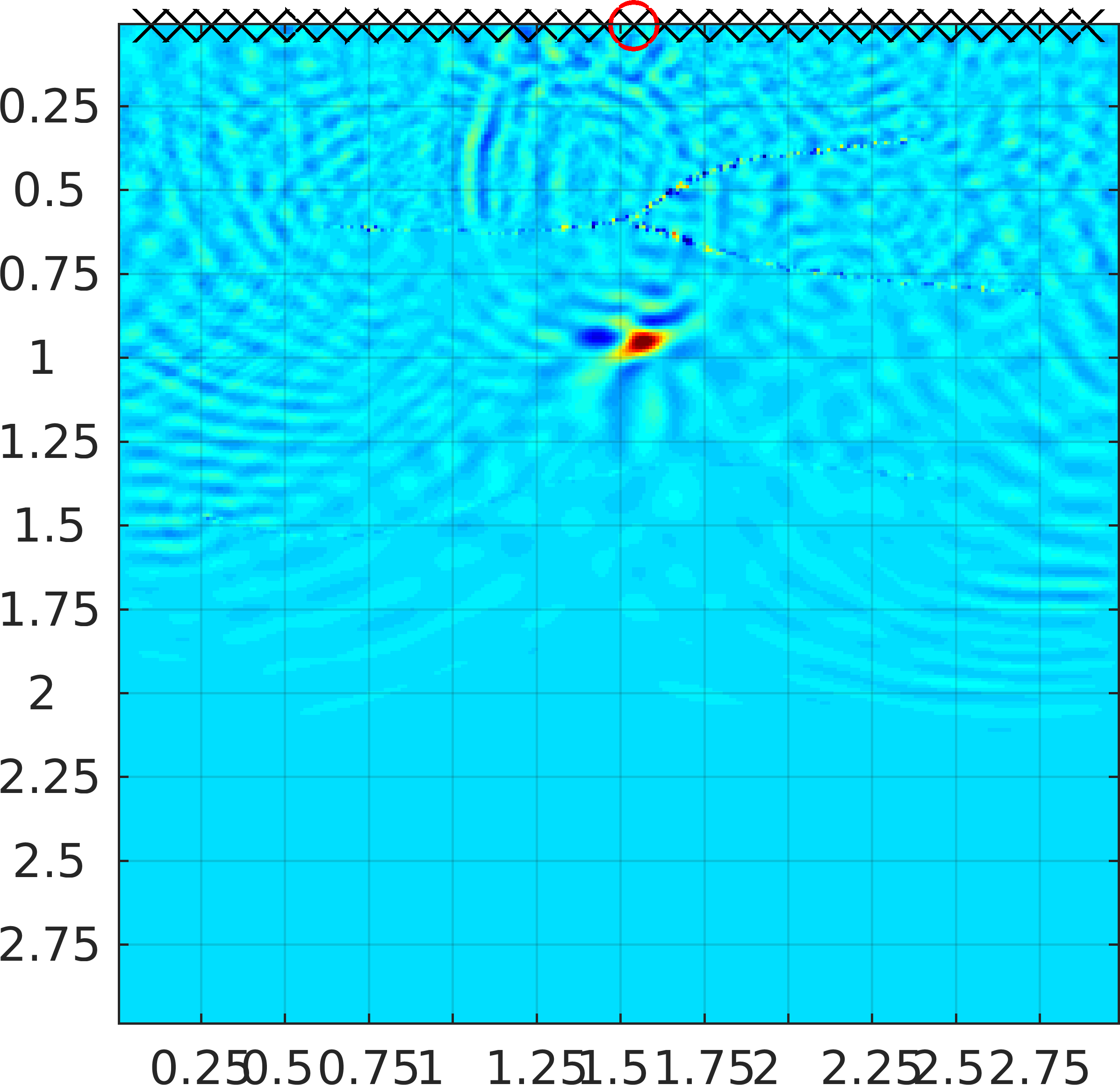} & 
\includegraphics[width=0.225\textwidth]{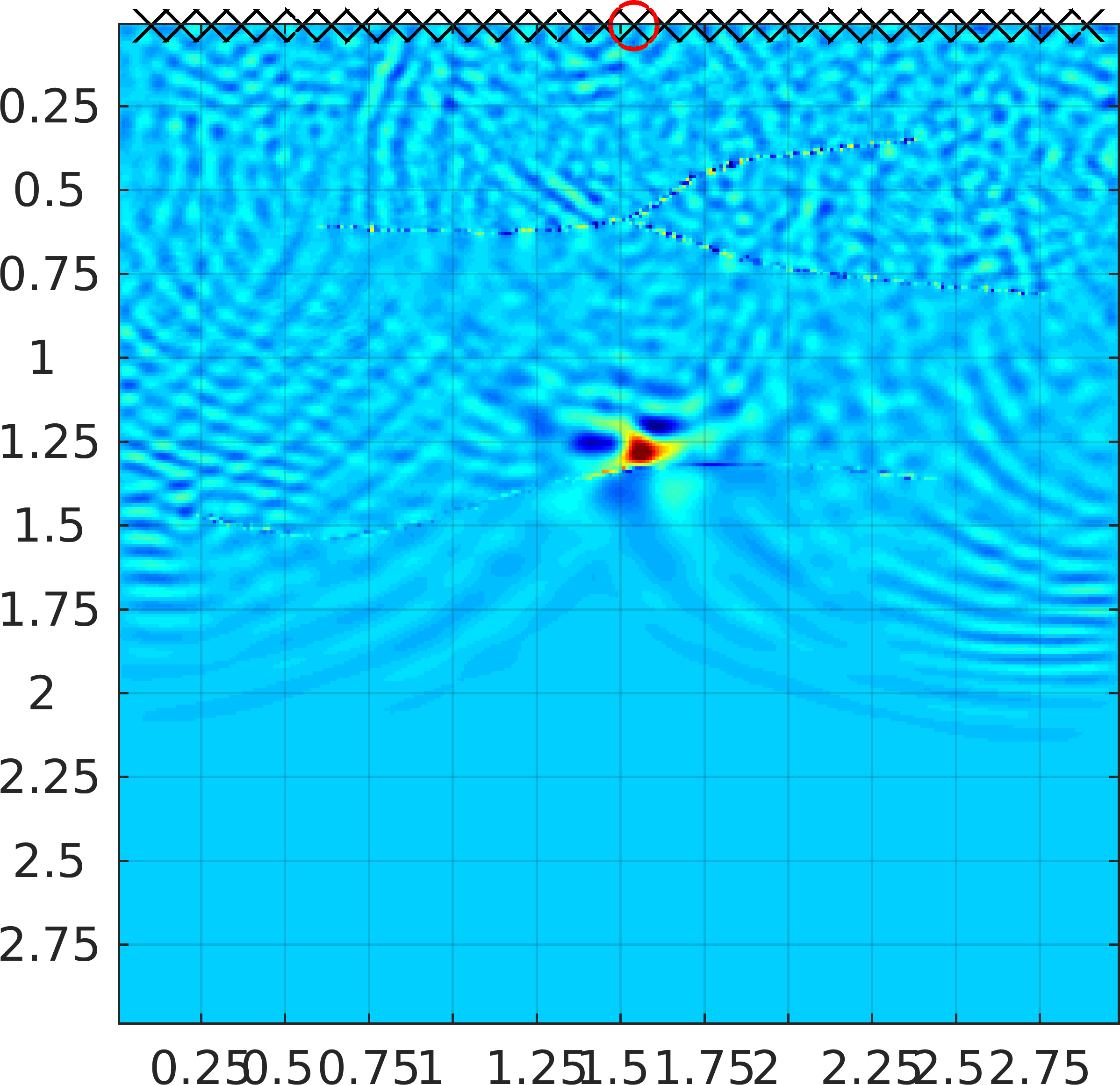}
\end{tabular}
\caption{Approximations of $\delta$-functions. 
Top row: $[\bV_{(o)} \bV^*_{(o)}](x_\alpha, y)$;
bottom row $[\bV_{(o)} \bV^*](x_\alpha, y)$, 
for $y \in \Omega$ and $x_\alpha = (1540\;m, z_\alpha)$, $\alpha=1,2,3,4$.}
\label{fig:delta}
\end{figure}

The imaging functional (\ref{eqn:imgfun}) can be written as
\be
\cI_{BP}(x) = [\bV_{(o)} \bV^*] (x, \;\cdot\;) \left[ P [\bV \bV^*_{(o)}] (\;\cdot\;, x) \right] - 
[\bV_{(o)} \bV^*_{(o)}] (x, \;\cdot\;) \left[ P_{(o)} [\bV_{(o)} \bV^*_{(o)}] (\;\cdot\;, x) \right],
\ee
which according to (\ref{eqn:greendotprod}) and (\ref{eqn:imgapprox}) should approximate
\be
g(x,x,\tau) - g_{(o)}(x,x,\tau) = \left< \delta_x, P \delta_x \right> - \left< \delta_x, P_{(o)} \delta_x \right>.
\ee
Thus, in order for the imaging method to work, the following approximations of delta functions
must hold
\be
[\bV_{(o)} \bV^*](x, y) = [\bV \bV^*_{(o)}](y, x) \approx \delta_{x}(y), \quad
[\bV_{(o)} \bV^*_{(o)}](x, y)  \approx \delta_{x}(y).
\label{eqn:vvtapproxdelta}
\ee

Given the focused nature of orthogonalized snapshots $\bv^k$ observed in section 
\ref{subsec:orthsnap}, we expect that their ``outer products'' in (\ref{eqn:vvtapproxdelta})
will be well localized, as expected from good approximations of delta functions. This conjecture is 
explored in Figure \ref{fig:delta}, where we plot $[\bV_{(o)} \bV^*_{(o)}](x_\alpha, y)$ and
$[\bV_{(o)} \bV^*](x_\alpha, y)$ for a few fixed points $x_\alpha \in \Omega$ as functions of
$y \in \Omega$. We observe indeed that both $[\bV_{(o)} \bV^*_{(o)}](x_\alpha, y)$ and
$[\bV_{(o)} \bV^*](x_\alpha, y)$ are well localized and have a clearly pronounced peak for
$y \approx x_\alpha$ while being close to zero elsewhere. Moreover, we note that the focusing 
of the outer products $\bV_{(o)} \bV^*_{(o)}$ and $\bV_{(o)} \bV^*$ is much tighter than that 
of the orthogonalized snapshots $\bv^k$. 

The outer products can be seen as nonlinear analogues of point spread functions quantifying the 
resolution of the imaging functional $\cI_{BP}$. The nonlinearity comes in as dependency of 
$\bV$ in $\bV_{(o)} \bV^*$ on the true medium $c(x)$. Note though that the dependency is
kinematic in nature, that is if $c_{(o)}(x)$ captures the kinematics of the true medium $c(x)$,
we expect
\begin{equation}
[\bV_{(o)} \bV^*](x, y) \approx [\bV_{(o)} \bV^*_{(o)}](x, y),
\end{equation}
as observed in Figure \ref{fig:delta}.

\subsection{Images and comparison to reverse time migration}
\label{subsec:imgcomprtm}

\begin{figure}
\centering
\begin{tabular}{cc}
$\cI_{BP}$ & $\cI_{RTM}$ \\
\includegraphics[width=0.44\textwidth]{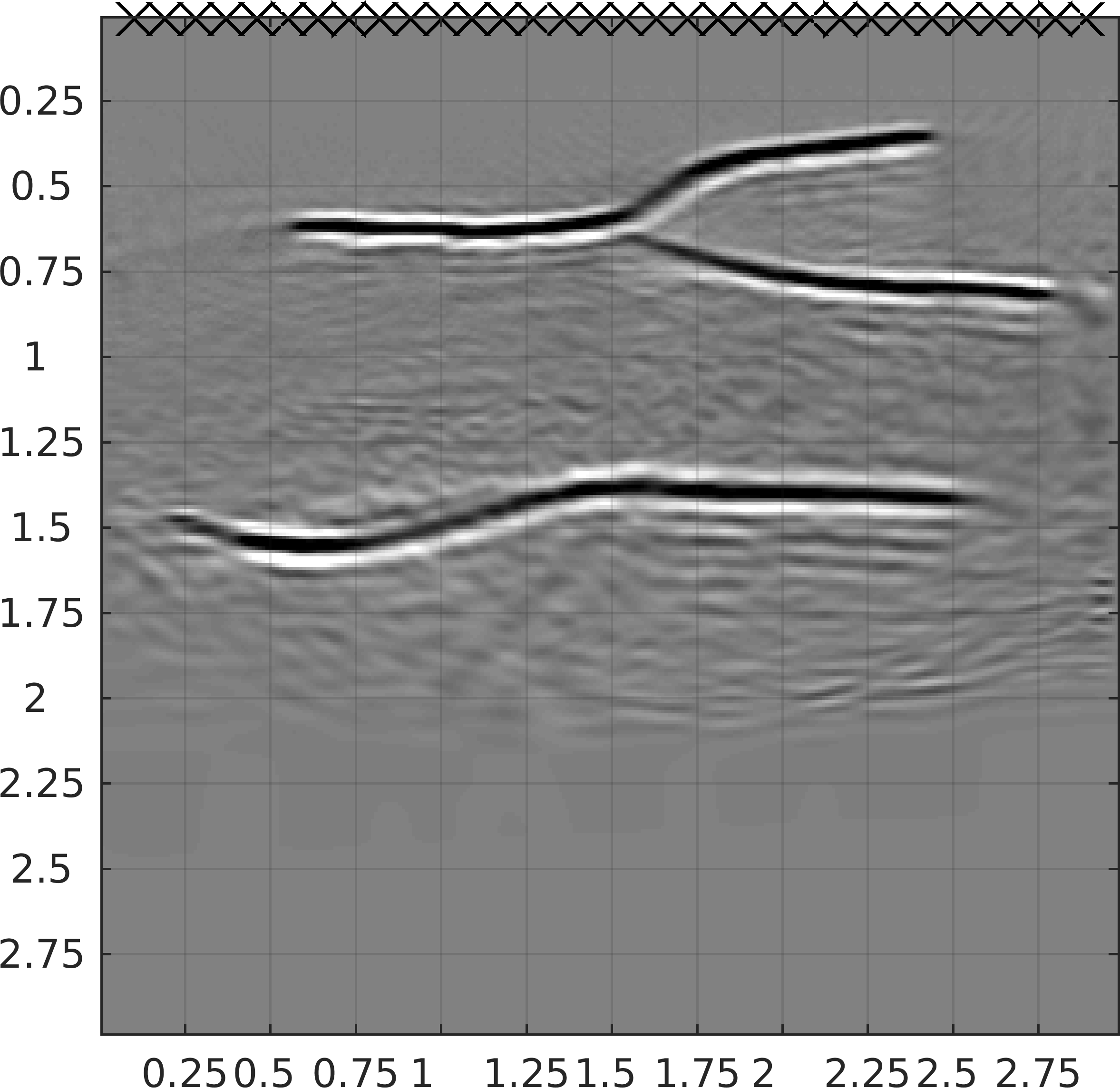} & 
\includegraphics[width=0.44\textwidth]{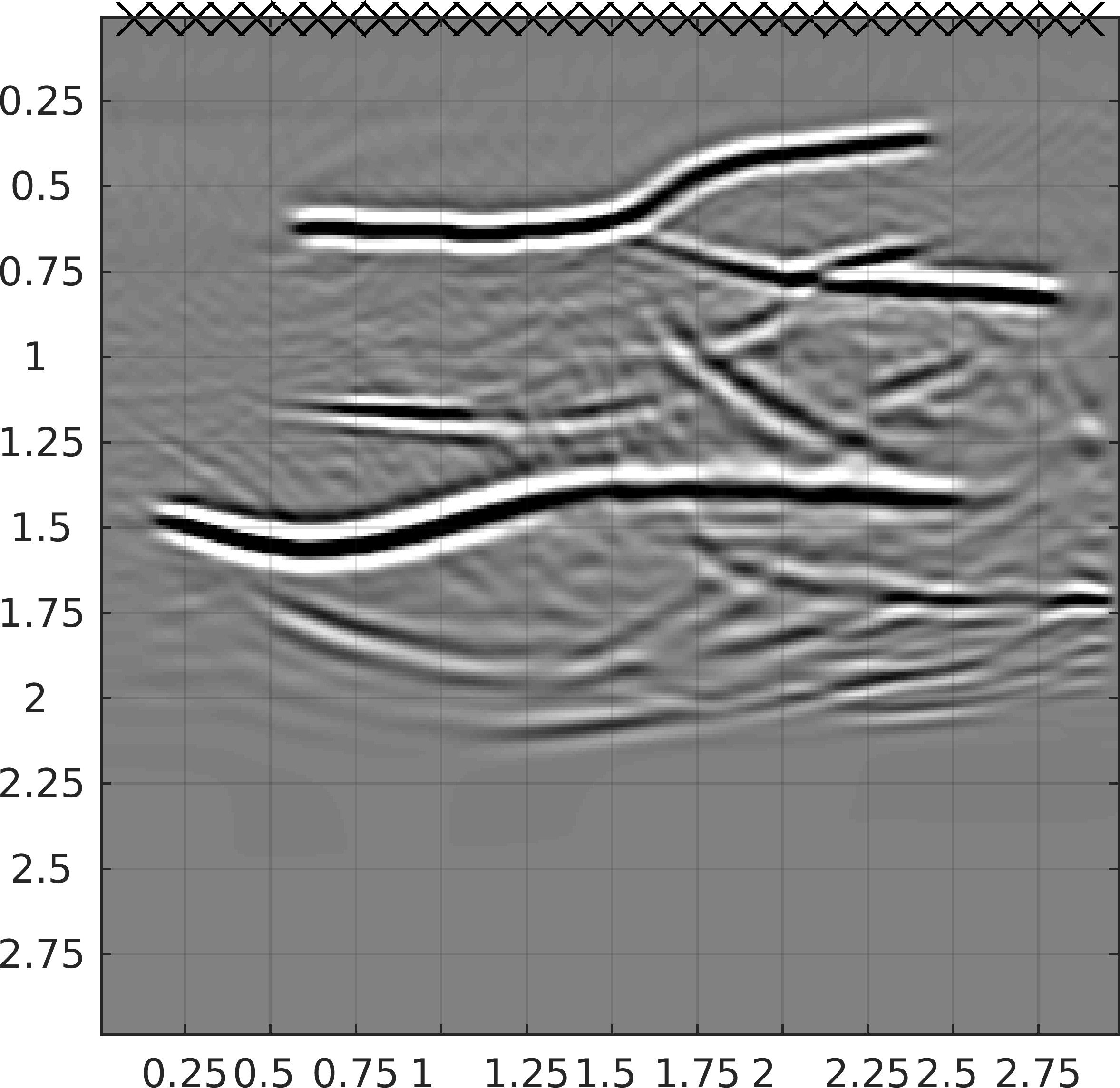} \\
\includegraphics[width=0.44\textwidth]{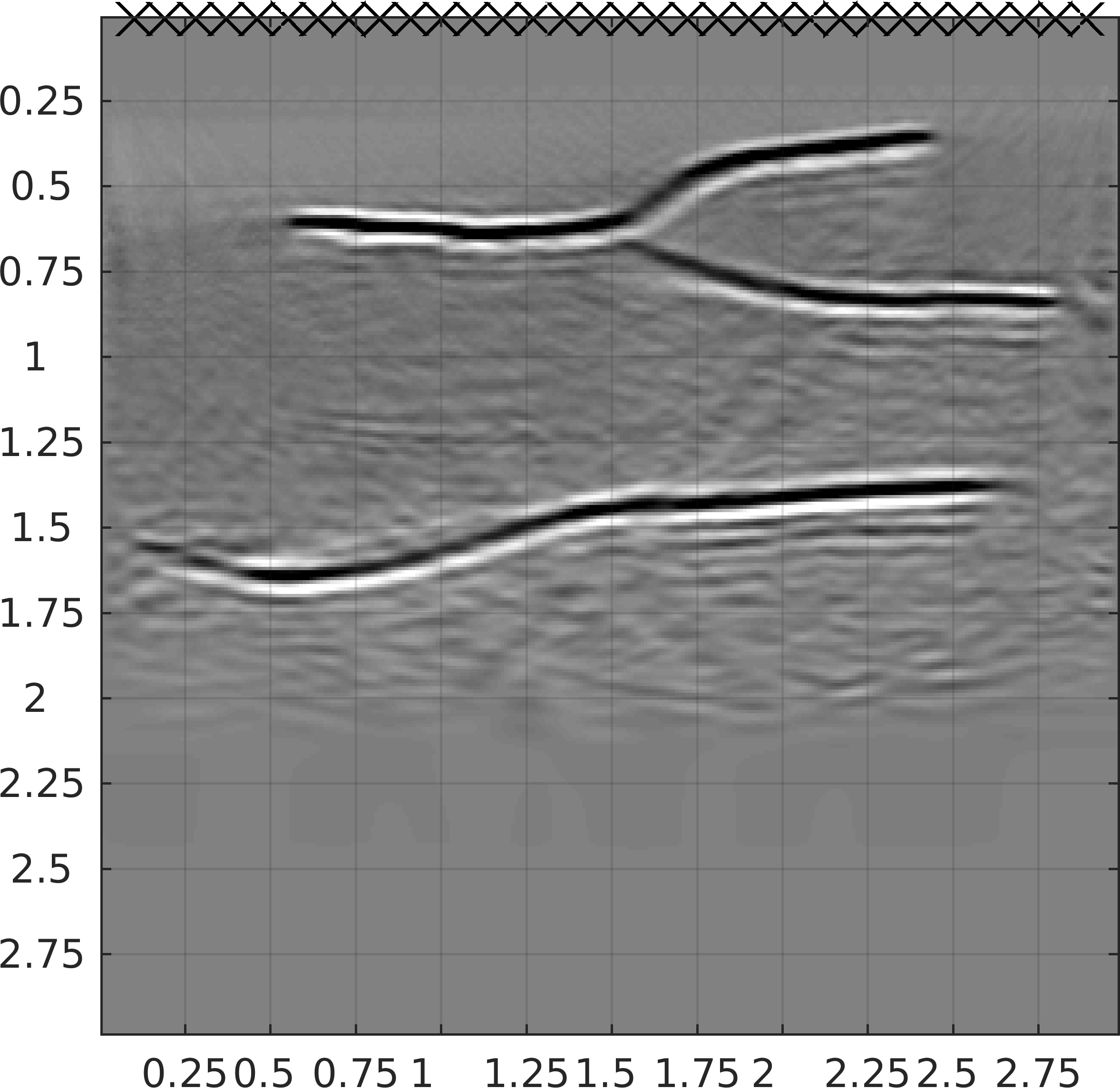} & 
\includegraphics[width=0.44\textwidth]{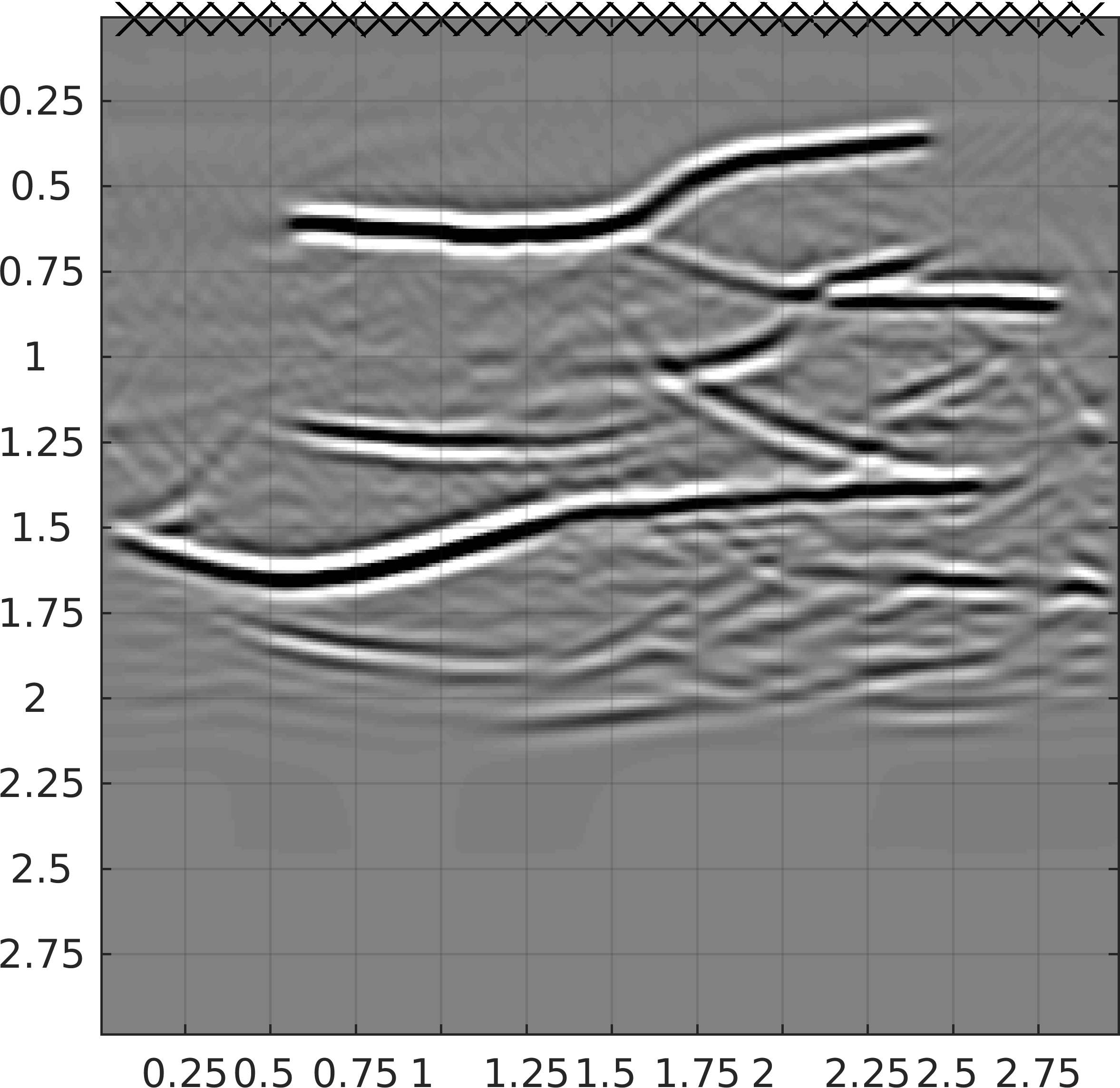}
\end{tabular}
\caption{Image comparison between backprojection and RTM functionals $\cI_{BP}$ (left column)
and $\cI_{RTM}$ (right column) respectively. Top row: kinematic model $c_{(o)}(x)$ from
Figure \ref{fig:cfrac4}; bottom row: constant kinematic model $c_{(o)} \equiv 2.5 \; km/s$.
Locations of transducers are black $\times$. All distances in $km$. 
}
\label{fig:imgfrac4}
\end{figure}

Once the focusing and localization properties of orthogonalized snapshots $\bv^k$ and their outer
product are established, we can finally assess the performance of the backprojection imaging 
functional $\cI_{BP}$ on the model described in section \ref{subsec:nummodel} 
(see Figure \ref{fig:cfrac4}). We compare our approach with reverse time migration (RTM)
\cite{baysal1983reverse}, a linear imaging approach that is the standard method for high quality
imaging with waves in geophysics \cite{chang1987elastic, symes2007reverse, yoon2004challenges}.
Here we use a particular version of RTM known as \emph{pre-stack} migration, which is more accurate 
than the computationally cheaper \emph{post-stack} migration. The pre-stack RTM imaging 
functional is hereafter referred to as $\cI_{RTM}$.

Note that both $\cI_{BP}$ and $\cI_{RTM}$ imaging functionals produce weaker images of 
reflectors that are further away from the array. This is particularly pronounced when the wave has 
to pass through reflectors closer to the array before reaching the far ones, thus losing energy in 
the process, as it gets reflected back to the array at each interface. To account for this behavior,
in what follows we modify both imaging functionals by scaling them with a multiplier that depends 
linearly on the distance from the imaging point to the transducer array. In a slight abuse of 
notation we still refer to such modified functionals as $\cI_{BP}$ and $\cI_{RTM}$.

In Figure \ref{fig:imgfrac4} we compare the images obtained with ROM backprojection formula
(\ref{eqn:imgfun}) and with conventional pre-stack RTM for the model described in section 
\ref{subsec:nummodel} with two kinematic models: a kinematic model $c_{(o)}(x)$ from 
Figure \ref{fig:cfrac4} consisting of the true smooth background without the reflectors; 
and a constant kinematic model with $c_{(o)} \equiv 2.5 \; km/s$. We observe that in both cases
our approach is clearly superior compared to its conventional linear counterpart. The main 
advantage of backprojection imaging is the automatic removal of \emph{multiple reflection}
artifacts. These artifacts appear inevitably in any linear imaging method due to the interaction
of reflectors with reflective boundaries (here mainly with the top part $\cB_A$ of the boundary) and
with each other. Each such interaction produces additional events in the data that linear migration
methods image as extra reflectors that are not actually present in the true medium. We observe 
in Figure \ref{fig:imgfrac4} that the RTM image has many such artifacts with magnitudes comparable
to the actual reflectors. Moreover, the artifacts change their location and magnitudes depending on 
the kinematic model used. Meanwhile, ROM backprojection image is entirely free of such artifacts
with the only features present being the actual reflectors.

We also observe in Figure \ref{fig:imgfrac4} that reflector images produced by $\cI_{BP}$ are thinner 
than those in RTM images. We conjecture that ROM backprojection has higher range resolution than 
RTM, but the detailed resolution analysis of ROM backprojection is outside the scope of this paper
and remains a topic of future research.

\subsection{Noisy data and regularization}
\label{subsec:noisereg}

\begin{figure}
\centering
\includegraphics[width=0.48\textwidth]{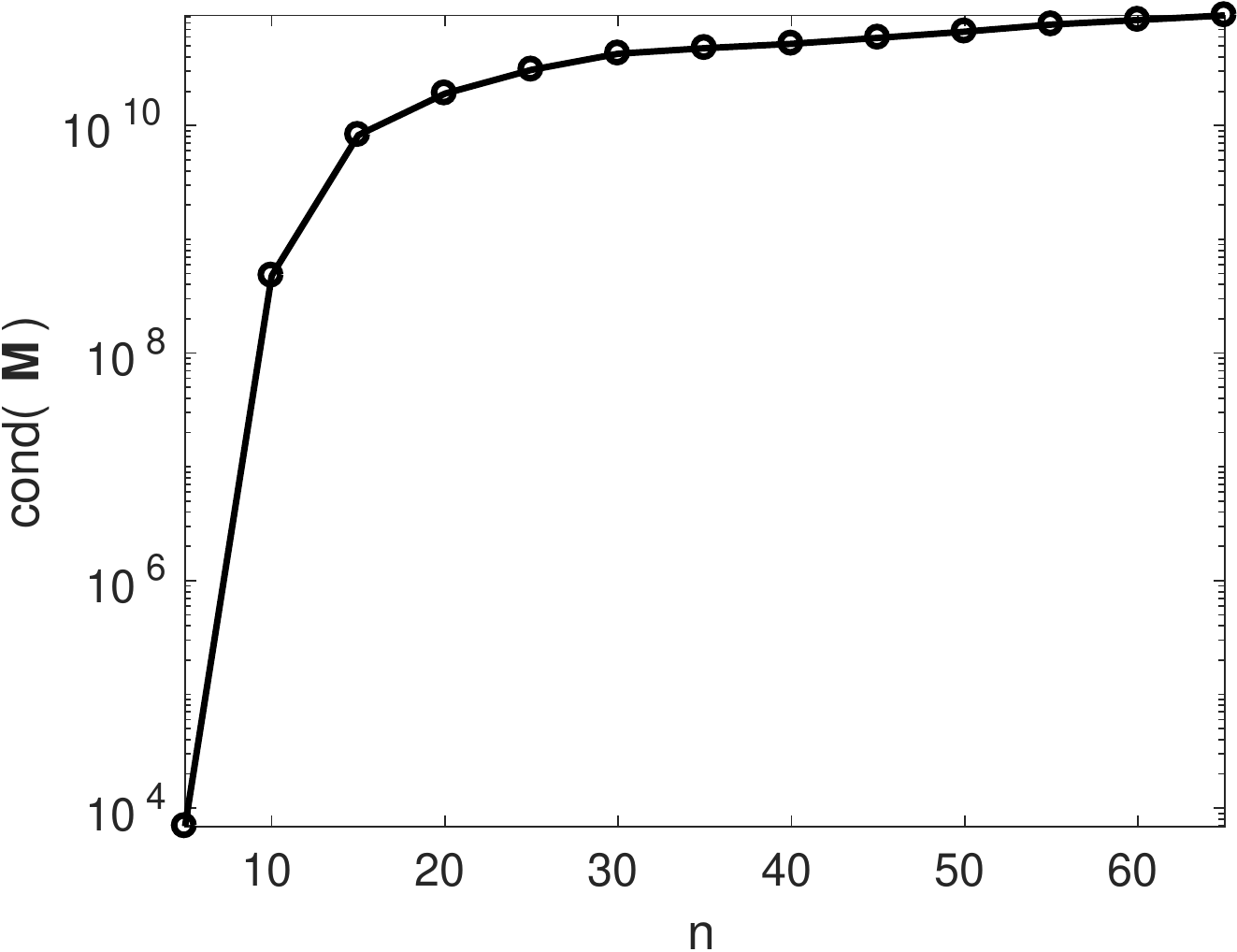}
\caption{Condition number of $\bM$ for the numerical example from section 
\ref{subsec:nummodel} for the values of $n$ from $5$ to $65$.}
\label{fig:condfrac4}
\end{figure}

\begin{figure}
\centering
\begin{tabular}{cc}
$\cI_{BP}$ & $\cI_{BP}^{(\mu)}$ \\
\includegraphics[width=0.44\textwidth]{ddbp_bp_modelfrac4_kintrue_r100_n32} & 
\includegraphics[width=0.44\textwidth]{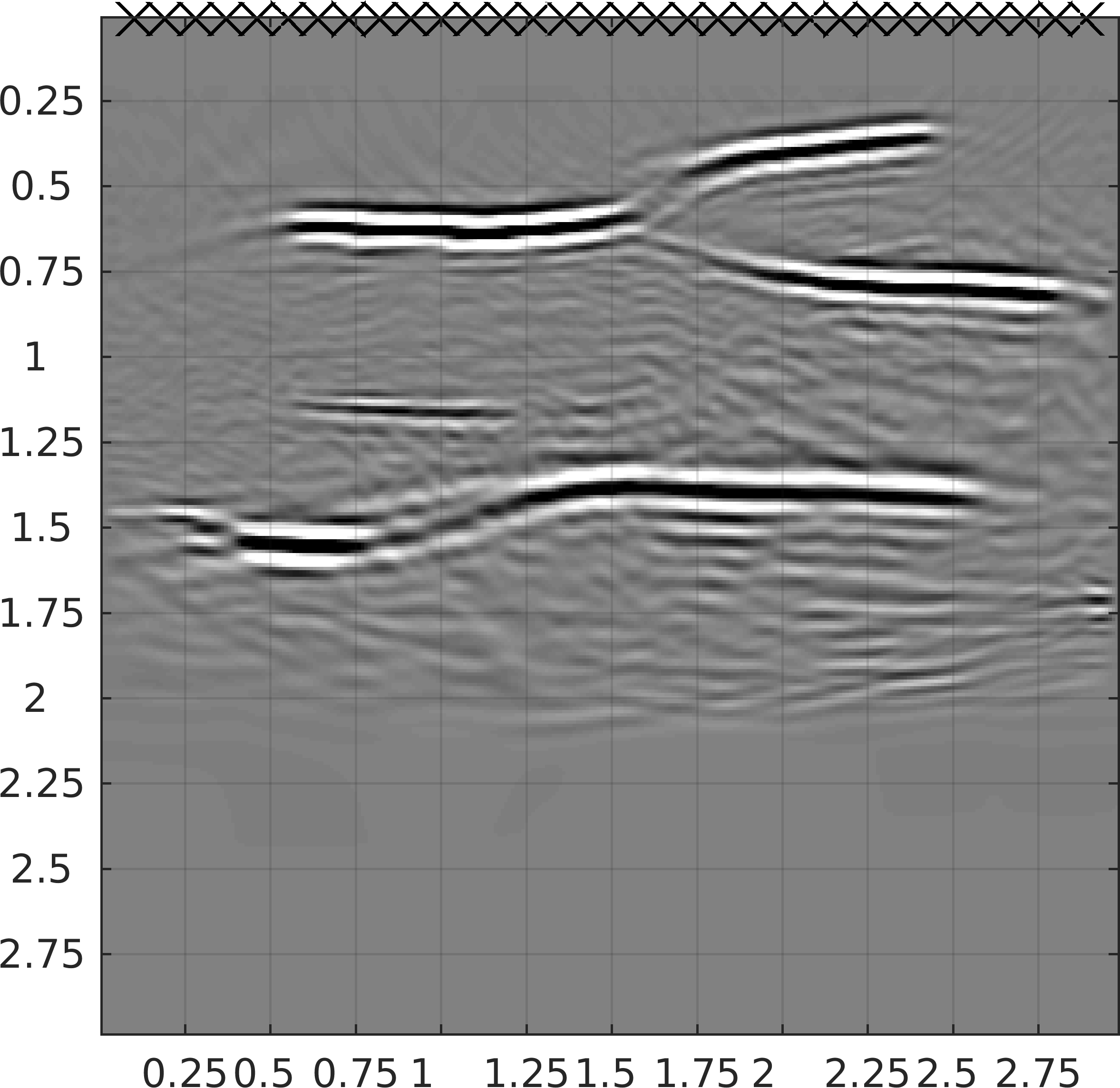} \\
\includegraphics[width=0.44\textwidth]{ddbp_bp_modelfrac4_kinconst_r100_n32} & 
\includegraphics[width=0.44\textwidth]{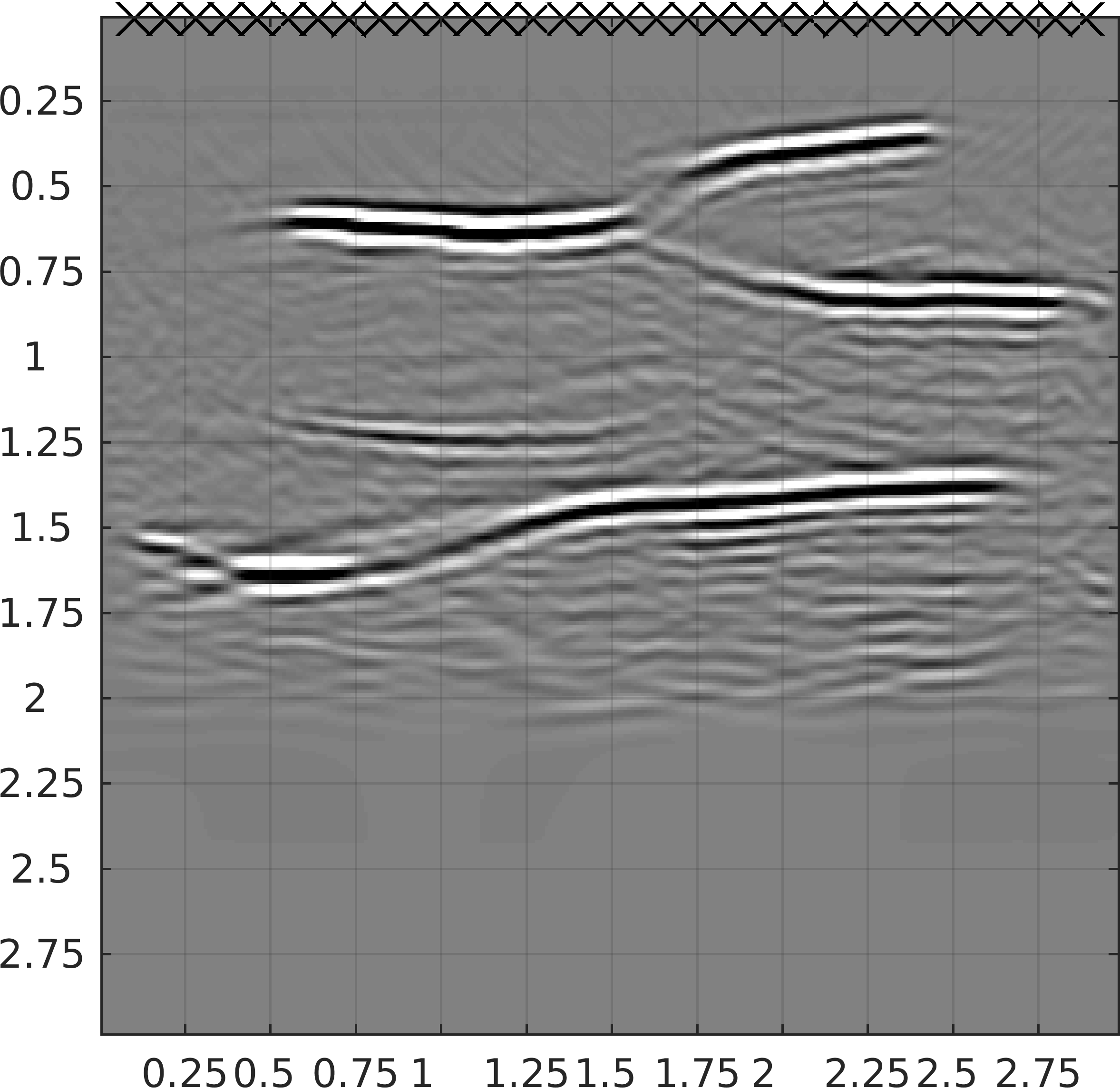} 
\end{tabular}
\caption{ROM backprojection images computed from noiseless data ($\cI_{BP}$) and from 
noisy data $\bD_{(\epsilon)}^k$ regularized with Algorithm \ref{alg:datareg} ($\cI_{BP}^{(\mu)}$).
Noise level $\epsilon = 10 \%$. Top row: kinematic model $c_{(o)}(x)$ from Figure \ref{fig:cfrac4}; 
bottom row: constant kinematic model $c_{(o)} \equiv 2.5 \; km/s$. Locations of transducers are 
black $\times$. All distances in $km$.}
\label{fig:imgfrac4reg}
\end{figure}

As mentioned above, construction of the reduced model $\bPt$ with Algorithm \ref{alg:romcompute} 
is a procedure nonlinear in the data $\bD^k$. While such approach allows for automatic removal of 
multiple reflection artifacts and improved resolution, as demonstrated in section
\ref{subsec:imgcomprtm}, it also has downsides compared to conventional linear migration. 

One particular issue that we address here is that of stability. The first nonlinear operation in 
Algorithm \ref{alg:romcompute} is block Cholesky factorization (\ref{eqn:blockcholromalg}) of the 
mass matrix $\bM$. Unfortunately, the mass matrix is ill-conditioned. In Figure \ref{fig:condfrac4} 
we display the condition number of $\bM$ for the numerical example from section \ref{subsec:nummodel}
as the half-number of data time samples $n$ increases from $5$ to $65$, the value used to compute
the images in section \ref{subsec:imgcomprtm}. The condition number increases rapidly from
$10^4$ to $10^{10}$ for $n=15$ and the subsequent increase is rather slow. 

Since $\mbox{cond}(\bM)$ is so large, the smallest eigenvalues of the mass matrix are very close to 
zero. Note that $\bM$ depends linearly on $\bD^k$ in (\ref{eqn:massmatdataalg}), and even a small 
amount of noise in the data can push the eigenvalues of the perturbed mass matrix into the negative
values. Thus, a regularization scheme is required for the imaging method to work with noisy data. 
Here we propose a simple regularization approach that is an analogue of Tikhonov regularization 
for linear inverse problems. The idea is to push the eigenvalues of the mass matrix 
$\bM_{(\epsilon)}$ computed from the noisy data $\bD_{(\epsilon)}^k$ back into the positive values.
However, this has to be done consistently with the algebraic structure of the problem. This is 
achieved by the following simple algorithm.

\begin{algorithm}[Data regularization]
\label{alg:datareg}~\\
Given noisy data $\bD_{(\epsilon)}^k$, $k=0,\ldots,2n-1$, and an initial regularization
parameter $\mu > 1$, to compute the regularized data $\bD_{\mu}^k$, $k=0,\ldots,2n-1$, 
perform the following steps:
\begin{enumerate}
\item Compute the regularized data
\be
\bD^k_{(\mu)} = \left\{
\begin{tabular}{rl}
$\mu \bD^{0}_{(\epsilon)}$, & if $k=0$ \\
$\bD^{k}_{(\epsilon)}$, & if $k > 0$
\end{tabular}
\right., \quad k=0,\ldots,2n-1.
\label{eqn:dataregalg}
\ee

\item Assemble the regularized mass matrix 
\be
\left[ \bM_{(\mu)} \right]_{k,l} = \frac{1}{2} \left( \bD^{k+l}_{(\mu)} + \bD^{|k-l|}_{(\mu)} \right),
\quad k,l = 0,\ldots,n,
\label{eqn:massregalg}
\ee
and find its smallest eigenvalue $\lambda_{min}(\bM_{(\mu)})$.

\item If $\lambda_{min}(\bM_{(\mu)}) < 0$, increase $\mu$ and go to step 1, otherwise quit.
\end{enumerate}
\end{algorithm}

We observe that the regularized data $\bD_{(\mu)}^k$ differs from the input noisy data 
$\bD^k_{(\epsilon)}$ at the initial time instant $k = 0$ only. 
Note that $\bD^0 = \bb^* \bb \in \mathbb{R}^{m \times m}$  is a positive definite matrix,
moreover, for the sources (\ref{eqn:edelta})--(\ref{eqn:qgauss}) it is diagonally dominant.
Since $\bD^0$ only enters the diagonal blocks of the mass matrix according to (\ref{eqn:massregalg}),
for sufficiently large $\mu$ the regularized mass matrix $\bM_{(\mu)}$ also becomes positive
definite. Obviously, the size of regularization parameter $\mu$ needed to achieve 
$\lambda_{min}(\bM_{(\mu)}) > 0$ depends on the amount of noise present in $\bD_{(\epsilon)}^k$.

It is important that Algorithm \ref{alg:datareg} modifies the data itself instead of just the mass matrix.
When the regularized data is fed to ROM computation Algorithm \ref{alg:romcompute}
it will not only affect the mass matrix (\ref{eqn:massmatdataalg}), but the stiffness matrix 
(\ref{eqn:stiffmatdataalg}) as well. This makes the regularization scheme more consistent with the 
algebraic structure of the problem.

We illustrate the performance of Algorithm \ref{alg:datareg} in Figure \ref{fig:imgfrac4reg}, where we
compare $\cI_{BP}$ computed from noiseless data to $\cI_{BP}^{(\mu)}$ computed from the 
noisy data that is regularized with Algorithm \ref{alg:datareg}. The noisy data $\bD^k_{(\epsilon)}$ is 
obtained from the simulated noiseless data $\bD^k$ using a multiplicative noise model from 
\cite{borcea2014model} with a noise level parameter $\epsilon$. Such model ensures that the 
signal-to-noise ratio is approximately $100\% / \epsilon$. Here we take $\epsilon = 10 \%$ which 
results in the value of regularization parameter $\mu = 1.55$, as the minimum value for which 
$\lambda_{min}(\bM_{(\mu)}) > 0$ robustly for many realizations of the noise.

We observe in Figure \ref{fig:imgfrac4reg} a certain deterioration of the image obtained from the
regularized noisy data. First, the reflectors become somewhat thicker which indicates a slight loss
of resolution. Second, $\cI_{BP}^{(\mu)}$ loses some of the contrast at the slanted parts 
of reflectors thus emphasizing the parts of reflectors that are nearly parallel to the transducer array.
However, ROM backprojection still suppresses the multiple reflections very well. While we notice
one weak artifact appearing above the bottom reflector, the image is still almost entirely artifact-free 
compared to RTM images in Figure \ref{fig:imgfrac4}.

\subsection{Large scale numerical examples}
\label{subsec:largescale}

\begin{figure}
\centering
\includegraphics[width=0.97\textwidth]{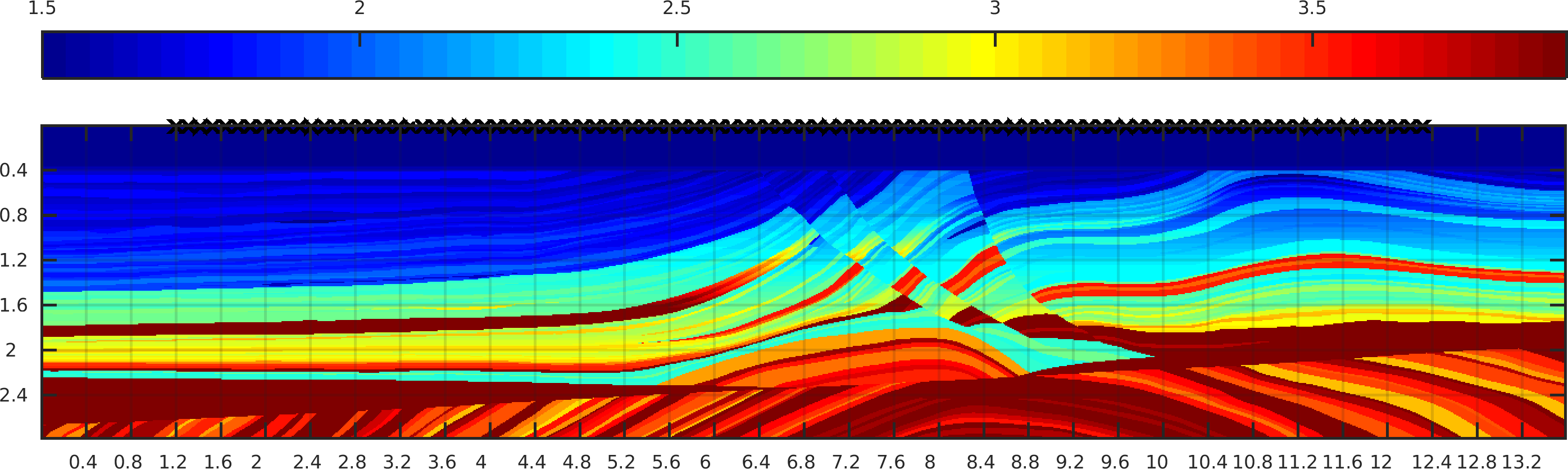}
\includegraphics[width=0.97\textwidth]{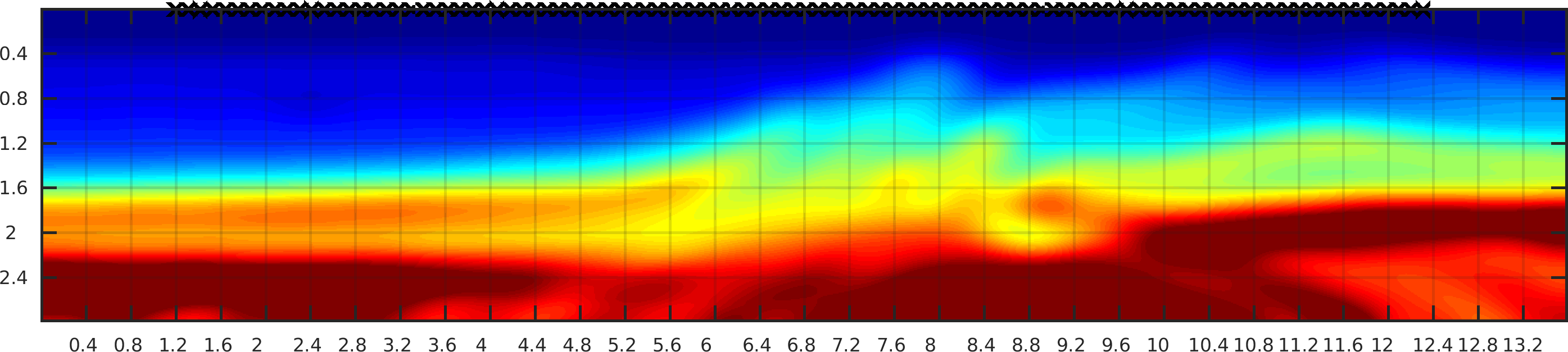}
\includegraphics[width=0.97\textwidth]{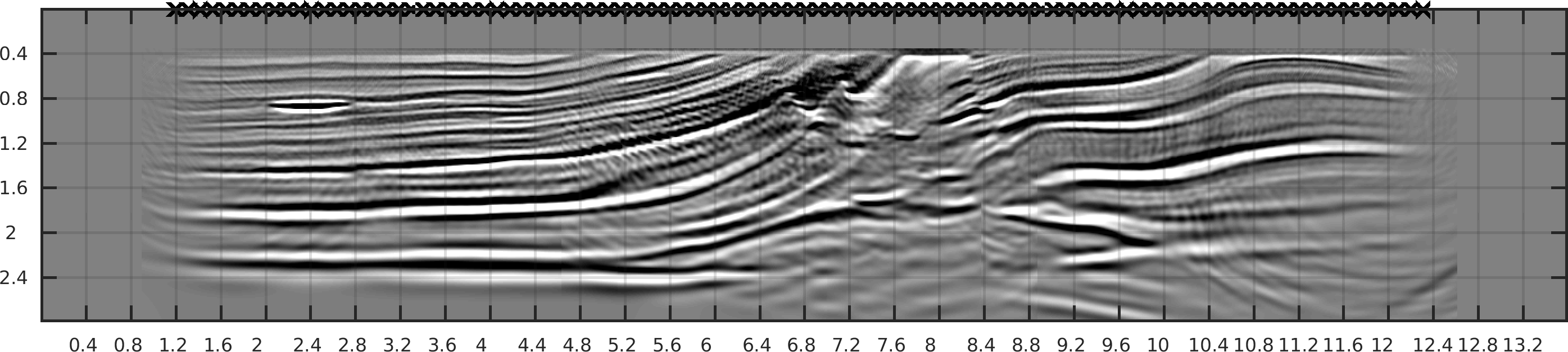}
\includegraphics[width=0.97\textwidth]{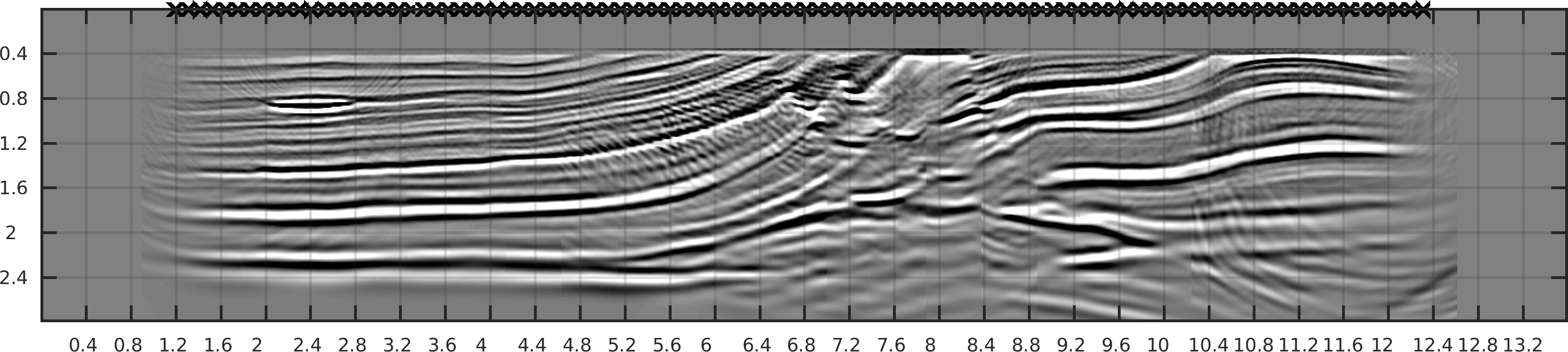}
\caption{Numerical experiment setting and results for Marmousi model. Top to bottom: 
true velocity $c(x)$; kinematic velocity model $c_o(x)$; composite ROM backprojection 
image $\cI_{BP}^s$ with $s = 5$; composite ROM image $\cI_{BP}^s$ with $s = 11$. 
Locations of transducers are black $\times$. All distances in $km$, velocities in $km/s$.}
\label{fig:marmousi}
\end{figure}


\begin{figure}
\centering
\begin{tabular}{ll}
$\quad \quad \quad \quad \quad \quad \quad \quad c(x)$ & 
$\quad \quad \quad \quad \quad \quad \quad \quad c_o(x)$ \\
\includegraphics[width=0.48\textwidth]{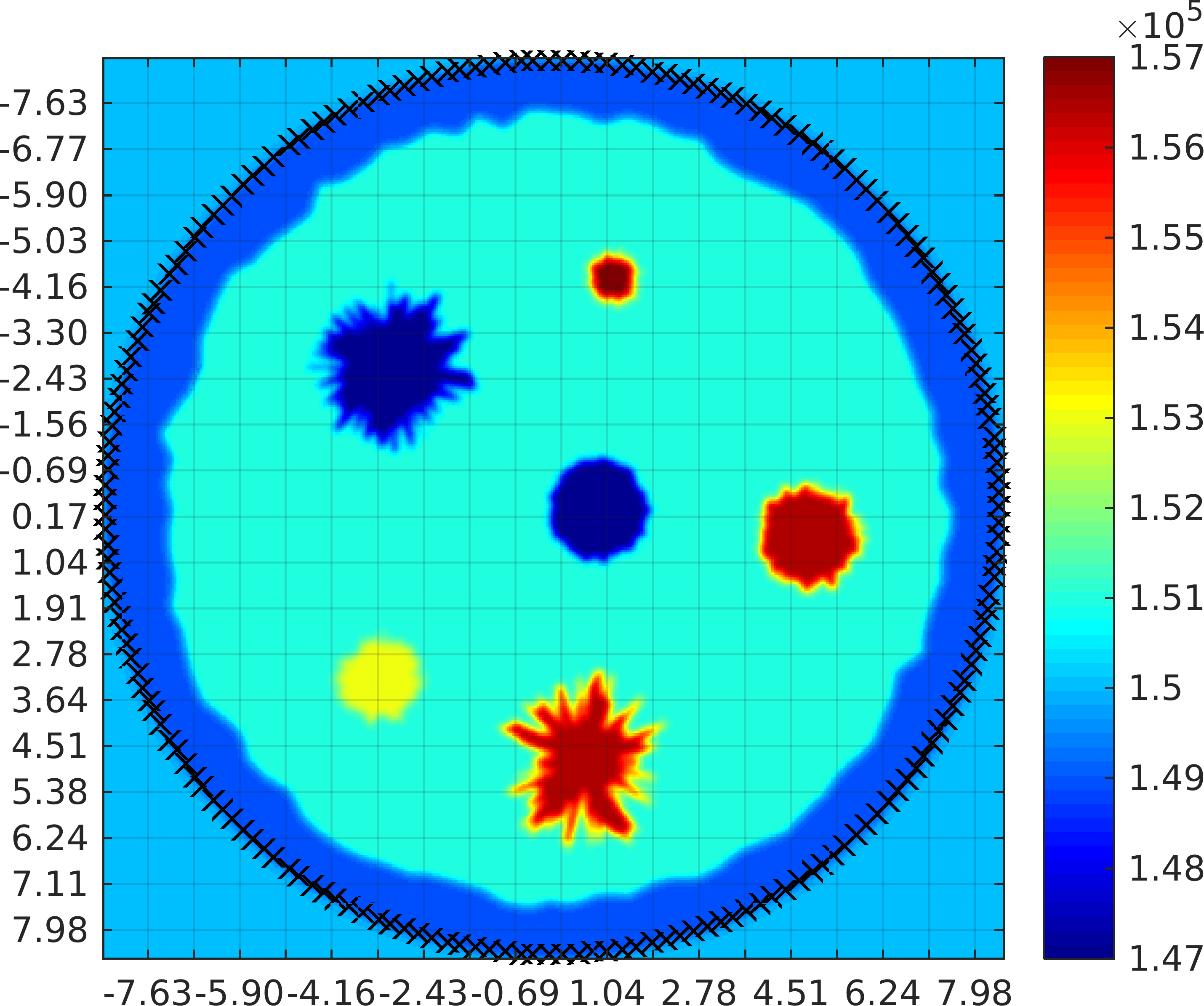} &
\includegraphics[width=0.48\textwidth]{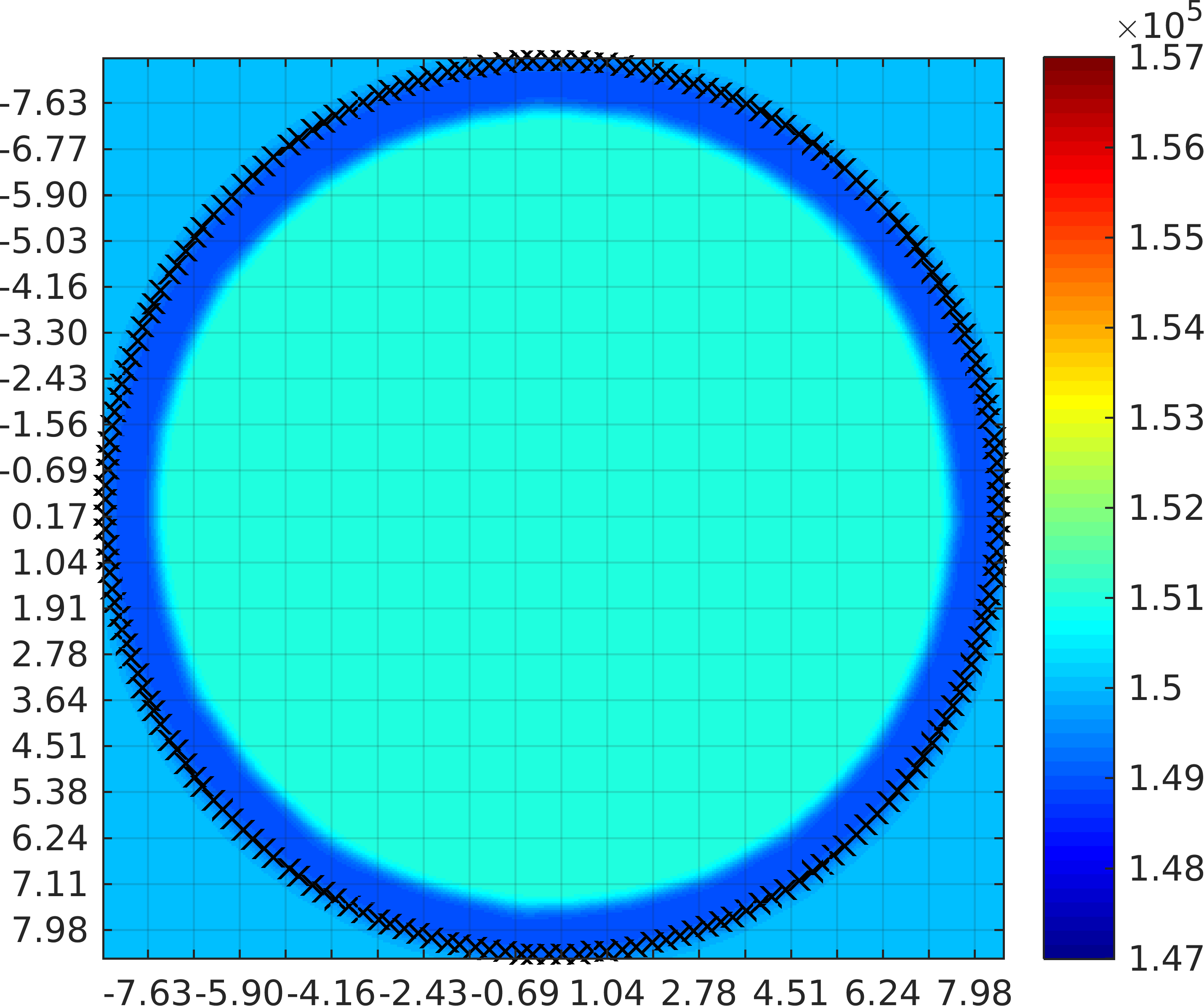} \\
$\quad \quad \quad \quad \quad \quad \quad \quad \quad \cI_{BP}^s$ & 
$\quad \quad \quad \quad \quad \quad \quad \quad \quad \cI_{RTM}$ \\
\includegraphics[width=0.45\textwidth]{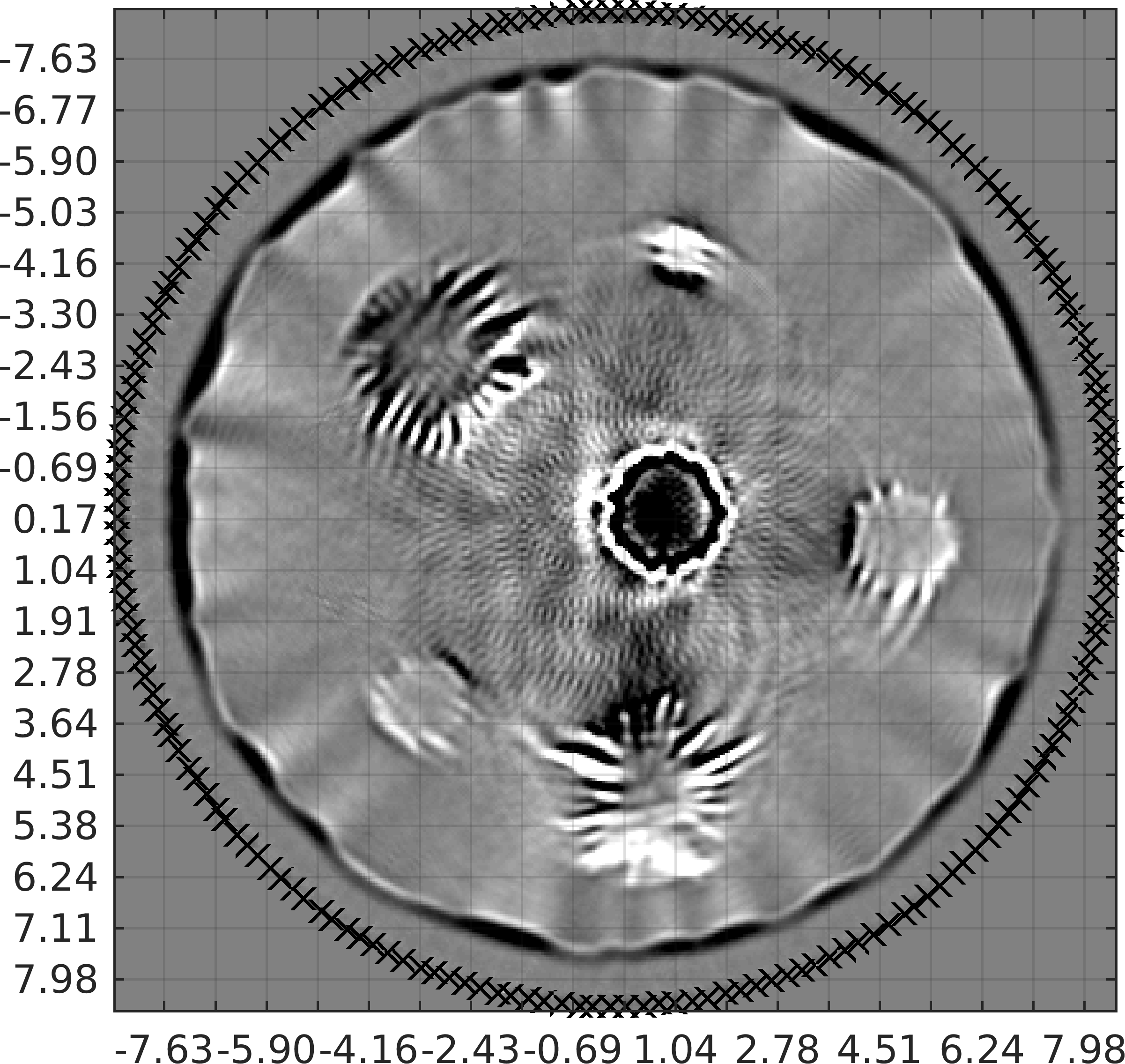} &
\includegraphics[width=0.45\textwidth]{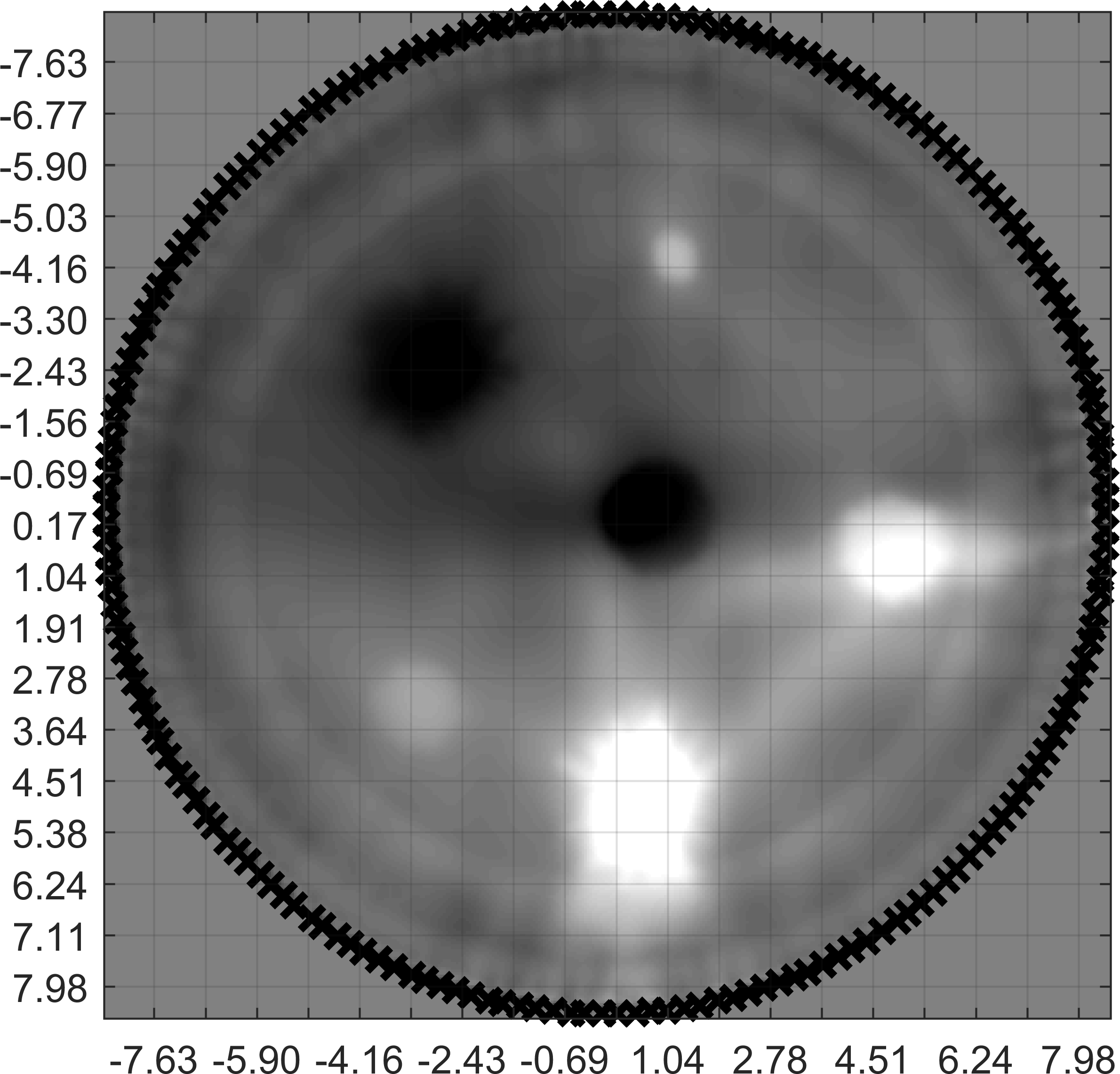}
\end{tabular}
\caption{Numerical experiment setting and results for ultrasound imaging of a breast 
phantom $c(x)$. Locations of transducers are black $\times$. All distances in $cm$, 
velocities in $cm/s$.}
\label{fig:ust}
\end{figure}

We conclude the numerical study with two large-scale examples inspired by the applications of
imaging with acoustic waves in geophysics and medical ultrasound tomography.

The first numerical example is the famous Marmousi model \cite{brougois1990marmousi} that has
become a standard synthetic test model in geophysics used to evaluate the performance of inversion
and imaging algorithms.  The model is discretized on a $10 \; m$ uniform grid with 
$N = 1360 \times 280 =$ $380,800$ nodes. The temporal sampling interval 
$\tau = 1.8 \cdot10^{-2} \; s$ with the corresponding sampling frequency $\omega_\tau = 55.55 \; Hz$. 
The data is measured for $2n = 130$ time instants at the array with $m = 102$ transducers spaced 
uniformly every $110 \; m$ over the accessible boundary $\cB_A = \{ (x, 0) \; | \; 1.2 \; km \leq x \leq 12.3 \; km \}$.
The kinematic model $c_o(x)$ is obtained by convolving the true velocity $c(x)$ with a Gaussian kernel 
of width $400 \; m$ and height of $280 \; m$, see Figure \ref{fig:marmousi}.

A large scale example like Marmousi allows us to illustrate a way to make ROM backprojection imaging
more computationally tractable for large problems. Instead of computing the ROM from the whole 
data set $\bD^k$, $k=0,\ldots,2n-1,$ at once, we split the transducer array $\{1,\ldots,m\}$ into $s$
(possibly overlapping) sub-arrays $S_i = \{m_i^{min},\ldots,m_i^{max} \}$, $i=1,\ldots,s$. 
For each sub-array $i$ we compute the imaging functional $\cI_{BP,i}$ from the restriction of the data 
to that sub-array, i.e. from $\bD^k_{S_i,S_i}$, $k=0,\ldots,2n-1$. Then the \emph{composite} 
backprojection image is simply a weighted sum of individual images
\be
\cI_{BP}^s = \sum_{i=1}^{s} \eta_i \cI_{BP,i},
\ee
with some positive weights $\eta_i$, $i=1,\ldots,s$. Note that such approach essentially discards the
large offset data, since $\bD^k_{S_i,S_i}$ is a block of the full data matrix $\bD^k$ containing the 
entries at most $m_i^{max} - m_i^{min}$ positions away from the main diagonal. Nevertheless,
if the sub-arrays are sufficiently large, the quality of such a composite image is still high, while the
computational cost is substantially lower. 

In Figure \ref{fig:marmousi} we display two composite images for the Marmousi model. The first 
is computed for $s=5$ overlapping sub-arrays of $m_i^{max} - m_i^{min} + 1 = 34$ transducers 
each, while the second one is for $s=11$ overlapping sub-arrays of $m_i^{max} - m_i^{min} + 1 = 17$ 
transducers each. We observe no visible deterioration of the image when using smaller sub-arrays 
(which discards more far offset data), in fact one may argue that the image for $s=11$ is somewhat
sharper than that for $s=5$.

The second large-scale numerical example is a model for medical ultrasound diagnostics of breast
cancer. The model consists of a circular phantom of the acoustic velocity distribution $c(x)$ in a 
cross-section of a human breast, see Figure \ref{fig:ust}. The phantom was designed and provided 
to the authors by the Computational Bioimaging Laboratory at the Washington University in St. Louis, 
Department of Biomedical Engineering, in collaboration with TomoWave Laboratories, Inc. (Houston, TX). 

The overall diameter of the phantom is $17 \; cm$. A homogeneous background corresponds to healthy 
breast tissue. Adjacent to the reflective boundary (skin surface) is a fatty layer with rough inner surface. 
Embedded into the background are six inclusions that represent both benign and malignant lesions. 
While four round-shaped inclusions represent benign lesions, two star-shaped inclusions 
(top left and bottom-most, see Figure \ref{fig:ust}) correspond to malignant tumors. The shape, 
dimensions and acoustic velocity of inclusions are chosen according to human physiology.

The model is discretized on a $4.34 \; mm$ uniform grid with $N = 392 \times 392 = 153,664$ nodes 
and a $784 \times 784$ fine grid ($N = 614,656$). The temporal sampling interval 
$\tau = 8.6 \cdot 10^{-7} \; s$ with the corresponding sampling frequency $\omega_\tau = 1.16 \; MHz$. 
The data is measured for $2n = 80$ time instants at the array with $m = 192$ transducers spaced 
uniformly over the whole boundary. A composite image is formed with the transducer array split into 
$s=16$ overlapping sub-arrays of $m_i^{max} - m_i^{min} + 1 = 24$ transducers each.

We observe in Figure \ref{fig:ust} that the backprojection image gives a reasonably good 
characterization of the shape of inclusions. The star-shaped (malignant) lesions are imaged to have 
much rougher boundaries than their round-shaped (benign) counterparts. This feature of makes ROM 
backprojection imaging promising for medical ultrasound diagnostics. By comparison, the shapes of 
the inclusions in the RTM image are substantially more blurred, which makes the distinction between 
star-shaped and round-shaped lesions difficult. Partially, this is a consequence of lower resolution 
of RTM compared to ROM backprojection for the same source frequency. Note that the medical ultrasound 
imaging problem (for soft tissues) exhibits substantially less nonlinearity than geophysical imaging 
problem. While the contrast in acoustic velocity in geophysical applications can be as high as 
$300 \%$, a typical velocity contrast in human soft tissues is within $5 \%$.

\section{Conclusions and future work}
\label{sec:conclude}

We presented a novel method for imaging the discontinuities of the acoustic velocity $c(x)$, 
a coefficient of the wave equation (\ref{eqn:wavep}), from discretely sampled time domain data
$\bD^k$ measured at an array of transducers. Unlike conventional migration methods, including 
time reversal approaches, the proposed ROM backprojection imaging method is nonlinear in the data.
The nonlinearity comes from an implicit orthogonalization of time-domain wavefield snapshots that
form a basis for the subspace on which the propagator is projected. This allows the method to handle
the nonlinear interactions between reflectors that often lead to artifacts in linear migration images.
In particular, ROM backprojection images are almost entirely free of multiple reflection artifacts 
that require special handling in conventional time reversal approaches. While ROM computation can
become unstable in the presence of noise in the data, we propose a computationally inexpensive 
regularization procedure similar to Tikhonov regularization for linear inverse problems.

Although the numerical study in section \ref{sec:num} demonstrate the potential of ROM backprojection
imaging approach and its advantages over conventional linear migration methods, there is a number of 
questions that must be investigated further for ROM backprojection imaging to become a viable
alternative to migration in real world applications. We identify the following directions for future research:
\begin{itemize}
\item \emph{Resolution analysis.} The numerical results in section \ref{subsec:imgcomprtm} imply that 
ROM backprojection imaging may achieve higher resolution in the range direction compared to 
conventional RTM. This question requires a careful separate consideration. One may consider first the
one-dimensional case from \cite{druskin2015direct}, where it might be possible to explicitly construct
the orthogonalized snapshots $\bV_{(o)}$ and their outer products $\bV_{(o)} \bV_{(o)}^*$
which characterize the resolution.
\item \emph{Advanced regularization techniques.} While the simple regularization procedure in 
Algorithm \ref{alg:datareg} provides reasonably good images from noisy data, more elaborate 
approaches can be studied. If one views the approach of Algorithm \ref{alg:datareg} as an analogue
of Tikhonov regularization for linear inverse problems, then it is natural to also consider an analogue 
of truncated SVD regularization \cite{hansen1987truncatedsvd}. In the context of ROM computation
this means projecting both mass $\bM$ and stiffness $\bS$ matrices on the appropriate subspace
of singular vectors of the mass matrix. The subspace should be chosen so that the projected mass 
matrix is positive definite. Note that since the backprojection imaging formula (\ref{eqn:imgfun})
involves the difference between the projected propagators for the unknown medium and the kinematic 
model, there must be consistency in the choice of singular vector subspaces for the unknown medium
and the kinematic model.
\item \emph{Non-collocated sources and receivers.} We assumed that the data is measured at an
array of transducers that can act both as sources and receivers, thus allowing us to write the data
model $\bD(t)$ in the symmetrized form (\ref{eqn:datafsym}). Consequently, this lead to both left
and right projection subspaces in interpolation relation (\ref{eqn:lemmainterp}) to be the same block 
Krylov subspace (\ref{eqn:krylov}), which simplifies the method. While an assumption of collocated 
sources and receivers is not unreasonable for medical ultrasound imaging, it is often not valid in 
geophysical exploration, where the sources (shots) and receivers (geo-/hydro-phones) are almost 
always distinct, with the number of shots being significantly smaller than the number of receivers.
This may require considering different left and right projection subspaces, which leads to implicit
bi-orthogonalization of wavefield snapshots instead of implicit orthogonalization (\ref{eqn:bVorth})
used currently.
\end{itemize}

\section*{Acknowledgments}

This material is based upon research supported in part by the U.S. Office of Naval Research under 
award number N00014-17-1-2057 to L. Borcea and A.V. Mamonov. The work of A.V. Mamonov was 
also partially supported by the National Science Foundation grant DMS-1619821 and by the 
New Faculty Research Program (2015-2016) of the University of Houston.

\appendix

\section{Proof of Lemma \ref{lemma:proj}}
\label{app:lemmaproj}

As mentioned in section \ref{sec:rom}, the choice of an orthonormal basis $\bV$ for $\cK_n(P, \bb)$ 
does not affect the interpolation condition (\ref{eqn:lemmainterp}). Thus, we can work with any basis 
for $\cK_n(P, \bb)$. In particular, for the purposes of this proof it is convenient to take
\be
\bV = \bU (\bU^* \bU)^{-1/2}  = \bU \bM^{-1/2},
\ee
which is obviously orthonormal in the sense of (\ref{eqn:matprodmn})--(\ref{eqn:bVorth}). In what 
follows we use extensively
\be
\bV \bV^* = \bU \bM^{-1} \bU^* \quad \mbox{and} \quad \buh^{k}  = \bU \bE^k,
\ee
where $\bE^k \in \mathbb{R}^{mn \times m}$ has $\bI_m$ block at the $k^{\mbox{th}}$
block position and zeros elsewhere.

We prove first that 
\be
\buh^k = \bV T_k(\bPt) \bBt, \quad k=0,\ldots,n-1.
\label{eqn:buhk}
\ee
The proof is by induction with two base cases being $k=0,1$.

If we notice that
\be
\bBt = \bV^* \bb = \bM^{-1/2} \bU^* \buh^0 = \bM^{1/2} (\bU^* \bU)^{-1} \bU^* \bU \bE^0 = 
\bM^{1/2} \bE^0,
\label{eqn:}
\ee
then the base case $k=0$ is established by
\be
\buh^0 = \bU \bE^0 = \bV \bM^{1/2} \bE^0 = \bV \bBt = \bV T_0(\bPt) \bBt.
\ee
For $k=1$ we have
\be
\begin{split}
\bV T_1(\bPt) \bBt & = \bV \bV^* P \bV \bV^* \bb = 
\bU \bM^{-1} \bU^* P \bU \bM^{-1} \bU^* \bU \bE^0 \\
& =  \bU \bM^{-1} \bU^* P \buh^0 = 
\bU \bM^{-1} \bU^* \buh^1 =
\bU \bM^{-1} \bU^* \bU \bE^1 = 
\buh^1.
\end{split}
\ee

Now we apply the three term recurrence relation for Chebyshev polynomials
\be
T_{k+1}(x) = 2 x T_k(x) - T_{k-1}(x),
\ee
to prove the general case using  (\ref{eqn:actionP}) and the induction hypothesis
\be
\begin{split}
\bV T_{k+1}(\bPt) \bBt & = 2 \bV \bPt T_k(\bPt) \bBt - \bV T_{k-1} (\bPt) \bBt =
2 \bV \bV^* P \bV T_k(\bPt) \bBt - \buh^{k-1} \\
& = 2 \bV \bV^* P \buh^k - \buh^{k-1} = 
\bU \bM^{-1} \bU^* (\buh^{k+1} + \buh^{k-1}) - \buh^{k-1} \\
& = \bU \bM^{-1} \bU^* \bU (\bE^{k+1} + \bE^{k-1}) - \buh^{k-1} = 
\buh^{k+1} + \buh^{k-1} - \buh^{k-1} = \buh^{k+1}.
\end{split}
\ee
This establishes (\ref{eqn:buhk}) which immediately implies the interpolation condition
\be
\bD^k = \bb^* \buh^k = \bb^* \bV T_{k}(\bPt) \bBt = \bBt^* T_{k}(\bPt) \bBt,
\label{eqn:proofinterp}
\ee
for $k=0,\ldots,n-1$. 

To prove the interpolation condition (\ref{eqn:lemmainterp}) for $k=n,\ldots,2n-1$
we consider the following. From the product formula for Chebyshev polynomials 
(\ref{eqn:chebprod}) it follows that
\be
 T_{n-1+l}(x) = 2 T_{n-1} (x) T_{l}(x) - T_{|n-1-l|}(x),
 \label{eqn:chebnm1pl}
\ee
thus 
\be
\buh^{n-1+l} = 2 T_{n-1}(P) \buh^{l} - \buh^{|n-1-l|},
\ee
and
\be
\bD^{n-1+l} = 2 (\buh^{n-1})^* \buh^{l} - \bD^{|n-1-l|},
\quad l=0,1,\ldots.
\label{eqn:bdnm1pl}
\ee

Now suppose that $1 \leq l \leq n-1$, then by induction proof above, (\ref{eqn:buhk}) holds for 
$\buh^{l}$ and also $|n-1-l| \leq n-1$, so (\ref{eqn:proofinterp}) holds for $\bD^{|n-1-l|}$. 
Substituting both into (\ref{eqn:bdnm1pl}) we find
\be
\begin{split}
\bD^{n-1+l} & = 2 \bBt^* T_{n-1}(\bPt) \bV^* \bV T_{l}(\bPt) \bBt - \bBt^* T_{|n-1-l|}(\bPt) \bBt \\
& = \bBt^*  \left( 2 T_{n-1}(\bPt) T_{l}(\bPt) - T_{|n-1-l|}(\bPt) \right) \bBt \\
& = \bBt^* T_{n-1+l}(\bPt) \bBt,
\end{split}
\ee
by (\ref{eqn:chebnm1pl}), which proves the interpolation condition (\ref{eqn:lemmainterp})
for $k = n-1+l = n,\ldots,2n-2$. So we have established (\ref{eqn:lemmainterp}) for all $k$
except $k = 2n-1$. 

The case $k = 2n-1$ follows from (\ref{eqn:bdnm1pl}) and (\ref{eqn:chebnm1pl}) because
\be
\begin{split}
\bD^{2n-1} & = 2 (\buh^{n-1})^* \buh^{n} - \bD^1 = 
2 (\buh^{n-1})^* ( 2P \buh^{n-1} - \buh^{n-2} ) - \bD^1 \\
& = 4 \bBt^* T_{n-1}(\bPt) \bV^* P \bV T_{n-1}(\bPt) \bBt - 
2 \bBt^* T_{n-1}(\bPt) \bV \bV^* T_{n-2}(\bPt) \bBt 
- \bBt^* T_1(\bPt) \bBt \\
& = \bBt^* \left( 2T_{n-1}(\bPt) \left( 2 \bPt T_{n-1}(\bPt) - T_{n-2}(\bPt) \right) - T_1(\bPt) \right) \bBt \\
& = \bBt^* \left( 2T_{n-1}(\bPt) T_{n}(\bPt) - T_1(\bPt) \right) \bBt =
\bBt^* T_{2n-1}(\bPt) \bBt.
\end{split}
\ee

\section{Proof of Lemma \ref{lemma:blocktri}}
\label{app:lemmablocktri}

Since $\bPt$ is self-adjoint, it is enough to consider 
\be
\bPt_{k+l, k} = (\bv^{k+l})^* (P \bv^{k}),
\quad k=0,\ldots,n-1, \quad 0 \leq l \leq n-1-k
\ee
to show that $\bPt$ is block tridiagonal.

Let us denote $\bR = \bL^{-*}$, then from block QR decomposition (\ref{eqn:blockqr}) we have
\be
\bV = \bU \bL^{-*} = \bU \bR.
\ee

Now consider
\be
\bPt_{k+l, k} = (\bv^{k+l})^* (P \bv^{k}) = (\bv^{k+l})^* \left( P \sum_{i=0}^{k} \buh^i \bR_{i,k} \right).
\label{eqn:pbtkpll1}
\ee
Using the action of the propagator on snapshots (\ref{eqn:actionP}) we can rewrite 
(\ref{eqn:pbtkpll1}) as
\be
\bPt_{k+l, k} = \frac{1}{2} \sum_{i=0}^{k} 
\left( (\bv^{k+l})^* \buh^{i+1} + (\bv^{k+l})^* \buh^{|i-1|}  \right) \bR_{i,k}.
\ee
But $\bV^* \bU = \bL^*$, hence $(\bv^{k+l})^* \buh^{i} = (\bL^*)_{k+l, i}$, 
thus
\be
\bPt_{k+l, k} = \frac{1}{2} \sum_{i=0}^{k} 
\left( (\bL^*)_{k+l, i+1} + (\bL^*)_{k+l, |i-1|} \right) \bR_{i,k}.
\ee

The blocks away from the main block diagonal and block sub-diagonal correspond to $l \geq 2$.
Since $i \leq k$ we have $k+l > i+1$ and $k+l > |i-1|$, but $\bL^*$ is block upper-triangular,
thus $(\bL^*)_{k+l, i+1} = (\bL^*)_{k+l, |i-1|} = \bzero$. So we conclude
\be
\bPt_{k, k+l}^* = \bPt_{k+l, k} = \bzero, \quad k=0,\ldots,n-1, \quad 2 \leq l \leq n-1-k,
\ee
i.e. $\bPt$ is block tridiagonal.

For $\bBt$ we simply have
\be
\bBt_k = (\bV^* \bb)_k = (\bv^{k})^* \buh^0 = (\bL^*)_{k,0} = \bzero, 
\quad \mbox{for } k>0,
\ee
which shows (\ref{eqn:lemmabbt}) and concludes the proof.

\bibliographystyle{siam}
\bibliography{biblio}

\begin{thebibliography}{10}

\bibitem{antoulas2010interpolatory}
{\sc Athanasios~C Antoulas, Christopher~A Beattie, and Serkan Gugercin}, {\em
  Interpolatory model reduction of large-scale dynamical systems}, in Efficient
  Modeling and Control of Large-Scale Systems, Springer, 2010, pp.~3--58.

\bibitem{antoulas2001survey}
{\sc Athanasios~C Antoulas, Danny~C Sorensen, and Serkan Gugercin}, {\em A
  survey of model reduction methods for large-scale systems}, Contemporary
  mathematics, 280 (2001), pp.~193--220.

\bibitem{baysal1983reverse}
{\sc Edip Baysal, Dan~D Kosloff, and John~WC Sherwood}, {\em {Reverse time
  migration}}, Geophysics, 48 (1983), pp.~1514--1524.

\bibitem{borcea2008electrical}
{\sc Liliana Borcea, Vladimir Druskin, and Fernando Guevara~Vasquez}, {\em
  Electrical impedance tomography with resistor networks}, Inverse Problems, 24
  (2008), p.~035013.

\bibitem{borcea2011resistor}
{\sc Liliana Borcea, Vladimir Druskin, Fernando Guevara~Vasquez, and
  Alexander~V Mamonov}, {\em Resistor network approaches to electrical
  impedance tomography}, Inverse Problems and Applications: Inside Out II,
  Math. Sci. Res. Inst. Publ, 60 (2012), pp.~55--118.

\bibitem{borcea2010circular}
{\sc Liliana Borcea, Vladimir Druskin, and Alexander~V Mamonov}, {\em Circular
  resistor networks for electrical impedance tomography with partial boundary
  measurements}, Inverse Problems, 26 (2010), p.~045010.

\bibitem{borcea2010pyramidal}
{\sc Liliana Borcea, Vladimir Druskin, Alexander~V Mamonov, and Fernando
  Guevara~Vasquez}, {\em Pyramidal resistor networks for electrical impedance
  tomography with partial boundary measurements}, Inverse Problems, 26 (2010),
  p.~105009.

\bibitem{borcea2014model}
{\sc Liliana Borcea, Vladimir Druskin, Alexander~V Mamonov, and Mikhail
  Zaslavsky}, {\em A model reduction approach to numerical inversion for a
  parabolic partial differential equation}, Inverse Problems, 30 (2014),
  p.~125011.

\bibitem{borcea2011study}
{\sc Liliana Borcea, Fernando Guevara~Vasquez, and Alexander~V Mamonov}, {\em
  Study of noise effects in electrical impedance tomography with resistor
  networks}, Inverse Problems and Imaging, 7 (2013), pp.~417--443.

\bibitem{borcea2016discrete}
\leavevmode\vrule height 2pt depth -1.6pt width 23pt, {\em {A discrete
  Liouville identity for numerical reconstruction of Schr\"{o}dinger
  potentials}}, arXiv preprint arXiv:1601.07603,  (2016).

\bibitem{borcea2002imaging}
{\sc Liliana Borcea, George Papanicolaou, Chrysoula Tsogka, and James
  Berryman}, {\em {Imaging and time reversal in random media}}, Inverse
  Problems, 18 (2002), p.~1247.

\bibitem{brougois1990marmousi}
{\sc A~Brougois, Marielle Bourget, Patriek Lailly, Michel Poulet, Patrice
  Ricarte, and Roelof Versteeg}, {\em Marmousi, model and data}, in EAEG
  Workshop-Practical Aspects of Seismic Data Inversion, 1990.

\bibitem{bunks1995multiscale}
{\sc Carey Bunks, Fatimetou~M Saleck, S~Zaleski, and G~Chavent}, {\em
  Multiscale seismic waveform inversion}, Geophysics, 60 (1995),
  pp.~1457--1473.

\bibitem{chang1987elastic}
{\sc Wen-Fong Chang and George~A McMechan}, {\em Elastic reverse-time
  migration}, Geophysics, 52 (1987), pp.~1365--1375.

\bibitem{claerbout1971toward}
{\sc Jon~F Claerbout}, {\em Toward a unified theory of reflector mapping},
  Geophysics, 36 (1971), pp.~467--481.

\bibitem{crase1990robust}
{\sc E~Crase, A~Pica, Mark Noble, J~McDonald, and A~Tarantola}, {\em Robust
  elastic nonlinear waveform inversion: Application to real data}, Geophysics,
  55 (1990), pp.~527--538.

\bibitem{derode1995robust}
{\sc Arnaud Derode, Philippe Roux, and Mathias Fink}, {\em Robust acoustic time
  reversal with high-order multiple scattering}, Physical review letters, 75
  (1995), p.~4206.

\bibitem{docherty1991brief}
{\sc Paul Docherty}, {\em {A brief comparison of some Kirchhoff integral
  formulas for migration and inversion}}, Geophysics, 56 (1991),
  pp.~1164--1169.

\bibitem{druskin2015direct}
{\sc Vladimir Druskin, Alexander Mamonov, Andrew~E Thaler, and Mikhail
  Zaslavsky}, {\em Direct, nonlinear inversion algorithm for hyperbolic
  problems via projection-based model reduction}, {SIAM Journal on Imaging
  Sciences}, 9 (2016), pp.~684--747.

\bibitem{druskin2013solution}
{\sc Vladimir Druskin, Valeria Simoncini, and Mikhail Zaslavsky}, {\em
  {Solution of the time-domain inverse resistivity problem in the model
  reduction framework Part I. One-dimensional problem with SISO data}}, SIAM
  Journal on Scientific Computing, 35 (2013), pp.~A1621--A1640.

\bibitem{edler2004use}
{\sc Inge Edler and Carl~Hellmuth Hertz}, {\em The use of ultrasonic
  reflectoscope for the continuous recording of the movements of heart walls.},
  Clinical physiology and functional imaging, 24 (2004), pp.~118--136.

\bibitem{grimme1997krylov}
{\sc Eric~James Grimme}, {\em Krylov projection methods for model reduction},
  PhD thesis, Citeseer, 1997.

\bibitem{hagedoorn1954process}
{\sc Johan~Gregorius Hagedoorn}, {\em A process of seismic reflection
  interpretation}, Geophysical Prospecting, 2 (1954), pp.~85--127.

\bibitem{hansen1987truncatedsvd}
{\sc Per~Christian Hansen}, {\em The truncatedsvd as a method for
  regularization}, BIT Numerical Mathematics, 27 (1987), pp.~534--553.

\bibitem{kremkau1986artifacts}
{\sc Frederick~W Kremkau and KJ~Taylor}, {\em Artifacts in ultrasound
  imaging.}, Journal of ultrasound in medicine, 5 (1986), pp.~227--237.

\bibitem{kuchment2013radon}
{\sc Peter Kuchment}, {\em The Radon transform and medical imaging}, SIAM,
  2013.

\bibitem{mamonov2015rombackproj}
{\sc Alexander~V Mamonov, Vladimir Druskin, and Mikhail Zaslavsky}, {\em
  {Nonlinear seismic imaging via reduced order model backprojection}}, SEG
  Technical Program Expanded Abstracts: 2015, pp.~4375--4379.

\bibitem{schneider1978integral}
{\sc William~A Schneider}, {\em Integral formulation for migration in two and
  three dimensions}, Geophysics, 43 (1978), pp.~49--76.

\bibitem{symes1995mathematics}
{\sc William~W Symes}, {\em Mathematics of reflection seismology}, Rice
  University,  (1995), pp.~1--85.

\bibitem{symes2007reverse}
\leavevmode\vrule height 2pt depth -1.6pt width 23pt, {\em Reverse time
  migration with optimal checkpointing}, Geophysics, 72 (2007),
  pp.~SM213--SM221.

\bibitem{verschuur1992adaptive}
{\sc Dirk~J Verschuur, AJ~Berkhout, and CPA Wapenaar}, {\em {Adaptive
  surface-related multiple elimination}}, Geophysics, 57 (1992),
  pp.~1166--1177.

\bibitem{virieux2009overview}
{\sc Jean Virieux and St{\'e}phane Operto}, {\em An overview of full-waveform
  inversion in exploration geophysics}, Geophysics, 74 (2009), pp.~WCC1--WCC26.

\bibitem{weglein1997inverse}
{\sc Arthur~B Weglein, Fernanda~Ara{\'u}jo Gasparotto, Paulo~M Carvalho, and
  Robert~H Stolt}, {\em {An inverse-scattering series method for attenuating
  multiples in seismic reflection data}}, Geophysics, 62 (1997),
  pp.~1975--1989.

\bibitem{yoon2004challenges}
{\sc Kwangjin Yoon, Kurt~J Marfurt, William Starr, et~al.}, {\em Challenges in
  reverse-time migration}, in 2004 SEG Annual Meeting, Society of Exploration
  Geophysicists, 2004.

\bibitem{ziskin1982comet}
{\sc MC~Ziskin, DI~Thickman, NJ~Goldenberg, MS~Lapayowker, and JM~Becker}, {\em
  The comet tail artifact.}, Journal of Ultrasound in Medicine, 1 (1982),
  pp.~1--7.

\end{thebibliography}

\end{document}